\documentclass[twoside,11pt,leqno]{article}
\usepackage{amssymb}
\usepackage{amsthm}
\usepackage{amsfonts}
\usepackage[tbtags]{amsmath}
\usepackage{doc}
\usepackage{latexsym}
\font\headd=cmr8

\font\Bbb=msbm10 \theoremstyle{change}

\begin{document}
\thispagestyle{plain}
 \markboth{}{}
\small{\addtocounter{page}{198} \pagestyle{plain}
\noindent{\scriptsize KYUNGPOOK Math J. 42(2002), 199-272}
\vspace{0.2in}\\
\noindent{\large\bf The Method of Orbits for Real Lie Groups}
\footnote{(Received: December 4, 2001. Revised: May 14, 2002.)\\
  \indent{}Key words and phrases: quantization, the Kirillov correspondence, nilpotent orbits, the Kostant-Sekiguchi
   correspondence, minimal representations, Heisenberg groups, the Jacobi group. \\
   This work was supported by INHA UNIVERSITY Research Grant.(INHA-21382)
   }
\vspace{0.15in}\\
\noindent{\sc Jae-Hyun Yang }
\newline
{\it Department of Mathematics, Inha University, Incheon 402-751,
Korea\\
e-mail} : {\verb|jhyang@math.inha.ac.kr|}
\vspace{0.2in}\\
\noindent{\footnotesize(2000 Mathematics Subject Classification:
Primary 22-XX, 20C35.)}
\vspace{0.15in}\\
{\footnotesize \indent In this paper, we outline a development of
the theory of orbit method for representations of real Lie groups.
In particular, we study the orbit method for representations of
the Heisenberg group and the Jacobi group.}
\vspace{0.2in}\\
\begin{center}
\bf{Contents}
\end{center}
 1. Introduction\\
 2. Quantization\\
 3. The Kirillov Correspondence\\
 4. Auslander-Kostant's Theorem\\
 5. The Obstacle for the Orbit Method\\
  \indent 5.1. Compact Lie Groups\\
  \indent 5.2. Semisimple Lie Groups\\
 6. Nilpotent Orbits and the Kostant-Sekiguchi Correspondence\\
  \indent 6.1. Jordan Decomposition\\
  \indent 6.2. Nilpotent Orbits\\
  \indent 6.3. The Kostant-Sekiguchi Correspondence\\
  \indent 6.4. The Quantization of the $K$-action (due to D. Vogan)\\
 7. Minimal Representations\\
 8. The Heisenberg Group $H_\mathbb {R}^{(g,h)}$\\
  \indent 8.1. Schr{\"o}dinger Representations\\
  \indent 8.2. The Coadjoint Orbits of Picture\\
 9. The Jacobi Group\\
  \indent 9.1. The Jacobi Group $G^J$\\
  \indent\qquad 9.1.1. The Standard Coordinates of the Jacobi Group
  $G^J$\\
  \indent\qquad 9.1.2. The Iwasawa Decomposition of the Jacobi Group $G^J$\\
  \indent 9.2. The Lie Algebra of the Jacobi Group $G^J$\\
  \indent 9.3. Jacobi Forms\\
  \indent 9.4. Characterization of Jacobi Forms as Functions on the Jacobi Group $G^J$\\
  \indent 9.5. Unitary Representations of the Jacobi Group $G^J$\\
  \indent 9.6. Duality Theorem\\
  \indent 9.7. Coadjoint Orbits for the Jacobi Group $G^J$
\vspace{0.2in}\\
\noindent{\bf 1. Introduction} \setcounter{equation}{0}
\renewcommand{\theequation}{1.\arabic{equation}}
\vspace{0.1in}\\
\indent Research into representations of Lie groups was motivated
on the one hand by physics, and on the other hand by the theory of
automorphic forms. The theory of unitary or admissible
representations of noncompact reductive Lie groups has been
developed systematically and intensively shortly after the end of
World War II. In particular, Harish-Chandra, R. Langlands, Gelfand
school and some other people made an enormous contribution to the
theory of unitary representations of noncompact reductive Lie
groups.
\pagestyle{myheadings}
 \markboth{\headd Jae-Hyun Yang $~~~~~~~~~~~~~~~~~~~~~~~~~~~~~~~~~~~~~~~~~~~~~~~~~~~~~~~~~~$}
 {\headd $~~~~~~~~~~~~~~~~~~~~~~~~~~~~~~~~~~~$The Method of Orbits for Real Lie Groups}

 \smallskip

 \indent Early in the 1960s A.A. Kirillov \cite{Ki1} first initiated
the orbit method for a nilpotent real Lie group attaching an
irreducible unitary representation to a coadjoint orbit (which is
a homogeneous symplectic manifold) in a perfect way. Thereafter
Kirillov's work was generalized to solvable groups of type I by L.
Auslander and B. Kostant\,\cite{AK} early in the 1970s in a nice
way. Their proof was based on the existence of complex
polarizations satisfying a positivity condition. Unfortunately
Kirillov's work fails to be generalized in some ways to the case
of compact Lie groups or semisimple Lie groups. Relatively simple
groups like $SL(2,\mathbb{R})$ have irreducible unitary
representations that do not correspond to any symplectic
homogeneous space. Conversely, P. Torasso\,\cite{To1} found that
the double cover of $SL(3,\mathbb{R})$ has a homogeneous
symplectic manifold corresponding to no unitary representations.
The orbit method for reductive Lie groups is a kind of a
philosophy but not a theorem. Many large families of orbits
correspond in comprehensible ways to unitary representations, and
provide a clear geometric picture of these representations. The
coadjoint orbits for a reductive Lie group are classified into
three kinds of orbits, namely, hyperbolic, elliptic and nilpotent
ones. The hyperbolic orbits are related to the unitary
representations obtained by the parabolic induction and on the
other hand, the elliptic ones are related to the unitary
representations obtained by the cohomological induction. However,
we still have no idea of attaching unitary representations to
nilpotent orbits. It is known that there are only finitely many
nilpotent orbits. In a certain case, some nilpotent orbits are
corresponded to the so-called {\it unipotent\ representations}.
For instance, a minimal nilpotent orbit is attached to a minimal
representation. In fact, the notion of unipotent representations
is not still well defined. The investigation of unipotent
representations is now under way. Recently D. Vogan\,\cite{Vo6}
presented a {\it new} method for studying the quantization of
nilpotent orbits in terms of the restriction to a maximal compact
subgroup even though it is not complete and is in a preliminary
stage. J.-S. Huang and J.-S. Li\,\cite{HuL} attached unitary
representations to spherical nilpotent orbits for the real
orthogonal and symplectic groups.
 \smallskip

 \indent In this article, we describe a development of the orbit method
for real Lie groups, and then in particular, we study the orbit
method for the Heisenberg group and the Jacobi group in detail.
This paper is organized as follows. In Section 2, we describe the
notion of geometric quantization relating to the theory of unitary
representations which led to the orbit method. The study of the
geometric quantization was first made intensively by A.A.
Kirillov\,\cite{Ki5}. In Section 3, we outline the beautiful
Kirillov's work on the orbit method for a nilpotent real Lie group
done early in the 1960s. In Section 4, we describe the work for a
solvable Lie group of type I done by L. Auslander and B.
Kostant\,\cite{AK} generalizing Kirillov's work. In Section 5, we
roughly discuss the cases of compact or semisimple Lie groups
where the orbit method does not work nicely. If $G$ is compact or
semisimple, the correspondence between $G$-orbits and irreducible
unitary representations of $G$ breaks down. In Section 6, for a
real reductive Lie group $G$ with Lie algebra $\mathfrak{g}$, we
present some properties of nilpotent orbits for $G$ and describe
the Kostant-Sekiguchi correspondence between $G$-orbits in the
cone of all nilpotent elements in $\mathfrak{g}$ and $K_
\mathbb{C}$-orbits in the cone of nilpotent elements on
$\mathfrak{p}_\mathbb{C}$, where $K_{\mathbb{C}}$ is the
complexification of a fixed maximal compact subgroup $K$ of $G$
and $\mathfrak{g}_\mathbb{C}=\mathfrak{k}_\mathbb{C} \oplus
\mathfrak{p}_\mathbb{C}$ is the Cartan decomposition of the
complexification $\mathfrak{g}_\mathbb{C}$ of $\mathfrak{g}$. We
do not know yet how to quantize a nilpotent orbit in general. But
for a maximal compact subgroup $K$ of $G$, D. Vogan attaches a
space with a representation of $K$ to a nilpotent orbit. We
explain this correspondence in a rough way. Most of the materials
in this section come from the article \cite{Vo6}. In Section 7, we
outline the notion of minimal orbits (that are nilpotent orbits),
and the relation of the minimal representations to the theory of
reductive dual pairs initiated first by R. Howe. We also discuss
the recent works for a construction of minimal representations for
various groups. For more detail, we refer to \cite{Li2}. In
Section 8, we study the orbit method for the Heisenberg group in
some detail. In Section 9, we study the unitary representations of
the Jacobi group and their related topics. The Jacobi group
appears in the theory of Jacobi forms. That means that Jacobi
forms are automorphic forms for the Jacobi group. We study the
coadjoint orbits for the Jacobi group.

\bigskip

\noindent{\bf Notation.}\quad We denote by
$\mathbb{Z},\mathbb{R}$, and $\mathbb{C}$ the ring of integers,
the field of real numbers, and the field of complex numbers
respectively. The symbol $\mathbb{C}_1^{\times}$ denotes the
multiplicative group consisting of all complex numbers $z$ with
$\vert z\vert =1$, and the symbol $Sp(n,\mathbb{R})$ the
symplectic group of degree $n$, $H_n$ the Siegel upper half plane
of degree $n$. The symbol ``:='' means that the expression on the
right hand side is the definition of that on the left. We denote
by $\mathbb{Z}^{+}$ the set of all positive integers, by
$F^{(k,l)}$ the set of all $k\times l$ matrices with entries in a
commutative ring $F$. For any $M\in F^{(k,l)},\ ^t\!M$ denotes the
transpose matrix of $M$. For $A\in F^{(k,k)},\ \sigma(A)$ denotes
the trace of $A$. For $A\in F^{(k,l)}$ and $B\in F^{(k,k)},$ we
set $B[A]={^t\!A}BA$. We denote the identity matrix of degree $k$
by $E_k$. For a positive integer $n,\ \text{Symm}\,(n,K)$ denotes
the vector space consisting of all symmetric $n\times n$ matrices
with entries in a field $K.$
\vspace{0.1in}\\
\noindent{\bf 2. Quantization }
\setcounter{equation}{0}
\renewcommand{\theequation}{2.\arabic{equation}}
\vspace{0.1in}\\
\indent The problem of quantization in mathematical physics is to
attach a quantum mechanical model to a classical physical system.
The notion of {\it geometric quantization} had emerged at the end
of the 1960s relating to the theory of unitary group
representations which led to the orbit method. The goal of
geometric quantization is to construct quantum objects using the
geometry of the corresponding classical objects as a point of
departure. In this paper we are dealing with the group
representations and hence the problem of quantization in
representation theory is to attach a unitary group representation
to a symplectic homogeneous space.
\vspace{0.1in}\\
\indent A classical mechanical system can be modelled by the phase
space which is a symplectic manifold. On the other hand, a quantum
mechanical system is modelled by a Hilbert space. Each state of
the system corresponds to a line in the Hilbert space.
\vspace{0.1in}\\
\noindent{\bf Definition 2.1.}\quad A pair $(M,\omega)$ is called
a {\it symplectic} manifold with a nondegenerate closed
differential 2-form $\omega.$ We say that a pair $(M,c)$ is a {\it
Poisson manifold} if $M$ is a smooth manifold with a bivector
$c=c^{ij}\partial_i\partial_j$ such that the Poisson brackets
\begin{equation}
\{ f_1,f_2 \}=c^{ij}\partial_if_1\partial_j f_2
\end{equation}
define a Lie algebra structure on $C^{\infty}(M)$. We define a
{\it Poisson} $G$-manifold as a pair $(M,\ f^{M}_{( \cdot )})$
where $M$ is a Poisson manifold with an action of $G$ and
 $f^{M}_{( \cdot )}: \mathfrak{g} \rightarrow C^{\infty}(M)\,( X\mapsto f^{M}_{X})$ is a
Lie algebra homomorphism  such that the following relation holds:
\begin{equation}
\text{s}\!-\!\text{grad}(f_X^M)=L_X,\quad X\in \mathfrak{g},\quad
X\in\mathfrak{g}.
\end{equation}
 Here $L_{X}$ is the Lie vector field on $M$ associated with $X\in
\mathfrak{g} $, and s-grad$(f)$ denotes the {\it skew gradient} of
a function $f$, that is, the vector field on $M$ such that
\begin{equation}
\text{s-grad}(f) g = \{f, g\}\quad \text{for all}\quad g\in
C^{\infty }(M).
\end{equation}
 \indent For a given Lie group $G$ the collection of all Poisson $G$-manifolds forms the category
$\mathcal{P}(G)$ where a morphism $\alpha : (M, f^{M}_{( \cdot )}
) \rightarrow (N, f^{N}_{( \cdot )} )$ is a smooth map from $M$ to
$N$ which preserves the Poisson brackets: $\{\alpha ^{*}(\phi ),\
\alpha ^{*}(\psi )\}=\alpha ^{*}(\{\phi , \psi \})$ and makes the
following relation holds:
\begin{equation}
\alpha^{*}(f_X^N)=f_X^M,\quad  X \in \mathfrak{g}.
\end{equation}
Observe that the last condition implies that $\alpha $ commutes
with the $G$-action.
\vspace{0.1in}\\
\indent First we explain the mathematical model of classical
mechanics in the Hamiltonian formalism.\\
\indent Let $(M,\omega)$ be a symplectic manifold of dimension
$2n.$ According to the Darboux theorem, the sympletic form
$\omega$ can always be written in the form
\begin{equation}
\omega = \Sigma_{k=1}^n dp_k \wedge dq_k
\end{equation}
in suitable canonical coordinates $p_1,\cdots,p_n,q_1,\cdots,q_n$.
However, these canonical coordinates are not uniquely determined.
\vspace{0.1in}\\
\indent The symplectic from $\omega$ sets up an isomorphism
between the tangent and cotangent spaces at each point of $M$. The
inverse isomorphism is given by a bivector $c$, which has the form
\begin{equation}
c = \Sigma_{k=1}^n {\partial \over {\partial p_k}}{\partial \over
{\partial q_k}}
\end{equation}
 in the same system of coordinates in which the equality (2.5) holds. In the general system of
coordinates the form $\omega$ and the bivecter $c$ are written in
the form
\begin{equation}
\omega = {\Sigma_{i<j}}{\omega}_{ij} dx_i \wedge dx_j , \quad
c={\Sigma_{i<j}}c^{ij} \partial_i  \partial_j
\end{equation}
 with mutually inverse skew-symmetric matrices $(\omega_{ij} )$ and $(c^{ij} )$. The set $C^{\infty} (M)$ of all
smooth functions on $M$ forms a commutative associative algebra
with respect to the usual multiplication. The Poisson bracket $\{\
,\ \}$ defined by
\begin{equation}
\{F,G\} := {\Sigma_{i,j}} {c^{ij}} {{\partial}_j} F \cdot
{\partial}_iG,\quad F,G \in {C^{\infty}}(M)
\end{equation}
defines a Lie algebra structure on ${C^{\infty}} (M)$. The Jacobi
identity for the Poisson bracket is equivalent to the condition
$d{\omega} =0$ and also to the vanishing of the Schouten bracket
\begin{equation}
{[c,c]^{ijk}} := {{\curvearrowleft}_{ijk}}{\Sigma_{m}}
{c^{im}}{{\partial}_m}{c^{jk}}
\end{equation}
where the sign ${\curvearrowleft}_{ijk}$ denotes the sum over the
cyclic permutations of the indices $i, j, k$.\\
\indent {\it Physical quantities} or {\it observables} are
identified with the smooth functions on $M$. A state of the system
is a linear functional on ${C^{\infty}} (M)$ which takes
non-negative values on non-negative functions and equals $1$ on
the function which is identically equal to $1$. The general form
of such a functional is a probability measure $\mu$ on $M$. By a
{\it pure\ state} is meant an extremal point of the set of
states.\\
\indent The dynamics of a system is determined by the choice of a
Hamiltonian function or energy, whose role can be played by an
arbitrary function $H \in {C^{\infty}} (M)$. The dynamics of the
system is described as follows. The states do {\it {not\ }} depend
on time, and the physical quantities are functions of the point of
the phase space and of time. If $F$ is any function on $M {\times}
\mathbb{R}$, that is, any observable, the equations of the motion
have the form
\begin{equation}
\dot{F} := {{\partial F} \over {\partial t}} = \{H,F\}.
\end{equation}
Here the dot denotes the derivative with respect to
time. In particular, applying (2.10) to the canonical variables
$p_k$,\,$q_k$, we obtain {\it Hamilton's\ equations}
\begin{equation}
\dot{q_k}= {{\partial H} \over {\partial p_k}} ,\quad \dot{p_k}=
-{{\partial H} \over {\partial q_k}}.
\end{equation}
 A set $\{ F_1, \cdots, F_m \}$ of physical quantities is called {\it complete} if
the conditions $\{ F_i, G \} =0\, (1 \leq i \leq m)$ imply that
$G$ is a constant.
\vspace{0.1in}\\
\noindent{\bf Definition 2.2.}\quad Let $(M, \omega )$ be a
symplectic manifold and $f \in {C^{\infty}} (M)$ a smooth function
on $M$. The Hamiltonian vector field $\xi_f$ of $f$ is defined by
\begin{equation}
{\xi}_f (g) = \{ f,g \}, \quad g \in C^{\infty} (M)
\end{equation}
where $\{ \ ,\ \}$ is the Poisson bracket on $C^{\infty} (M)$
defined by (2.8). Suppose $G$ is a Lie group with a smooth action
of $G$ on $M$ by symplectomorphisms. We say that $M$ is a {\it
Hamiltonian} $G$-{\it space} if there exist a linear map
\begin{equation}
\tilde{\mu} : \mathfrak{g} \longrightarrow C^{\infty} (M)
\end{equation}
with the following properties (H1)-(H3) :\\
\noindent (H1) $\tilde{\mu}$ intertwines the adjoint action of $G$
on $\mathfrak{g}$ with its action on $C^{\infty}(M)$ ;\\
\noindent (H2) For each $Y \in \mathfrak{g}$, the vector field by
which $Y$ acts on $M$ is $\xi_{\tilde{\mu}(Y)}$;\\
\noindent (H3) $\tilde{\mu}$ is a Lie algebra homomorphism.\\
\indent The above definition can be formulated in the category of
Poisson manifolds, or even of possibly singular Poisson algebraic
varieties. The definition is due to  A. A. Kirillov \cite{Ki2} and B.
Kostant \cite{Ko1}.
\vspace{0.1in}\\
\indent The natural quantum analogue of a Hamiltonian $G$-space is
simply a unitary representation of $G$.
\vspace{0.1in}\\
\noindent{\bf Definition 2.3.}\quad Suppose $G$ is a Lie group. A
{\it unitary\ representation} of $G$ is a pair $(\pi , {\mathcal
H})$ with a Hilbert space $\mathcal{H}$, and $$ \pi : G
\longrightarrow U(\mathcal{H})$$
 a homomorphism from $G$ to the group of unitary operations on $\mathcal{H}$.\\

\smallskip

\indent We would like to have a notion of quantization passing
from Definition 2.2 to Definition 2.3\,: that is, from Hamiltonian
$G$-space to unitary representations.

\smallskip

\indent Next we explain the mathematical model of quantum
mechanics. In quantum mechanics the physical quantities or
observables are self-adjoint linear operators on some complex
Hilbert space $\mathcal{H}$. They form a linear space on which two
bilinear operations are defined\,:
\begin{equation}
A\circ B := {1 \over 2}\,(AB+BA) \quad(\text{Jordan
multiplication})
\end{equation}

\begin{equation}
[A,B]_{\hbar}:= {{2\pi i}\over {\hbar}}\,(AB-BA) \quad(\text{the\
commutator})
\end{equation}
 where $\hbar$ is the Planck's constant.\\

\smallskip

 \indent With respect to (2.14) the set of
observables forms a commutative but not associative algebra. With
respect to (2.15) it forms a Lie algebra. These two operations
(2.14) and (2.15) are the quantum analogues of the usual
multiplication and the Poisson bracket in classical mechanics. The
{\it phase\ space} in quantum mechanics consists of the
non-negative definite operators $A$ with the property that
$\text{tr}\,A=1$. The {\it pure\ states} are the one-dimensional
projection operators on $\mathcal{H}$. The {\it dynamics} of the
system is defined by the {\it energy\ operator} ${\hat{H}}$. If
the states do not depend on time but the quantities change, then
we obtain the {\it Heisenberg\ picture.} The equations of motion
are given by
\begin{equation}
\dot{\hat{A}} = [\hat{H}, \hat{A}]_{\hbar},\quad
\text{(Heisenberg's\ equation)}.
\end{equation}
 The integrals of the system are all the operators which commute with
$\hat{H}$. In particular, the energy operator itself does not
change with time.\\
\indent The other description of the system is the so-called
Schr{\"o}dinger picture. In this case, the operators corresponding
to physical quantities do not change, but the states change. A
pure state varies according to the law
\begin{equation}
\dot{\Psi}= {{2\pi i} \over {\hbar}} \hat{H} \Psi \quad
\text{(Schr{\"o}dinger\ equation)}.
\end{equation}
The eigenfunctions of the Schr{\"o}dinger operator give the
stationary states of the system. We call a set of quantum physical
quantities $\hat{A_1}, \cdots ,\hat{A_m}$ {\it complete} if any
operator $\hat{B}$ which commutes with $\hat{A_i}\, (1 \leq i \leq
m)$ is a multiple of the identity. One can show that this
condition is equivalent to the irreducibility of the set
$\hat{A_1}, \cdots , \hat{A_m}.$

\smallskip

\indent Finally we describe the quantization problem relating to
the orbit method in group representation theory. As I said
earlier, the problem of geometric quantization is to construct a
Hilbert space $\mathcal{H}$ and a set of operators on
$\mathcal{H}$ which give the quantum analogue of this system from
the geometry of a symplectic manifold which gives the model of
classical mechanical system. If the initial classical system had a
symmetry group $G$, it is natural to require that the
corresponding quantum model should also have this symmetry. That
means that on the Hilbert space $\mathcal{H}$ there
should be a unitary representation of the group $G$.\\

\smallskip

\indent We are interested in homogeneous symplectic manifolds on
which a Lie group $G$ acts transitively. If the thesis is true
that every quantum system with a symmetry group $G$ can be
obtained by quantization of a classical system with the same
symmetry group, then the irreducible representations of the group
must be connected with homogeneous symplectic $G$-manifolds. The
orbit method in representation theory intiated first by A. A.
Kirillov early in the 1960s relates the unitary representations of
a Lie group $G$ to the coadjoint orbits of $G$. Later L.
Auslauder, B. Kostaut, M. Duflo, D. Vogan, P. Torasso etc
developed the theory of the orbit method to more general cases.\\
\indent For more details, we refer to \cite{Ki1}- \cite{Ki5},
\cite{Ki9}, \cite{Ko1}-\cite{Ko2} and \cite{Vo3}-\cite{Vo4},
\cite{Vo6}.
\vspace{0.2in}\\
\noindent{\bf 3. The Kirillov Correspondence}
\setcounter{equation}{0}
\renewcommand{\theequation}{3.\arabic{equation}}
\vspace{0.1in}\\
\indent In this section, we review the results of Kirillov on
unitary representations of a nilponent real Lie group. We refer to
\cite{Ki1}-\cite{Ki5}, \cite{Ki9} for more detail.\\
\indent Let $G$ be a  simply connected real Lie group with its Lie
algebra $\mathfrak{g}$. Let ${Ad}_G : G \longrightarrow
GL(\mathfrak{g})$ be the adjoint representation of $G$. That is,
for each $g \in G$, ${{Ad}_G}(g)$ is the differential map of $I_g$
at the identity $e$, where $I_g : G \longrightarrow G $ is the
conjugation by $g$ given by
$$I_g (x):= gxg^{-1} \quad\text{for }\ x \in G.$$
Let $\mathfrak{g}^*$ be the dual space of the vector space
$\mathfrak{g}$. Let $Ad^*_G : G \longrightarrow
GL({\mathfrak{g}}^*)$ be the contragredient of the adjoint
representation $Ad_G$. $Ad^*_G$ is called the coadjoint
representation of $G$. For each $\ell \in {\mathfrak{g}}^*$, we
define the alternating bilinear form $B_{\ell}$ on $\mathfrak{g}$
by
\begin{equation}
B_{\ell} (X,Y) =<[X,Y],\ell>, \quad X,Y \in \mathfrak{g}.
\end{equation}
\vspace{0.1in}\\
\noindent {\bf Definition 3.1.}\quad (1) A Lie subalgebra
$\mathfrak{h}$ of $\mathfrak{g}$ is said to be {\it subordinate}
to $\ell \in {\mathfrak{g}}^*$ if $\mathfrak{h}$ forms a totally
isotropic vector space of $\mathfrak{g}$ relative to the
alternating bilinear form $B_{\ell}$ define by (3.1), i.e.,
$B_{\ell} |_{\mathfrak{h} \times \mathfrak{h}} =0.$\\

\smallskip

\noindent(2) A Lie subalgebra $\mathfrak{h}$ of $\mathfrak{g}$
subordinate to $\ell \in {\mathfrak{g}}^*$ is called a {\it
polarization} of $\mathfrak{g}$ for $\ell$ if $\mathfrak{h}$ is
maximal among the totally isotropic vector subspaces of
$\mathfrak{g}$ relative to $B_{\ell}$. In other words, if $P$ is a
vector subspace of $\mathfrak{g}$ such that $\mathfrak{h}\subset
P$ and $B_{\ell}|_{P\times P}=0,$ then we have $\mathfrak{h} =
P$.\\

\smallskip

\noindent(3) Let $\ell \in \mathfrak{g}^*$ and let $\mathfrak{h}$
be a polarization of $\mathfrak{g}$ for $\ell$. We let $H$ the
simply connected closed subgroup of $G$ corresponding to the Lie
subalgebra $\mathfrak{h}$. We define the unitary character
$\chi_{\ell,\mathfrak{h}}$ of $H$ by
\begin{equation}
\chi_{\ell,\mathfrak{h}} (\text{exp}_H (X)) = e^{2\pi i<X,\ell>}
,\quad X \in \mathfrak{h},
\end{equation}
where $\text{exp}_H : \mathfrak {h} \longrightarrow H$ denotes the
exponential mapping of $\mathfrak {h}$ to $H$. It is known that
$\text{exp}_H$ is surjective.\\

\smallskip

Using the Mackey machinary, Dixmier and Kirillov proved the
following important theorem.
\vspace{0.1in}\\
\noindent{\bf Theorem 3.1
(Dixmier-Kirillov,\,\cite{D1},\cite{Ki1}).}\quad {\it A simply
connected real nilpotent Lie group $G$ is monomial, that is, each
irreducible unitary representation of $G$ can be unitarily induced
by a unitary character of some closed subgroup of $G$}.
\vspace{0.1in}\\
\noindent {\bf Remark 3.2.}\quad More generally, it can be proved
that a simply connected real Lie group whose exponential mapping
is a
diffeomorphism is monomial.\\
Now we may state Theorem 3.1 explicitly.
\vspace{0.1in}\\
\noindent {\bf Theorem\ 3.3 (Kirillov, \cite{Ki1}).}\quad {\it Let
$G$ be a simply connected nilpotent real Lie group with its Lie
algebra $\mathfrak {g} $. Assume that there is given an
irreducible unitary representation $\pi$ of $G$. Then there exist
an element $\ell \in \mathfrak {g} ^*$ and a polarization
$\mathfrak {h}$ of $\mathfrak {g} $ for $\ell$ such that $\pi
\cong \text{Ind}_H^G\chi_{\ell,\mathfrak {h}}$, where
$\chi_{\ell,\mathfrak {h}}$ is the unitary character of $H$
defined by (3.2).}
\vspace{0.1in}\\
\noindent {\bf Theorem\ 3.4 (Kirillov, \cite{Ki1}).}\quad {\it Let
$G$ be a simply connected nilpotent real Lie group with its Lie
algebra $\mathfrak {g} $. If $\ell \in \mathfrak {g} ^*$, there
exists a polarization $\mathfrak {h}$ of $\mathfrak {g} $ for
$\ell$ such that the monomial representation
$\text{Ind}^G_H\chi_{\ell,\mathfrak {h}}$ is irreducible and of
trace class. If $\ell'$ is  an element of $\mathfrak {g} ^*$ which
belongs to the coadjoint orbit $Ad^*_G (G)\ell$ and $\mathfrak
{h}'$ a polarization of $\mathfrak {g} $ for $\ell'$, then the
monomial representations $\text{Ind}_H^G\chi_{\ell,\mathfrak {h}}$
and $\text{Ind}_{H'}^G{\chi}_{\ell',\mathfrak {h}'}$ are unitarily
equivalent. Here $H$ and $H'$ are the simply connected closed
subgroups corresponding to the Lie subalgebras $\mathfrak {h}$ and
$\mathfrak {h}'$ respectively. Conversely, if $\mathfrak {h}$ and
$\mathfrak {h}'$ are polarizations of $\mathfrak {g} $ for $\ell
\in {\mathfrak {g} }^*$ and $\ell' \in {\mathfrak {g} }^*$
respectively such that the monomial representations
$\text{Ind}^G_H\chi_{\ell,\mathfrak {h}}$ and
$\text{Ind}^G_{H'}{\chi}_{\ell',\mathfrak {h}'}$ of $G$  are
unitarily equivalent, then $\ell$ and $\ell'$ belong to the same
coadjoint orbit of $G$ in $\mathfrak {g} ^*$. Finally, for each
irreducible unitary representation $\tau$ of $G$, there exists a
unique coadjoint orbit $\Omega$ of $G$ in $\mathfrak {g} ^*$ such
that for any linear from $\ell \in \Omega$ and each polarization
$\mathfrak {h}$ of $\mathfrak {g} $ for $\ell$, the
representations $\tau$ and $\text{Ind}_H^G\chi_{\ell,\mathfrak
{h}}$ are unitarily equivalent. Any irreducible unitary
representation of $G$ is strongly trace class.}
\vspace{0.1in}\\
\noindent {\bf Remark 3.5.}\quad (a) The bijection of the space
$\mathfrak {g} ^* / G$ of coadjoint orbit of $G$ in $\mathfrak {g}
^*$ onto the unitary dual $\hat G$ of $G$ given by Theorem 3.4 is
called the {\it Kirillov\ correspondence} of $G$. It provides a
parametrization of $\hat G$ by means of the coadjoint orbit space.

\smallskip

\noindent (b) The above Kirillov's work was generalized
immediately to the class known as exponential solvable groups,
which are characterized as those solvable group $G$ whose
simply-connected cover ${\tilde G}$ is such that the exponential
map $\text{exp}:{\tilde{\mathfrak {g} }}\longrightarrow {\tilde
G}$ is a diffeomorphism. For exponential solvable groups, the
bijection between coadjoint orbits and representations holds, and
can be realized using induced representations by an explicit
construction using a polarization just as in the case of a
nilpotent real Lie group. However, two difficulties arise\,:
Firstly not all polarizations yield the same representation, or
even an irreducible representation, and secondly not all
representations are strongly trace class.
\vspace{0.1in}\\
\noindent {\bf Theorem\ 3.6 (I.D.\ Brown).}\quad {\it Let $G$ be a
connected simply connected nilpotent Lie group with its Lie
algebra $\mathfrak {g} $. The Kirillov correspondence
$${\hat G}\longrightarrow {\mathcal {O} }(G)=\mathfrak {g} ^*/G$$ is a
homeomorphism.}
\vspace{0.1in}\\
\noindent {\bf Theorem\ 3.7.} {\it Let $G$ be a connected simply
connected nilpotent real Lie group with Lie algebra $\mathfrak {g}
.$ Let $\mathfrak {h}$ be a subalgebra of $\mathfrak {g} $. Let
$p:\mathfrak {g} ^*\longrightarrow \mathfrak {h}^*$ be the natural
projection. Let $H$ be the simply connected subgroup of $G$ with
its Lie algebra $\mathfrak {h}$. The following (a),(b) and (c)
hold.\\
\indent (a) Let $\pi$ be an irreducible unitary representation of
$G$ corresponding to a coadjoint orbit $\Omega \subset \mathfrak
{g} ^*$ of $G$ via the Kirillov correspondence. Then
$\text{Res}^G_H \pi$ decomposes into the direct integral of
irreducible representations of $H$ corresponding to a coadjoint
orbit $\omega(\subset \mathfrak {h}^*)$ of $H$ such that
$\omega\subset p(\Omega).$

\smallskip

\indent (b) Let $\tau$ be an irreducible unitary representation of
$H$ corresponding to a coadjoint orbit $\omega \subset \mathfrak
{h}^*$ of $H$. Then the induced representation $\text{Ind}_H^G
\tau$ decomposes into the direct integral of irreducible
representations $\pi_{\Omega}$ of $G$ corresponding to coadjoint
orbits $\Omega\subset \mathfrak {g} ^*$ such that
$p(\Omega)\supset \omega.$

\smallskip

\indent (c) Let $\pi_1$ and $\pi_2$ be the irreducible unitary
representations of $G$ corresponding to coadjoint orbits
$\Omega_1$ and $\Omega_2$ respectively. Then the tensor product
$\pi_1\otimes \pi_2$ decomposes into the direct integral of
irreducible representations of $G$ corresponding to coadjoint
orbits $\Omega\subset \mathfrak {g} ^*$ such that $\Omega\subset
\Omega_1 + \Omega_2.$}
\vspace{0.1in}\\
\noindent {\bf Theorem\ 3.8.}\quad {\it Let $G$ be a connected
simply connected nilpotent real Lie group with its Lie algebra
$\mathfrak {g} $. Let $\pi$ be an irreducible unitary
representation of $G$ corresponding to an orbit $\Omega\subset
\mathfrak {g} ^*$. Then the character $\chi_{\pi}$ is a
distribution on ${\mathcal S}(G)$ and its Fourier transform
coincides with the canonical measure on $\Omega$ given by the
symplectic structure. Here ${\mathcal {S} }(G)$ denotes the
Schwarz space of rapidly decreasing functions on $G$.}

\smallskip

\indent A.A. Kirillov gave an explicit formula for the Plancherel
measure on ${\hat G}$. We observe that for a nilpotent Lie group
$G$, we may choose a subspace $Q$ of $\mathfrak {g} ^*$ such that
generic coadjoint orbits intersect $Q$ exactly in one point. We
choose a basis $x_1,\cdots,x_l,y_1,\cdots,y_{n-l}$ in $\mathfrak
{g} $ so that $y_1,\cdots,y_{n-l}$, considered as linear
functionals on $\mathfrak {g} ^*$, are constant on $Q$. Then
$x_1,\cdots,x_l$ are coordinates on $Q$ and hence on an open dense
subset of ${\mathcal {O} }(G)$. For every $f\in Q$ with
coordinates $x_1,\cdots,x_l$, we consider the skew-symmetric
matrix $A=(a_{ij})$ with entries
$$a_{ij}=<[y_i,y_j],f>,\quad 1\leq i,j\leq n-l.$$ We denote by
$p(x_1,\cdots,x_l)$ the Pfaffian of the matrix $A$. We note that
$p(x_1,\cdots,x_l)$ is a homogeneous polynomial of degree
${{n-l}\over 2}$.
\vspace{0.1in}\\
\noindent {\bf Theorem\ 3.9.}\quad {\it Let $G$ be a connected
simply connected nilpotent Lie group. Then the Plancherel measure
on ${\hat G}\cong {\mathcal {O} }(G)$ is concentrated on the set
of generic orbits and it has the form
$$\theta = p(x_1,\cdots,x_l)dx_1\wedge\cdots\wedge dx_l$$ in the
coordinates $x_1,\cdots,x_l.$}
\vspace{0.2in}\\
\noindent{\bf 4. Auslander-Kostant's Theorem }
\setcounter{equation}{0}
\renewcommand{\theequation}{4.\arabic{equation}}
\vspace{0.1in}\\
\indent In this section, we present the results obtained by L.
Auslander and B. Kostant in \cite{AK} together with some complements
suggested by I.M. Shchepochkina \cite{Sh1}-\cite{Sh2}. Early in the 1970s L.
Auslander and B. Kostant described the unitary dual of all
solvable Lie groups of type I.
\vspace{0.1in}\\
\noindent {\bf Theorem\ 4.1.}\quad {\it A connected, simply
connected solvable Lie group $G$ belongs to type I if and only if
the orbit space ${\mathcal {O} }(G)=\mathfrak {g} ^*/G$ is a
$T_0$-space and the canonical symplectic form $\sigma$ is exact on
each orbit.}
\vspace{0.1in}\\
\noindent {\bf Remark\ 4.2.}\quad (a) Let $G$ be a real Lie group
with its Lie algebra $\mathfrak {g} $. Let
$\Omega_{\ell}:=Ad^*(G)\ell$ be the coadjoint orbit containing
$\ell$. From now on, we write $Ad^*$ instead of $Ad^*_G$. Then
$\Omega_{\ell}$ is simply connected if a Lie group $G$ is
exponential. We recall that a Lie group $G$ is said to be {\it
exponential} if the exponential mapping $\text{exp}:\mathfrak {g}
\longrightarrow G$ is a diffeomorphism. But if $G$ is solvable,
$\Omega_{l}$ is not necessarily simply connected.

\smallskip

\noindent (b) Let $G$ be a connected, simply connected solvable
Lie group and for $\ell\in \mathfrak {g} ^*,$ we let $G_{\ell}$ be
the stabilizer at $\ell.$ Then for any $\ell\in \mathfrak {g} ^*,$
we have $\pi_1(\Omega_{l})\cong G_{\ell}/G_{\ell}^0,$ where
$G_{\ell}^0$ denotes the identity component of $G_{\ell}$ in $G.$
\def\Gl{G_{\ell}}
\def\lg{\ell \in \mathfrak {g} ^*}
\def\Adl{Ad^*(G)\ell}
\def\rig{\mathfrak {g} _{\text{rigg}}^*}
\def\OG{{\mathcal {O} }_{\text{rigg}}(G)}
\def\lchi{(\ell,\chi)}

\medskip

\indent Let $G$ be a Lie group with Lie algebra $\mathfrak {g} $.
A pair $(\ell,\chi)$ is called a {\it rigged\ momentum} if $\lg$,
and $\chi$ is a unitary character of $\Gl$ such that $d\chi_e=2\pi
i \ell|_{\mathfrak {g} _{\ell}},$ where $d\chi_e$ denotes the
differential of $\chi$ at the identity element $e$ of $G_{\ell}$.
We denote by $\rig$ the set of all rigged momenta. Then $G$ acts
on $\rig$ by
\begin{equation}
g\cdot \lchi:=(Ad^*(g)\ell, \chi\circ I_{g^{-1}})=(\ell\circ
Ad(g^{-1}),\chi\circ I_{g^{-1}})
\end{equation}
for all $g\in G$ and $\lchi\in \rig.$ Here $I_g$ denotes the inner
automorphism of $G$ defined by $I_g(x)=gxg^{-1}\,(x\in G).$ We
note that $\chi\circ I_{g^{-1}}$ is a unitary character of
$G_{\Adl}=gG_{\ell}g^{-1}.$ We denote by $\OG$ the set of all
orbits in $\rig$ under the action (4.1).
\vspace{0.1in}\\
\noindent {\bf Proposition\ 4.3.}\quad {\it Let $G$ be a
connected, simply connected solvable Lie group. Then the following
(a) and
(b) hold.\\
\indent (a) The $G$-action commutes with the natural projection
\begin{equation}
\pi :\rig\longrightarrow \mathfrak {g} ^*\quad\quad \lchi\mapsto
\ell.
\end{equation}
 \indent (b) For a solvable Lie group
$G$ of type I, the projection $\pi$ is surjective and the fiber
over a point $\lg$ is a torus of dimension equal to the first
Betti number $b_1(\Omega_{\ell})$ of $\Omega_{\ell}.$}

\smallskip

Now we mention the main theorem obtained by L. Auslander and B.
Kostant in \cite{AK}.
\vspace{0.1in}\\
\noindent {\bf Theorem 4.4 (Auslander-Kostant).}\quad {\it Let $G$
be a connected, simply connected solvable Lie group of type I.
Then there is a natural bijection between the unitary dual ${\hat
G}$ and the orbit space $\OG$. The correspondence between ${\hat
G}$ and $\OG$ is given as follows. Let $\lchi\in\rig$. Then there
always exists a complex subalgebra $\mathfrak {p}$ of $\mathfrak
{g}_{\mathbb{C}}$ subordinate to $\ell$. We let
$L(G,\ell,\chi,\mathfrak {p})$ be the space of complex valued
functions $\phi$ on $G$ satisfying the following conditions
\begin{equation}
\phi(hg)=\chi(h)\phi(g),\quad h\in \Gl
\end{equation}
and
\begin{equation}
(L_X+2\pi i <\ell,X>)\phi=0,\quad X\in \mathfrak {p},
\end{equation}
where $L_X$ is the right invariant complex vector field on $G$
defined by $X\in\mathfrak {g} _{\mathbb{C}}.$ Then we have the
representation $T$ of $G$ defined by
\begin{equation}
(T(g_1)\phi)(g):=\phi(gg_1),\quad g,g_1\in G.
\end{equation}
We can show that under suitable conditions on $\mathfrak {p}$
including the Pukanszky condition
\begin{equation}
p^{-1}(p(\ell))=\ell +\mathfrak {p}^{\perp}\subset \Omega_{\ell}
\end{equation}
and the condition
\begin{equation}
\text{codim}_{\mathbb{C}}\,\mathfrak {p}={\frac
12}\text{rank}\,B_{\ell},
\end{equation}
the representation $T$ is irreducible and its equivalence class
depends only on the rigged orbit $\Omega$ containing $\lchi$. Here
$p:\mathfrak {g} _{\mathbb{C}}^*\longrightarrow \mathfrak {p}^*$
denotes the natural projection of $\mathfrak {g} _{\mathbb {C}}^*$
onto $\mathfrak {p}^*$ dual to the inclusion $\mathfrak
{p}\hookrightarrow \mathfrak {g} _{\mathbb {C}}.$ We denote by
$T_{\Omega}$ the representation $T$ of $G$ obtained from
$\lchi\in\rig$ and $\mathfrak {p}$. The correspondence between
$\OG$ and ${\hat G}$ is given by
$$\lchi\in \Omega \mapsto T_{\Omega}.$$}
\vspace{0.1in}\\
\noindent {\bf Definition\ 4.5.}\quad Let $H$ be a closed subgroup
of a Lie group $G$. We say that a rigged orbit $\Omega'\in
{\mathcal {O} }_{\text{rigg}}(H)$ lies under a rigged orbit
$\Omega\in \OG$\,(or equivalently, $\Omega$ lies over $\Omega'$)
if there exist rigged momenta $\lchi\in \Omega$ and
$(\ell',\chi')\in\Omega'$ such that the following conditions are
satisfied:
\begin{equation}
p(\ell)=\ell',\quad \chi=\chi'\quad\text{on}\ H\cap \Gl.
\end{equation}
 We define the {\it sum\ of\ rigged\ orbits} $\Omega_1$
and $\Omega_2$ as the set of all $\lchi\in \OG$ for which there
exist $(\ell_i,\chi_i)\in\Omega_i,\ i=1,2$, such that
\begin{equation}
\ell=\ell_1+\ell_2,\quad \chi=\chi_1\chi_2\ \ \text{on}\
G_{\ell_1}\cap G_{\ell_2}.
\end{equation}
\quad I. M. Shchepochkina \cite{Sh1}-\cite{Sh2} proved the following.
\vspace{0.1in}\\
\noindent {\bf Theorem\ 4.6.}\quad {\it Let $G$ be a connected,
simply connected solvable Lie group and $H$ a closed subgroup of
$G$. Then

\smallskip

\indent (a) The spectrum of $\text{Ind}_H^G S_{\Omega'}$ consists
of those $T_{\Omega}$ for which $\Omega$ lies over $\Omega'$,
where $S_{\Omega'}$ is an irreducible unitary representation of
$H$ corresponding to a rigged orbit $\Omega'$ in $\mathfrak
{h}_{\text{rigg}}^*$ by Theorem 4.4.

\smallskip

\indent (b) The spectrum of $\text{Res}^G_H T_{\Omega}$ consists
of those $S_{\Omega'}$ for which $\Omega'$ lies under $\Omega.$

\smallskip

\indent (c) The spectrum of $T_{\Omega_1}\otimes T_{\Omega_2}$
consists of those $T_{\Omega}$ for which $\Omega$ lies in
$\Omega_1+\Omega_2.$}
\vspace{0.2in}\\
\noindent{\bf 5. The Obstacle for the Orbit Method}
\setcounter{equation}{0}
\renewcommand{\theequation}{5.\arabic{equation}}
\vspace{0.1in}\\
\indent In this section, we discuss the case where the
correspondence between irreducible unitary representations and
coadjoint orbits breaks down. If $G$ is a compact Lie group or a
semisimple Lie group, the correspondence breaks down.

\smallskip

\def\SG{{\mathcal {S} }(G)}
\def\SGP{{\mathcal {S} }'(G)}
\def\CF{{\Cal F}}
\indent First we collect some definitions. Let $G$ be a Lie group
with Lie algebra $\mathfrak {g} $ and let ${\mathcal {S} }(G)$ be
the Schwarz space of rapidly decreasing functions on $G$. We
define the Fourier transform ${\mathcal {F}}_f$ for $f\in
{\mathcal {S}}(G)$ by
\begin{equation}
\mathcal {F}_f(\ell)=\int_{\mathfrak {g}
}f(\text{exp}\,X)\,e^{2\pi i \lambda(X)}dX,\quad \lambda \in
{\mathfrak {g} }^*.
\end{equation}
Then (5.1) is a well-defined function on $\mathfrak {g} ^*.$ As
usual, we define the Fourier transform $\mathcal{F}_{\chi}$ of a
distribution $\chi\in \SGP$ by
\begin{equation}
<\mathcal{F}_{\chi},\mathcal{F}_f>=<\chi, f>,\quad f\in\SG.
\end{equation}
For $f\in \SG$ and an irreducible unitary representation $T$ of
$G$, we put
$$T(f)=\int_{\mathfrak {g} }f(\text{exp}\,X)\,T(\text{exp}\,X)dX.$$
Then we can see that for an irreducible unitary representation of
a {\it nilpotent} Lie group $G$, we obtain the following formula
\begin{equation}
\text{tr}\,T(f)=\int_{\Omega}\mathcal{F}_f(\lambda)d_{\Omega}\lambda,
\end{equation}
where $\Omega$ is the coadjoint orbit in $\mathfrak {g} ^*$
attached to $T$ under the Kirillov correspondence and
$d_{\Omega}\lambda$ is the measure on $\Omega$ with dimension $2k$
given by the form ${1\over{k!}}\, B_{\Omega}\wedge\cdots\wedge
B_{\Omega}\,(k$ factors) with the canonical symplectic form
$B_{\Omega}$ on $\Omega.$
\vspace{0.1in}\\
\noindent {\bf Definition\ 5.1.}\quad Let $G$ be a Lie group with
Lie algebra $\mathfrak {g} .$ A coadjoint orbit $\Omega$ in
$\mathfrak {g} ^*$ is called {\it integral} if the two dimensional
cohomology class defined by the canonical two form $B_{\Omega}$
belongs to $H^2(\Omega,\mathbb{Z})$, namely, the integral of
$B_{\Omega}$ over a two dimensional cycle in $\Omega$ is an
integer.
\vspace{0.2in}\\
\noindent{\bf 5.1.\ Compact Lie Groups}
\vspace{0.1in}\\
\indent Let $G$ be a connected and simply connected Lie group with
Lie algebra $\mathfrak {g} $. Then the $G$-orbits in $\mathfrak
{g} ^*$ are simply connected and have K{\"a}hler structures(not
unique). These K{\"a}hler manifolds are called {\it flag\
manifolds} because their elements are realized in terms of flags.
Let $T$ be a maximal abelian subgroup of $G$. Then $X=G/T$ is
called the {\it full\ flag\ manifold} and other flag manifolds are
called {\it degenerate} ones. From the exact sequence
$$\cdots\longrightarrow
\pi_k(G)\longrightarrow \pi_k(X)\longrightarrow
\pi_{k-1}(T)\longrightarrow\pi_{k-1}(G)\longrightarrow\cdots$$ and
the fact that $\pi_1(G)=\pi_2(G)$ under the assumption that $G$ is
simply connected, we obtain
\begin{equation}
H_2(X,\mathbb {Z})\cong \pi_2(X)\cong \pi_1(T)\cong
\mathbb{Z}^{\text{dim}\,T}.
\end{equation}
Let $\Omega$ be a coadjoint orbit in $\mathfrak {g} ^*$. We
identify $\mathfrak {g} ^*$ with $\mathfrak {g} $ and $\mathfrak
{g} _{\mathbb{C}}^*$ with $\mathfrak {g} _{\mathbb{C}}$ via the
Killing form so that $\mathfrak{t}^*_{\mathbb {C}}$ goes to
$\mathfrak {t}_{\mathbb {C}}$ and the weight $P\subset \mathfrak
{t}^*_{\mathbb {C}}$ corresponds to a lattice in $i\mathfrak
{t}^*\subset i\mathfrak {g} ^*\cong i\mathfrak {g} .$ Then
$\Omega\cap \mathfrak {t}^*$ is a finite set which forms a single
$W$-orbit, where $W$ is the Weyl group defined as
$W=N_G(T)/Z_G(T).$
\vspace{0.1in}\\
\noindent {\bf Proposition\ 5.2.}\quad {\it Let $\Omega_{\lambda}$
be the orbit passing through the point $i\lambda\in\mathfrak
{t}^*.$ Then

\smallskip

\indent (1) the orbit $\Omega_{\lambda}$ is integral if and only
if $\lambda\in P.$

\smallskip

\indent (2) $\text{dim}\,\Omega_{\lambda}$ is equal to the number
of roots non-orthogonal to $\lambda$.}

\bigskip

\indent Let $\Omega$ be an integral orbit of maximal dimension in
$\mathfrak {g} ^*.$ Let $\lambda\in\Omega$ and let $\mathfrak {h}$
be a positive admissible polarization for $\lambda$. Here the
admissibility for $\lambda$ means that $\mathfrak {h}$ satisfies
the following
conditions:\\
\indent (A1) $\mathfrak {h}$ is invariant under the action of
$G_{\lambda},$\\
\indent (A2) $\mathfrak {h}+{\bar{\mathfrak {h}}}$ is a subalgebra
of $\mathfrak {g} _{\mathbb{C}}$.\\
\noindent Here $G_{\lambda}$ denotes the stabilizer of $G$ at
$\lambda$.

\smallskip

\indent Let $\chi_{\lambda}$ be the unitary character of
$G_{\lambda}$ defined by
\begin{equation}
\chi_{\lambda}(\text{exp}\,X)=e^{2\pi i\lambda (X)},\quad X\in
\mathfrak {g} _{\lambda},
\end{equation}
 where $\mathfrak {g} _{\lambda}$  is the Lie algebra of $G_{\lambda}.$ We
note that $G_{\lambda}$ is connected. Let $L_{\lambda}$ be the
hermitian line bundle over $\Omega=G/G_{\lambda}$ defined by the
unitary character $\chi_{\lambda}$ of $G_{\lambda}$. Then $G$ acts
on the space $\Gamma(L_{\lambda})$ of holomorphic sections of
$L_{\lambda}$ as a representation of $G$. A. Borel and A. Weil
proved that $\Gamma(L_{\lambda})$ is non-zero and is an
irreducible unitary representation of $G$ with highest weight
$\lambda$. This is the so-called {\it Borel}-{\it Weil\ Theorem}.
Thereafter this theorem was generalized by R. Bott in the late
1950s as follows.
 \vspace{0.1in}\\
 \noindent {\bf Theorem\ 5.3 (R. Bott).}\quad {\it Let $\rho$ be the half sum of positive roots of the
root system for $(\mathfrak {g} _{\mathbb{C}},\mathfrak
{t}_{\mathbb{C}})$. Then the cohomology space $H^k(X,L_{\lambda})$
is non-zero precisely when
\begin{equation}
\rho-i\lambda=w(\mu+\rho)
\end{equation}
for some $\mu\in P_+,\ w\in W$ and $k=l(w),$ the length of $w$. In
this case the representation of $G$ in $H^k(\Omega,L_{\lambda})$
is equivalent to $\pi_{-i\mu}.$}

 \smallskip

 \indent We note that the Borel-Weil
Theorem strongly suggests relating $\pi_{\lambda}$ to
$\Omega_{\lambda}$ and, on the other hand, the Bott's Theorem
suggests the correspondence $\pi_{\lambda}\leftrightarrow
\Omega_{\lambda + \rho}$ which is a bijection between the unitary
dual ${\hat G}$ of $G$ and the set of all integral orbits of
maximal dimension. It is known that
\begin{equation}
\text{dim}\,\pi_{\lambda}=\text{vol}(\Omega_{\lambda + \rho}).
\end{equation}
The character formula (5.3) is valid for a compact Lie group $G$
and provides an integral representation of the character\,:
\begin{equation}
\chi_{\lambda}(\text{exp}\,X)={1\over
{p(X)}}\,\int_{\Omega_{\lambda}}e^{2\pi
i\lambda(X)}d_{\Omega}\lambda.
\end{equation}
\indent In 1990 N.J. Wildberger \cite{Wi} proved the following.
\vspace{0.1in}\\
\noindent {\bf Theorem\ 5.4.}\quad {\it Let
$\Phi:C^{\infty}(\mathfrak {g} )'\longrightarrow C^{\infty}(G)'$
be the transform defined by
\begin{equation}
<\Phi(\nu),f>=<\nu, p\cdot (f\circ \text{exp})>,\quad \nu\in
C^{\infty}(\mathfrak {g} )',\ f\in C^{\infty}(G)'.
\end{equation}
Then for $Ad(G)$-invariant distributions the convolution operators
on $G$ and $\mathfrak {g} $ are related by the transform above\,:}
\begin{equation}
\Phi(\mu)*_G \Phi(\nu)=\Phi(\mu *_{\mathfrak {g} }\nu).
\end{equation}
 \indent The above theorem says that $\Phi$ straightens the group
convolution, turning it into the abelian convolution on $\mathfrak
{g} $. This implies the following geometric fact.
\vspace{0.1in}\\
\noindent {\bf Corollary\ 5.5.}\quad {\it For any two coadjoint
orbits $\Omega_1,\Omega_2\subset \mathfrak {g} $, we let
$C_1=\text{exp}\,\Omega_1$ and $C_2=\text{exp}\,\Omega_2$. Then
the following holds.}
\begin{equation}
C_1\cdot C_2\subset \text{exp}\,(\Omega_1+\Omega_2).
\end{equation}
\vspace{0.2in}\\
\noindent{\bf 5.2.\ Semisimple Lie Groups}
\vspace{0.1in}\\
\indent The unitary dual ${\hat G}_u$ of a semisimple Lie group
$G$ splits into different series, namely, the principal series,
degenerate series, complimentary series, discrete series and so
on. These series may be attached to different types of coadjoint
orbits. The principal series were defined first for complex
semisimple Lie groups and for the real semisimple Lie group $G$
which admit the {\it split} Cartan subalgebra $\mathfrak
{h}\subset \mathfrak {g} .$ These series are induced from
characters of the Borel subgroup $B\supset H=\text{exp}\,\mathfrak
{h}$. The degenerate series are obtained by replacing $B$ by a
parabolic subgroup $P\supset B.$ All these series are in a perfect
correspondence with the family of coadjoint orbits which have a
non-empty intersection with $\mathfrak {h}.$ An irreducible
unitary representation $\pi$ of $G$ is said to be a {\it discrete\
series} if it occurs as a direct summand in the regular
representation $R$ of $G$ on $L^2(G,dg)$. According to
Harish-Chandra, if $G$ is a real semisimple Lie group, ${\hat
G}_d\neq 0$ if and only if $G$ has a compact Cartan subgroup. Here
${\hat G}_d$ denotes the set of equivalent classes of discrete
series of $G$. There is an interesting complimentary series of
representations which are {\it not \ weakly\ contained} in the
regular representation $R$ of $G$. These can be obtained from the
principal series and degenerate series by analytic
continuation.

\smallskip

\indent The principal series are related to the semisimple orbits.
On the other hand, the nilpotent orbits are related to the
so-called {\it unipotent} representations if they exist. In fact,
the hyperbolic orbits are related to the representations obtained
by the {\it parabolic\ induction} and the elliptic orbits are
connected to the representations obtained by the {\it
cohomological\ parabolic\ induction}. The notion of unipotent
representations are not still well defined and hence not
understood well. Recently J.-S. Huang and J.-S. Li \cite{HuL}
attached unitary representations to spherical nilpotent orbits for
the real orthogonal and symplectic groups. The study of unipotent
representations is under way. For some results and conjectures on
unipotent representations, we refer to \cite{A}, [16], \,\
\cite{HuL} and \cite{Vo3}-\cite{Vo4}.
\vspace{0.2in}\\
\noindent{\bf 6. Nilpotent Orbits and the Kostant-Sekiguchi
Correspondence}
\setcounter{equation}{0}
\renewcommand{\theequation}{6.\arabic{equation}}
\vspace{0.1in}\\
\indent In this section, we present some properties of nilpotent
orbits for a reductive Lie group $G$ and describe the
Kostant-Sekiguchi correspondence. We also explain the work of D.
Vogan that for a maximal compact subgroup $K$ of $G$, he attaches
a space with a $K$-action to a nilpotent orbit. Most of the
materials in this section are based on the article \cite{Vo6}.
\vspace{0.2in}\\
\noindent{\bf 6.1. Jordan Decomposition}
\vspace{0.1in}\\
\noindent {\bf Definition\ 6.1.1.}\quad Let $GL(n)$ be the group
of nonsingular real or complex $n\times n$ matrices. The {\it
Cartan\ involution} of $GL(n)$ is the automorphism conjugate
transpose inverse\,:
\begin{equation}
\theta(g)=\,^t{\bar g}^{-1},\quad g\in GL(n).
\end{equation}
 A {\it linear\ reductive\ group} is a closed
subgroup $G$ of some $GL(n)$ preserved by $\theta$ and having
finitely many connected components. A {\it reductive} Lie group is
a Lie group ${\tilde G}$ endowed with a homomorphism $\pi:{\tilde
G}\longrightarrow G$ onto a linear reductive group $G$ so that the
kernel of $\pi$ is finite.
\vspace{0.1in}\\
\noindent {\bf Theorem\ 6.1.2 (Cartan\ Decomposition).}\quad {\it
Let ${\tilde G}$ be a reductive Lie group with $\pi:{\tilde
G}\longrightarrow G$ as in Definition 6.1.1. Let
$$K=G^{\theta}=\left\{ g\in G\,\vert\ \theta(g)=g\,\right\}$$ be a
maximal compact subgroup of $G$. We write ${\tilde
K}=\pi^{-1}(K),$ a compact subgroup of ${\tilde G}$, and use
$d\pi$ to identify the Lie algebras of ${\tilde G}$ and $G$. Let
$\mathfrak {p}$ be the $(-1)$-eigenspace of $d\theta$ on the Lie
algebra $\mathfrak {g} $ of $G$. Then the map
\begin{equation}
{\tilde K}\times \mathfrak {p}\longrightarrow {\tilde G},\quad
({\tilde k},X)\mapsto {\tilde k}\cdot \text{exp}\,X,\quad {\tilde
k}\in {\tilde K},\ X\in \mathfrak {p}
\end{equation}
is a diffeomorphism from ${\tilde K}\times \mathfrak {p}$ onto
${\tilde G}.$ In particular, ${\tilde K}$ is maximal among the
compact subgroups of ${\tilde G}$.}

 \smallskip

\indent Suppose $\tilde{G}$ is a reductive Lie group. We define a
map $\theta : \tilde{G} \longrightarrow \tilde{G}$ by
\begin{equation}
\theta(\tilde{k} \cdot \exp X)=\tilde{k} \cdot \exp(-X),\quad
\tilde{k}\in\tilde{K},\ X \in \mathfrak {p}.
\end{equation}
Then $\theta$ is an involution, that is, the Cartan involution of
$\tilde{G}$. The group of fixed points of $\theta$ is $\tilde{K}.$

\smallskip

\indent The following proposition makes us identify the Lie
algebra of a reductive Lie group with its dual space.
\vspace{0.1in}\\
\noindent {\bf Proposition\ 6.1.3.}\quad {\it Let $G$ be a
reductive Lie group. Identify $\mathfrak {g} $ with a Lie algebra
of $n\times n$ matrices (cf. Definition 6.1.1). We define a real
valued symmetric bilinear form on $\mathfrak {g} $ by
\begin{equation}
<X,Y>=\text{Re\,tr}(XY),\quad X,Y\in \mathfrak {g} .
\end{equation}
Then the following (a),(b) and (c) hold: \\
\indent (a) The form $<\, ,\, >$ is invariant under $Ad(G)$ and
the Cartan involution $\theta.$

\smallskip

\indent (b) The Cartan decomposition $\mathfrak {g}
=\mathfrak{k}+\mathfrak {p}$ is orthogonal with respect to the
form $<\, ,\, >$, where $\mathfrak{k}$ is the Lie algebra of
$K$\,(the group of fixed points of $\theta$) which is the
$(+1)$-eigenspace of $\theta:=d\theta$ on $\mathfrak {g} $ and
$\mathfrak {p}$ is the (-1)-eigenspace of $\theta$ on $\mathfrak
{g} $. The form $<\, ,\,>$ is negative definite on $\mathfrak{k}$
and positive definite on $\mathfrak {p}$. And hence the form $<\,
,\,>$ is nondegenerate on $\mathfrak {g}$.

\smallskip

\indent (c) There is a $G$-equivariant linear isomorphism
$$\mathfrak {g} ^*\cong \mathfrak {g} ,\quad \lambda\mapsto X_{\lambda}$$ characterized
by}
\begin{equation}
\lambda(Y)=<X_{\lambda},Y>,\quad Y\in\mathfrak {g} .
\end{equation}
\vspace{0.1in}\\
\noindent {\bf Definition\ 6.1.4.}\quad Let $G$ be a reductive Lie
group with Lie algebra $\mathfrak {g} $ consisting of $n\times n$
matrices. An element $X\in\mathfrak {g} $ is called {\it
nilpotent} if it is nilpotent as a matrix. An element
$X\in\mathfrak {g} $ is called {\it semisimple} if the
corresponding complex matrix is diagonalizable. An element
$X\in\mathfrak {g} $ is called {\it hyperbolic} if it is
semisimple and its eigenvalues are real. An element $X\in\mathfrak
{g} $ is called {\it elliptic} if it is semisimple and its
eigenvalues are purely imaginary.
\vspace{0.1in}\\
\noindent {\bf Proposition\ 6.1.5\,(Jordan\ Decomposition).}\quad
{\it Let $G$ be a reductive Lie group with its Lie algebra
$\mathfrak {g} $ and let $G=K\cdot\text{exp}\,\mathfrak {p}$ be
the Cartan decomposition of $\mathfrak {g} $\,(see Theorem 6.1.2).
Then the
following (1)-(5) hold:\\

\smallskip

\indent (1) Any element $X\in\mathfrak {g} $ has a unique
decomposition
$$X=X_h+X_e+X_n$$ characterized by the conditions that $X_h$ is
hyperbolic, $X_e$ is elliptic, $X_n$ is nilpotent and
$X_h,\,X_e,\,X_n$ commute with each other.

\smallskip

\indent (2) After replacing $X$ by a conjugate under $Ad(G)$, we
may assume that $X_h\in\mathfrak {p},\ X_e\in\mathfrak{k}$ and
that $X_n=E$ belongs to a standard $\mathfrak{s}{\mathfrak
{l}}(2)$ triple. We recall that a triple $\{
H,E,F\}\subset\mathfrak {g} $ is called a standard
$\mathfrak{s}{\frak l}(2)$ triple, if they satisfy the following
conditions
\begin{equation}
\theta(E)=-F,\quad\theta(H)=-H,\quad [H,E]=2E,\quad [E,F]=H.
\end{equation}

\smallskip

\indent (3) The $Ad(G)$ orbits of hyperbolic elements in
$\mathfrak {g} $ are in one-to-one correspondence with the $Ad(K)$
orbits in $\mathfrak {p}$.

 \smallskip

 \indent (4) The $Ad(G)$ orbits of elliptic
elements in $\mathfrak {g} $ are in one-to-one correspondence with
the $Ad(K)$ orbits in $\mathfrak{k}$.

 \smallskip

 \indent (5) The $Ad(G)$ orbits of nilpotent orbits are in one-to-one correspondence with the
$Ad(K)$ orbits of standard $\mathfrak{s}{\mathfrak {l}}(2)$
triples in $\mathfrak {g} $.}
\vspace{0.2in}\\
\noindent{\bf 6.2. Nilpotent Orbits}
\def\sl{{\mathfrak {s}}{\mathfrak {l}}(2,\mathbb  {C})}
\def\SL{SL(2,\mathbb {C})}
\def\SR{SL(2,\mathbb{R})}
\def\th{\theta}
\def\s{\sigma}
\def\sr{{\frak s}{\frak l}(2,\mathbb{R})}
\def\SL{SL(2,\mathbb  {C})}
\vspace{0.1in}\\
\indent Let $G$ be a real reductive Lie group with Lie algebra
$\mathfrak {g} $. Let $\mathfrak {g} _{\mathbb  {C}}$ be the
complexification of $\mathfrak {g} $. We consider the complex
special linear group $SL(2,\mathbb {C})$ which is the
complexification of $SL(2,\mathbb{R}).$ We define the involution
$\theta_{0} : SL(2, \mathbb{C}) \longrightarrow SL(2, \mathbb{C})$
by
\begin{equation}
\theta_0(g)=\,^tg^{-1}, \quad  g \in SL(2, \mathbb{C}).
\end{equation}
\noindent We denote its differential by the same letter
\begin{equation}
\theta_0(Z)=-\,^tZ,\quad Z\in \mathfrak{s}{\mathfrak {l}}(2,
\mathbb{C}).
\end{equation}
\noindent The complex conjugation $\s_0$ defining the real form
$\SR$ is just the complex conjugation of matrices.

\smallskip

\indent We say that a triple $\{ H,X,Y\}$ in a real or complex Lie
algebra is a {\it standard triple} if it satisfies the following
conditions:
\begin{equation}
[H,X]=2X,\quad [H,Y]=-2Y,\quad [X,Y]=H.
\end{equation}
\noindent We call the element $H$ (resp. $X,\,Y$) a {\it neutral}
(resp. {\it nilpositive,\ nilnegative}) element of a standard
triple $\{ H,X,Y\}$.\\
\indent We consider the standard basis $\{ H_0,E_0,F_0\}$ of
$\mathfrak{s}{\mathfrak {l}}(2, \mathbb{R})$ given by
\begin{equation}
H_0=\begin{pmatrix} 1 & 0 \\ 0 & -1
\end{pmatrix} ,\quad E_0=\begin{pmatrix} 0 & 1 \\ 0 & 0
\end{pmatrix},\quad F_0=\begin{pmatrix} 0 & 0 \\ 1 & 0
\end{pmatrix}.
\end{equation}
\noindent Then they satisfy
\begin{equation}
[H_0,E_0]=2E_0,\quad [H_0,F_0]=-2F_0,\quad [E_0,F_0]=H_0
\end{equation}
 and
\begin{equation}
\th_0(H_0)=-H_0,\quad \th_0(E_0)=-F_0,\quad \th_0(F_0)=-E_0.
\end{equation}
\noindent We fix a Cartan decomposition $\mathfrak {g}
=\mathfrak{k}+\mathfrak {p}$ and let $\th$ be the corresponding
Cartan involution. We say that a standard triple $\{ H,X,Y\}$ is a
{\it Cayley\ triple} in $\mathfrak {g} $ if it satisfies the
conditions:
\begin{equation}
\th(H)=-H,\quad \th(X)=-Y,\quad \th(Y)=-X.
\end{equation}
\noindent According to (6.11) and (6.12), the triple $\{
H_0,E_0,F_0\}$ is a Cayley triple in $\sr$.

\smallskip

\indent For a real Lie algebra $\mathfrak {g} $, we have the
following theorems.
\vspace{0.1in}\\
\noindent {\bf Theorem\ 6.2.1.}\quad {\it Given a Cartan
decomposition $\th$ on $\mathfrak {g} $, any triple $\{ H,X,Y\}$
in $\mathfrak {g} $ is conjugate under the adjoint group $Ad(G)$
to a Cayley triple $\{ H',X',Y'\}$ in $\mathfrak {g} $.}
\vspace{0.1in}\\
\noindent {\bf Theorem\ 6.2.2 (Jacobson-Morozov).}\quad {\it Let
$X$ be a nonzero nilpotent element in $\mathfrak {g} .$ Then there
exists a standard triple $\{ H,X,Y\}$ in $\mathfrak {g} $ such
that $X$ is nilpositive.}
\vspace{0.1in}\\
\noindent {\bf Theorem\ 6.2.3 (Kostant).}\quad {\it Any two
standard triples $\{ H,X,Y\},\ \{ H',X,Y'\}$ in $\mathfrak {g} $
with the same nilpositive element $X$ are conjugate under $G^X$,
the
centralizer of $X$ in the adjoint group of $G$.}\\
\indent Let $\{ H,X,Y\}$ be a Cayley triple in $\mathfrak {g} .$
We are going to look for a semisimple element in $\mathfrak {g} .$
For this, we need to introduce an auxiliary standard triple
attached to a Cayley triple. We put
\begin{equation}
H'=i(X-Y),\quad X'={\frac 12}(X+Y+iH),\quad Y'={\frac 12}(X+Y-iH).
\end{equation}
 Then the triple $\{ H',X',Y'\}$ in $\mathfrak {g} _{\mathbb
{C}}$ is a standard triple, called the {\it Cayley\ transform} of
$\{ H,X,Y\}$.

\smallskip

\indent Since $$ \theta(H')=H',\quad \theta(X')=-X',\quad
\theta(Y')=-Y',$$ we have
\begin{equation}
H'\in\mathfrak{k}_{\mathbb {C}}\quad \text{and}\quad
X',Y'\in\mathfrak {p}_{\mathbb{C}}.
 \end{equation}
Therefore the subalgebra $\mathbb {C} <H',X',Y'>$ of $\mathfrak
{g} c$ spanned by $H',X',Y'$ is stable under the action of $\th.$
A standard triple in $\mathfrak {g} c$ with the property (6.15) is
called {\it normal}.
\vspace{0.1in}\\
\noindent {\bf Theorem\ 6.2.4.}\quad {\it Any nonzero nilpotent
element $X\in\mathfrak {p}_{\mathbb {C}}$ is the nilpositive
element of a normal triple (\,see Theorem 6.2.2).}
\vspace{0.1in}\\
\noindent {\bf Theorem\ 6.2.5.}\quad {\it Any two normal triples
$\{ H,X,Y\},\ \{ H',X,Y'\}$ with the same nilpositive element $X$
is $K_{\mathbb {C}}^X$-conjugate, where $K_{\mathbb  {C}}^X$
denotes the centralizer of $X$ in the complexification $K_{\mathbb
{C}}$ of a maximal compact subgroup $K$ corresponding to the Lie
algebra $\mathfrak{k}.$}
\vspace{0.1in}\\
\noindent {\bf Theorem\ 6.2.6.}\quad {\it Any two normal triples
$\{ H,X,Y\},\ \{ H,X',Y'\}$ with the same neutral element $H$ are
$K_{\mathbb {C}}^H$-conjugate.}
 \vspace{0.1in}\\
 \noindent {\bf Theorem\ 6.2.7\,(Rao).}\quad {\it Any two standard triples $\{
H,X,Y\},\ \{ H',X',Y'\}$ in $\mathfrak {g} $ with $X-Y=X'-Y'$ are
conjugate under $G^{X-Y}$, the centralizer of $X-Y$ in $G$. In
fact, $X-Y$ is a semisimple element which we are looking for.}

 \smallskip

Let ${\mathcal A}_{triple}$ be the set of all $Ad(G)$-conjugacy
classes of standard triples in $\mathfrak {g} .$ Let ${\mathcal
{O} }_ {\mathcal N}$ be the set of all nilpotent orbits in
$\mathfrak {g} .$ We define the map
\begin{equation}
\Omega :{\mathcal A}_{triple}\longrightarrow {\mathcal {O}
}_{\mathcal N}^{\times}:={\mathcal {O} }_{\mathcal N}-\{ 0\}
\end{equation}
 by
 \begin{equation}
 \Omega([\{ H,X,Y\}]):={\mathcal {O} }_X,\quad {\mathcal {O} }_X:=Ad(G)\cdot
X,
\end{equation}
where $[\{ H,X,Y\}]$ denotes the $G$-conjugacy class of a standard
triple $\{ H,X,Y\}$. According to Theorem 6.2.2\,(Jacobson-Morozov
Theorem) and Theorem 6.2.3\,(Kostant's Theorem), the map $\Omega$
is bijective.

\smallskip

We put
 \begin{equation}
h_0=\begin{pmatrix} 0 & i\\ -i & 0
\end{pmatrix},
\quad x_0={\frac 12}
 \begin{pmatrix} 1 & -i\\ -i & -1
 \end{pmatrix},\quad y_0={\frac 12}
 \begin{pmatrix} 1 & i\\ i & -1
 \end{pmatrix}.
\end{equation}
 It is easy to see that the
triple $\{ h_0,x_0,y_0\}$ in $\sl$ is a normal triple. The complex
conjugation $\s_0$ acts on the triple $\{ h_0,x_0,y_0\}$ as
follows:
\begin{equation}
\s_0(h_0)=-h_0,\quad \s_0(x_0)=y_0,\quad \s_0(y_0)=x_0.
\end{equation}
 \indent We introduce some notations. We denote by
$\text{Mor}(\sl,\mathfrak {g}_{\mathbb {C}})$ the set of all
nonzero Lie algebra homomorphisms from $\sl$ to $\mathfrak
{g}_{\mathbb{C}}$. We define
$$\text{Mor}^{\mathbb{R}}(\sl,\mathfrak {g}_{\mathbb{C}})=\{\,\phi\in
\text{Mor}(\sl,\mathfrak {g}_{\mathbb{C}})\,\vert\ \phi\ \text{is\
defined\ over}\ \mathbb{R}\,\},$$ $$\text{Mor}^{\th}(\sl,\mathfrak
{g}_{\mathbb{C}})=\{\,\phi\in \text{Mor}(\sl,\mathfrak
{g}_{\mathbb{C}})\,\vert\ \th\circ \phi=\phi\circ \th_0 \,\},$$
\begin{equation}
\text{Mor}^{\mathbb{R},\th}(\sl,\mathfrak {g}_{\mathbb{C}})=
\text{Mor}^{\mathbb{R}}(\sl,\mathfrak {g}_{\mathbb{C}})\cap
\text{Mor}^{\th}(\sl,\mathfrak {g}_{\mathbb{C}}),
\end{equation}
$$\text{Mor}^{\s}(\sl,\mathfrak {g}_{\mathbb{C}})=\{\,\phi\in
\text{Mor}(\sl,\mathfrak {g}_{\mathbb{C}})\,\vert\ \s\circ
\phi=\phi\circ \s_0 \,\},$$
$$\text{Mor}^{\s,\th}(\sl,\mathfrak {g}_{\mathbb{C}})= \text{Mor}^{\s}(\sl,\mathfrak {g}_{\mathbb{C}})\cap
\text{Mor}^{\th}(\sl,\mathfrak {g}_{\mathbb{C}}).$$ We observe
that $\text{Mor}^{\mathbb{R}}(\sl,\mathfrak {g}_{\mathbb{C}})$ is
naturally isomorphic to $\text{Mor}(\sr,\mathfrak {g} ),$ the set
of all nonzero Lie algebra real homomorphisms from $\sr$ to
$\mathfrak {g} .$
 \vspace{0.1in}\\
 \noindent {\bf Proposition\ 6.2.8.}\quad {\it Suppose $\phi$ be a nonzero Lie algebra homomorphism from $\sl$ to
a complex reductive Lie algebra $\mathfrak {g} _{\mathbb  {C}}$.
Write
$$H=\phi(H_0),\quad E=\phi(E_0),\quad F=\phi(F_0)\quad(\text{see}\
(6.2.10)).$$ Then the following hold.

\smallskip

\noindent (1)$$ \mathfrak {g} _{\mathbb
{C}}=\sum_{k\in\mathbb{Z}}\mathfrak {g} _{\mathbb {C}}(k),$$
 where
  $$\mathfrak {g}_{\mathbb{C}}(k)=\left\{ X\in\mathfrak {g}_{\mathbb{C}}\vert \ [H,X]=kX\,\right\},\quad k\in
\mathbb{Z}.$$

\smallskip

\noindent (2) If we write $$\mathfrak{l}=\mathfrak
{g}_{\mathbb{C}}(0)\quad\text{and}\quad
\mathfrak{u}=\sum_{k>0}\mathfrak {g}_{\mathbb{C}}(k),$$ then
$\mathfrak{q}=\mathfrak{l}+\mathfrak{u}$ is a Levi decomposition
of a parabolic subalgebra of $\mathfrak {g}_{\mathbb{C}}.$\

\smallskip

\noindent (3) The centralizer of $E$ is graded by the
decomposition in (1). More precisely,
$$\mathfrak {g}_{\mathbb{C}}^E=\mathfrak{l}^E+\sum_{k>0}\mathfrak {g}_{\mathbb{C}}(k)^E=\mathfrak{l}^E+\mathfrak{u}^E.$$

\smallskip

\noindent (4) The subalgebra $\mathfrak{l}^E=\mathfrak
{g}_{\mathbb{C}}^{H,E}$ is equal to $\mathfrak
{g}_{\mathbb{C}}^{\phi}$, the centralizer in $\mathfrak
{g}_{\mathbb{C}}$ of the image of $\phi.$ It is a reductive
subalgebra of $\mathfrak {g}_{\mathbb{C}}.$ Consequently the
decomposition in (3) is a Levi decomposition of $\mathfrak
{g}_{\mathbb{C}}^E.$}

 \smallskip

 Parallel results hold if $\{ H,E,F \}$ are replaced by
$$h=\phi(h_0),\quad x=\phi(x_0),\quad y=\phi(y_0).$$
\vspace{0.1in}\\
\noindent {\bf Proposition\ 6.2.9.}\quad {\it Suppose $G$ is a
real reductive Lie group, and let $\phi_{\mathbb{R}}$ be an
element of $\text{Mor}^{\s}(\sl,\mathfrak {g}_{\mathbb{C}})$.
Define $E_{\mathbb{R}},H_{\mathbb{R}},F_{\mathbb{R}}$ by
$$H_{\mathbb{R}}=\phi_{\mathbb{R}}(H_0),\quad E_{\mathbb{R}}=\phi_{\mathbb{R}}(E_0),\quad
F_{\mathbb{R}}=\phi_{\mathbb{R}}(F_0).$$ Then the following hold.

\smallskip

\indent (1) $E_{\mathbb{R}},\ F_{\mathbb{R}}$ are nilpotent, and
$H_{\mathbb{R}}$ is hyperbolic.

\smallskip

\indent (2) If we define $L=G^{H_{\mathbb{R}}}$ to be the
isotrophy group of the adjoint action at $H_{\mathbb{R}}$ and
$U=\text{exp}(\mathfrak{u}\cap\mathfrak {g} ),$ then $Q=LU$ is the
parabolic subgroup of $G$ associated to $H_{\mathbb{R}}.$

\smallskip

\indent (3) The isotrophy group $G^{E_{\mathbb{R}}}$ of the
adjoint action at $E_{\mathbb{R}}$ is contained in $Q$, and
respects the Levi decomposition :
$$G^{E_{\mathbb{R}}}=\left(L^{E_{\mathbb{R}}}\right)\left(U^{E_{\mathbb{R}}}\right).$$
\indent (4) The subgroup $L^{E_{\mathbb{R}}}=G^{H,E_{\mathbb{R}}}$
is equal to $G^{\phi_{\mathbb{R}}},$ the centralizer in $G$ of the
image image of $\phi_{\mathbb{R}}.$ It is a reductive subgroup of
$G$. The The subgroup $U^{E_{\mathbb{R}}}$ is simply connected
unipotent.

 \smallskip

 \indent (5) Suppose that $\phi_{{\mathbb{R}},\th}$ is an
element of $\text{Mor}^{\s,\th}(\sl,\mathfrak {g}_{\mathbb{C}})$.
Then $G^{\phi_{{\mathbb{R}},\th}}$ is stable under the action of
$\th$, and we may take $\th$ as a Cartan involution on this
reductive group. In particular, $G^{E_{\mathbb{R}}}$ and
$G^{\phi_{{\mathbb{R}},\th}}$ have a common maximal compact
subgroup
$$ K^{\phi_{\mathbb{R}, \theta}} = (L \cap K)^{E_{\mathbb{R}}}.$$ }
The above proposition provides good information about the action
of $G$ on the cone ${\mathcal N}_{\mathbb{R}}$ of all nilpotent
orbits in $\mathfrak {g} .$

\smallskip

 The following proposition gives information about the
action of $K_{\mathbb {C}}$ on the cone ${\mathcal {N}}_{\th}$ of
all nilpotent elements in $\mathfrak {p}_{\mathbb  {C}}.$
\vspace{0.1in}\\
\noindent {\bf Proposition\ 6.2.10.}\quad {\it Suppose $G$ is a
real reductive Lie group with Cartan decomposition $G=K\cdot
\text{exp}\,\mathfrak {p}.$ Let $\phi_{\th}\in
\text{Mor}^{\th}(\sl,\mathfrak {g}_{\mathbb{C}}).$ We define
$$h_{\th}=\phi_{\th}(h_0),\quad x_{\th}=\phi_{\th}(x_0),\quad
y_{\th}=\phi_{\th}(y_0).$$ Then we have the following results.

\smallskip

\indent (1) $x_{\th}$ and $y_{\th}$ are nilpotent elements in
$\mathfrak {p}c,\ h_{\th}\in\mathfrak{k}c$ is hyperbolic and
$ih_{\th}\in \mathfrak{k}c$ is elliptic.\\
\indent (2) The parabolic subalgebra
$\mathfrak{q}=\mathfrak{l}+\mathfrak{u}$ constructed as in
Proposition 6.2.8 using $h_{\th}$ is stable under $\theta.$\\
\indent (3) If  we  define $L_K:=\left( K_{\mathbb
C}\right)^{h_{\th}}$ and $U_K=\text{exp}(\mathfrak{u}\cap
\mathfrak{k}c),$ then $Q_K=L_K U_K$ is the parabolic subgroup of
$K_{\mathbb  {C}}$ associated to
$h_{\th}$. \\
\indent (4) $K_{\mathbb {C}}^{x_{\th}}\subset Q_K$ and $K_{\mathbb
{C}}^{x_{\th}}$ respects the Levi decomposition
$$K_{\mathbb  {C}}^{x_{\th}}=\left(L_K^{x_{\th}}\right)\left(
U_K^{x_{\th}}\right).$$ \\
\indent (5) $L_K^{x_{\th}}=K_{\mathbb {C}}^{h_{\th},x_{\th}}$ is
equal to $K_{\mathbb {C}}^{\phi_{\th}}$, the centralizer in
$K_{\mathbb {C}}$ of the image of $\phi_{\th}$. It is a reductive
algebraic subgroup of $K_{\mathbb  {C}}$. The subgroup $U_K^x$ is
simply connected unipotent. In particular, the decomposition of
(4) is a Levi decomposition of $K_{\mathbb  {C}}^x.$ \\
\indent (6) Let $\phi_{\mathbb{R},\th}$ be an element of
$\text{Mor}^{\s,\th}(\sl,\mathfrak {g}_{\mathbb{C}})$. Then
$K_{\mathbb {C}}^{\phi_{\mathbb{R},\th}}$ is stable under $\s$,
and we may take $\s$ as complex conjugation for a compact real
form of this reductive algebraic group. In particular, $K_{\mathbb
{C}}^x$ and  $K_{\mathbb {C}}^{\phi_{\mathbb{R},\th}}$ have a
common maximal compact subgroup}
$$K_{\mathbb {C}}^{\phi_{\mathbb{R},\th}}=L_K^{x_{\th}}\cap K.$$
\vspace{0.2in}\\
\noindent {\bf 6.3. The Kostant-Sekiguchi Correspondence}
\vspace{0.1in}\\
\indent J. Sekiguchi \cite{Se} and B. Kostant (unpublished)
established a bijection between the set of all nilpotent
$G$-orbits in $\mathfrak {g} $ on the one hand and, on the other
hand, the set of all nilpotent $K_{\mathbb {C}}$-orbits in
$\mathfrak {p}_{\mathbb {C}}.$ The detail is as follows.
\vspace{0.1in}\\
\noindent {\bf Theorem\ 6.3.1.}\quad {\it Let $G$ be a real
reductive Lie group with Cartan involution $\th$ and its
corresponding maximal compact subgroup $K$. Let $\s$ be the
complex conjugation on the complexification $\mathfrak
{g}_{\mathbb{C}}$ of $\mathfrak {g} $. Then the following sets are
in one-to-one correspondence.

\smallskip

\indent (a) $G$-orbits on the cone ${\mathcal N}_{\mathbb{R}}$ of
nilpotent elements in $\mathfrak {g} $.

 \smallskip

 \indent (b) $G$-conjugacy
classes of Lie algebra homomorphisms $\phi_{\mathbb{R}}$ in
$\text{Mor}^{\s}(\sl,\mathfrak {g}_{\mathbb{C}})$.

\smallskip

\indent (c) $K$-conjugacy classes of Lie algebra homomorphisms
$\phi_{\mathbb{R},\th}$ in $\text{Mor}^{\s,\th}(\sl,\mathfrak
{g}_{\mathbb{C}})$.

\smallskip

\indent (d) $K_{\mathbb {C}}$-conjugacy classes of Lie algebra
homomorphisms $\phi_{\th}$ in $\text{Mor}^{\th}(\sl,\mathfrak
{g}_{\mathbb{C}})$.

\smallskip

 \indent (e) $K_{\mathbb {C}}$-orbits on the
cone ${\mathcal N}_{\th}$ of nilpotent elements in $\mathfrak
{p}_\mathbb{C}$.}

\medskip

\indent Here $G$ acts on $\text{Mor}^{\s}(\sl,\mathfrak
{g}_{\mathbb{C}})$ via the adjoint action of $G$ in
$\mathfrak{g}$:
\begin{equation}
(g \cdot \phi_{\mathbb{R}})(\zeta)=
Ad(g)(\phi_{\mathbb{R}}(\zeta)),\quad g\in G\   \text{and} \quad
\phi_{\mathbb {R}} \in  \text{Mor}^{\s}(\sl,\mathfrak
{g}_{\mathbb{C}}).
\end{equation}

\noindent Similarly $K$ and $K_{\mathbb {C}}$ act on
$\text{Mor}^{\s,\th}(\sl,\mathfrak {g}_{\mathbb{C}})$ and
$\text{Mor}^{\th}(\sl,\mathfrak {g}_{\mathbb{C}})$ like (6.21)
respectively. The correspondence between (a) and (e) is called the
{\it Kostant}-{\it Sekiguchi\ correspondence} between the
$G$-orbits in ${\mathcal{N}}_{\mathbb{R}}$ and the $K_{\mathbb
{C}}$-orbits in ${\mathcal {N}}_{\th}$. If $\phi_{\mathbb{R},\th}$
is an element in $\text{Mor}^{\s,\th}(\sl,\mathfrak
{g}_{\mathbb{C}})$ as in (c) , then the correspondence is given by
\begin{equation}
E=\phi_{\mathbb{R},\th}(E_0)\leftrightsquigarrow
x=\phi_{\mathbb{R},\th}(x_0).\quad (\text{see}\ (6.18))
\end{equation}
\noindent The proof of the above theorem can be found in \cite{Vo3},
pp. 348-350. \\
\indent M. Vergne \cite{Ve} showed that the orbits $G\cdot E=Ad(G)E$
and $K_{\mathbb {C}}\cdot x=Ad(K_{\mathbb {C}})x$ are
diffeomorphic as manifolds with $K$-action under the assumption
that they are in the Kostant-Sekiguchi correspondence. Here $E$
and $x$ are given by (6.22).
\vspace{0.1in}\\
\noindent {\bf Theorem\ 6.3.2\,(M.\ Vergne).}\quad {\it Suppose
$G=K\cdot \text{exp}(\mathfrak {p})$ is a Cartan decomposition of
a real reductive Lie group $G$ and $E\in \mathfrak {g} ,\
x\in\mathfrak {p}c$ are nilpotent elements. Assume that the orbits
$G\cdot E$ and $K_{\mathbb {C}}\cdot x$ correspond under the
Kostant-Sekiguchi correspondence. Then there is a $K$-equivariant
diffeomorphism from $G\cdot E$ onto $K_{\mathbb {C}}\cdot x$.}
\vspace{0.1in}\\
\noindent {\bf Remark\ 6.3.3.}\quad The Kostant-Sekiguchi
correspondence sends the zero orbit to the zero orbit, and the
nilpotent orbit through the nilpositive element of a Cayley triple
in $\mathfrak {g} $ to the orbit through the nilpositive element
of its Cayley transform.
\vspace{0.1in}\\
\noindent {\bf Remark\ 6.3.4.}\quad Let $G\cdot E$ and $K_{\mathbb
{C}}\cdot x$ be in the Kostant-Sekiguchi correspondence, where
$E\in {\mathcal {N}}_{\mathbb{R}}\subset \mathfrak {g} $ and $x\in
{\mathcal {N}}_{\th}\subset \mathfrak {p}_\mathbb{C}.$ Then the
following hold.

 \medskip

\indent (1) $G_{\mathbb {C}}\cdot E=G_{\mathbb  {C}}\cdot x,$
where $G_{\mathbb  {C}}$ denotes the complexification of $G.$

\medskip

\indent (2) $\text{dim}_{\mathbb {C}}\,(K_{\mathbb {C}}\cdot
x)={\frac 12}\text{dim}_{\mathbb{R}}\,(G\cdot E)={\frac
12}\text{dim}_{\mathbb {C}}\,(G_{\mathbb {C}}\cdot x).$

\medskip

\indent (3) The centralizers $G^E,\ K_{\mathbb  {C}}^x$ have a
common maximal compact subgroup $K^{E,x}$ which is the centralizer
of the span of $E$ and $x$ in $K$.
\vspace{0.1in}\\
\noindent {\bf Remark\ 6.3.5.}\quad Let $\pi$ be an irreducible,
admissible representation of a reductive Lie group $G$. Recently
Schmid and Vilonen gave a new geometric description of the
Kostant-Sekiguchi correspondence\,(cf.\,\cite{ScV1}, Theorem 7.22)
and then using this fact proved that the associated cycle
$\text{Ass}(\pi)$ of $\pi$ coincides with the wave front cycle
$\text{WF}(\pi)$ via the Kostant-Sekiguchi
correspondence\,(cf.\,\cite{ScV2}, Theorem 1.4).
\vspace{0.2in}\\
\noindent {\bf 6.4.\ The Ouantization of the $K$-action}
\vspace{0.1in}\\
\indent It is known that a hyperbolic orbit could be quantized by
the method of a parabolic induction, and on the other hand an
elliptic orbit may also be quantized by the method of
cohomological induction. However, we do not know yet how to
quantize a nilpotent orbit. But D. Vogan attached a space with a
representation of $K$ to a nilpotent orbit.

\smallskip

\indent We first fix a nonzero nilpotent element $\lambda_{n}\in
\mathfrak {g} ^*.$ Let $E$ be the unique element in $\mathfrak {g}
$ given from $\lambda_{n}$ via (6.5). According to
Jacobson-Morosov Theorem or Theorem 6.3.1, there is a non-zero Lie
algebra homomorphism $\phi_{\mathbb{R}}$ from $\sr$ to $\mathfrak
{g} $ with $\phi_{\mathbb{R}}(E_0)=E.$ We recall that $E_0$ is
given by (6.10). $\phi_{\mathbb{R}}$ extends to an element
$\phi_{\mathbb{R}}\in \text{Mor}^{\s}(\sl,\mathfrak
{g}_{\mathbb{C}}).$ After replacing $\lambda_{n}$ by a conjugate
under $G$, we may assume that
$\phi_{\mathbb{R}}=\phi_{\mathbb{R},\th}$ intertwines $\th_0$ and
$\th$, that is, $\phi_{\mathbb{R},\th}\in
\text{Mor}^{\s,\th}(\sl,\mathfrak {g}_{\mathbb{C} }).$

 \smallskip

 \indent We define
 \begin{equation}
  x=\phi_{\mathbb{R},\th}(x_0)\in
\mathfrak {p}_\mathbb{C},\quad x_0={\frac 12}\,\begin{pmatrix} 1 &
-i\\ -i & -1
\end{pmatrix}.
\end{equation}
\noindent The isomorphism in (c), Proposition 6.1.3 associates to
$x$ a linear functional
\begin{equation}
\lambda_{\th}\in \mathfrak {p}_\mathbb{C}^*,\quad
\lambda_{\th}(Y)=<x,Y>,\quad Y\in \mathfrak {g} .
\end{equation}
\noindent We note that the element $\lambda_{\th}$ is not uniquely
determined by $\lambda_{n}$, but the orbit $K_{\mathbb {C}}\cdot
\lambda_{\th}$ is determined by $G\cdot \lambda_{n}.$ According to
Theorem 6.3.2, we get a $K$-equivariant diffeomorphism
\begin{equation}
G\cdot \lambda_{n}\cong K_{\mathbb {C}}\cdot \lambda_{\th}.
\end{equation}
\vspace{0.1in}\\
\noindent {\bf Definition\ 6.4.1.}\quad (1) Let $\mathfrak {g}
_{im}^*$ be the space of purely imaginary-valued linear
functionals on $\mathfrak {g} .$ We fix an element
$\lambda_{im}\in\mathfrak {g} _{im}^*.$ We denote the $G$-orbit of
$\lambda_{im}$ by
\begin{equation}
{\mathcal {O} }_{im}:=G\cdot \lambda_{im}=Ad^*(G)\cdot
\lambda_{im}.
\end{equation}
\noindent We may define an imaginary-valued symplectic form
$\omega_{im}$ on the tangent space
\begin{equation}
T_{\lambda_{im}}({\mathcal {O} }_{im})\cong \mathfrak {g}
/\mathfrak {g} ^{\lambda_{im}},
\end{equation}
\noindent where $\mathfrak {g} ^{\lambda_{im}}$ is the Lie algebra
of the isotropy subgroup $G^{\lambda_{im}}$ of $G$ at
$\lambda_{im}.$ We denote by $Sp(\omega_{im})$ the group of
symplectic real linear transformations of the tangent space
(6.27). Then the isotropy action gives a natural homomorphism
\begin{equation}
j:G^{\lambda_{im}}\longrightarrow Sp(\omega_{im}).
\end{equation}
\noindent On the other hand, we let $Mp(\omega_{im})$ be the
metaplectic group of $Sp(\omega_{im})$. That is, we have the
following exact sequence
\begin{equation}
1\longrightarrow \{ 1,\epsilon\}\longrightarrow
Mp(\omega_{im})\longrightarrow Sp(\omega_{im})\longrightarrow 0.
\end{equation}
\noindent Pulling back (6.29) via (6.28), we have the so-called
{\it metaplectic\ double\ cover} of the isotropy group
$G^{\lambda_{im}}$:
\begin{equation}
1\longrightarrow \{ 1,\epsilon\}\longrightarrow {\tilde
G}^{\lambda_{im}}\longrightarrow G^{\lambda_{im}}.
\end{equation}
\noindent That is, ${\tilde G}^{\lambda_{im}}$ is defined by
\begin{equation}
{\tilde G}^{\lambda_{im}}=\{ (g,m)\in G^{\lambda_{im}}\times
Mp(\omega_{im})\vert\ j(g)=p(m)\,\}.
\end{equation}
\noindent A representation $\chi$ of ${\tilde G}^{\lambda_{im}}$
is called {\it genuine} if $\chi(\epsilon)=-I.$ We say that $\chi$
is {\it admissible} if it is genuine, and the differential of
$\chi$ is a multiple of $\lambda_{im}$: namely, if
\begin{equation}
\chi(\text{exp}\,x)=\text{exp}(\lambda_{im}(x))\cdot I,\quad x\in
\mathfrak {g} ^{\lambda_{im}}.
\end{equation}
\noindent If admissible representations exist, we say that
$\lambda_{im}$ (or the orbit ${\mathcal {O} }_{im}$) is {\it
admissible}. A pair $(\lambda_{im},\chi)$ consisting of an element
$\lambda_{im}\in \mathfrak {g} _{im}^*$ and an irreducible
admissible representation $\chi$ of ${\tilde G}^{\lambda_{im}}$ is
called an {\it admissible} $G$-{\it orbit\ datum}. Two such are
called {\it equivalent} if they are conjugate by $G$.

\smallskip

\indent We observe that if $G^{\lambda_{im}}$ has a finite number
of connected components, an irreducible admissible representation
of ${\tilde G}^{\lambda_{im}}$ is unitarizable. The notion of
admissible $G$-orbit data was introduced by M. Duflo \cite{Du4}.

\medskip

 \indent (2) Suppose $\lambda_{\th}\in \mathfrak {p}_{\mathbb{C}}^*$ is a
non-zero nilpotent element. Let $K_{\mathbb {C}}^{\lambda_{\th}}$
be the isotropy subgroup of $K_{\mathbb {C}}$ at $\lambda_{\th}.$
Define $2\rho$ to be the algebraic character of $K_{\mathbb
{C}}^{\lambda_{\th}}$ by which it acts on the top exterior power
of the cotangent space at $\lambda_{\th}$ to the orbit:
\begin{equation}
2\rho(k):=\text{det}\left(
Ad^*(k)\vert_{(\mathfrak{k}_{\mathbb{C}}/\mathfrak{k}_{\mathbb{C}}^{\lambda_{\th}})^*}\right),\quad
k\in K_{\mathbb  {C}}^{\lambda_{\th}}.
\end{equation}
\noindent The differential of $2\rho$ is a one-dimensional
representation of $\mathfrak{k}_{\mathbb{C}}^{\lambda_{\th}}$,
which we denote also by $2\rho.$ We define $\rho\in
(\mathfrak{k}_{\mathbb{C}}^{\lambda_{\th}})^*$ to be the half of
$2\rho.$ More precisely,
\begin{equation}
\rho(Z)={\frac 12}\,\text{tr}\left(
ad^*(Z)\vert_{(\mathfrak{k}_{\mathbb{C}}/\mathfrak{k}_{\mathbb{C}}^{\lambda_{\th}})^*}\right),\quad
Z\in \mathfrak{k}_{\mathbb{C}}^{\lambda_{\th}}.
\end{equation}
\noindent A {\it nilpotent\ admissible} $K_{\mathbb  {C}}$-{\it
orbit\ datum} at $\lambda_{\th}$ is an irreducible algebraic
representation $(\tau,V_{\tau})$ of $K_{\mathbb
{C}}^{\lambda_{\th}}$ whose differential is equal to $\rho\cdot
I_{\tau},$ where $I_{\tau}$ denotes the identity map on
$V_{\tau}$. The nilpotent element $\lambda_{\th}\in \mathfrak
{p}_{\mathbb{C}}^*$ is called {\it admissible} if a nilpotent
admissible $K_{\mathbb {C}}$-orbit datum at $\lambda_{\th}$
exists. Two such data are called {\it equivalent} if they are
conjugate.
\vspace{0.1in}\\
\noindent {\bf Theorem\ 6.4.2\,(J.\ Schwarz).}\quad {\it Suppose
$G$ is a real reductive Lie group, $K$ is a maximal compact
subgroup, and $K_{\mathbb {C}}$ is its complexification. Then
there is a natural bijection between equivalent classes of
nilpotent admissible $G$-orbit data and equivalent classes of
nilpotent admissible $K_{\mathbb {C}}$-orbit data.}

\smallskip

\indent Suppose $\lambda_n\in\mathfrak {g} ^*$ is a non-zero
nilpotent element. Let $\lambda_{\th}\in \mathfrak
{p}_{\mathbb{C}}^*$ be a nilpotent element which corresponds under
the Kostant-Sekiguchi correspondence. We fix a nilpotent
admissible $K_{\mathbb {C}}$-orbit datum $(\tau,V_{\tau})$ at
$\lambda_{\th}.$ We let
\begin{equation}
{\mathcal V}_{\tau}:=K_{\mathbb {C}}\times_{K_{\mathbb
{C}}^{\lambda_{\th}}} V_{\tau}
\end{equation}
\noindent be the corresponding algebraic vector bundle over the
nilpotent orbit $K_{\mathbb {C}}\cdot \lambda_{\th}\cong
K_{\mathbb {C}}/K_{\mathbb {C}}^{\lambda_{\th}}.$ We note that a
complex structure on ${\mathcal V}_{\tau}$ is preserved by $K$ but
not preserved by $G.$ We now {\it assume} that the boundary of
$\overline{ K_{\mathbb {C}}\cdot \lambda_{\th} }$ (that is,
$\overline{ K_{\mathbb {C}}\cdot \lambda_{\th} }-K_{\mathbb
{C}}\cdot \lambda_{\th}$) has a complex codimension at least two.
We denote by denote by $X_K(\lambda_n,\tau)$ the space of
algebraic sections of ${\mathcal V}_{\tau}$. Then
$X_K(\lambda_n,\tau)$ is an algebraic representation of
$K_{\mathbb {C}}$. That is, if $k_1\in K_{\mathbb {C}}$ and $s\in
X_K(\lambda_n,\tau),$ then $(k_1\cdot
s)(k\lambda_{\th})=s((k_1^{-1}k)\lambda_{\th}).$ We call the
representation $(K_{\mathbb {C}},\,X_K(\lambda_n,\tau))$ of
$K_{\mathbb {C}}$ the {\it quantization\ of\ the} $K$-{\it action}
on $G\cdot \lambda_n$ for the admissible orbit datum
$(\tau,V_{\tau}).$

\smallskip

\indent What this definition amounts to is a desideratum for the
quantization of the $G$-action on $G\cdot \lambda_{\th}$. That is,
whatever a unitary representation $\pi_G(\lambda_n,\tau)$ we
associate to these data, we hope that we have
\begin{equation}
K\!-\!\text{finite\ part\ of}\ \pi_G(\lambda_n,\tau)\cong
X_K(\lambda_n,\tau).
\end{equation}
\noindent When $G$ is a complex Lie group, the coadjoint orbit is
a complex symplectic manifold and hence of real dimension $4m$.
Consequently the codimension condition is automatically satisfied
in this case.
\vspace{0.1in}\\
\noindent {\bf Remark\ 6.4.3.}\quad Nilpotent admissible orbit
data may or may not exist. When they exist, there is a
one-dimensional admissible datum $(\tau_0,V_{\tau_0})$. In this
case, all admissible data are in one-to-one correspondence with
irreducible representations of the group of connected components
of $K_{\mathbb {C}}^{\lambda_{\th}}$; the correspondence is
obtained by tensoring with $\tau_0.$ If $G$ is connected and
simply connected, then this component group is just the
fundamental group of the nilpotent orbit $ K_{\mathbb {C}}\cdot
\lambda_{\th}$.
\vspace{0.2in}\\
\noindent{\bf 7.\ Minimal Representations}
\setcounter{equation}{0}
\renewcommand{\theequation}{7.\arabic{equation}}
\vspace{0.1in}\\
\indent Let $G$ be a real reductive Lie group. Let $\pi$ be an
admissible representation of $G$. Let $\mathfrak {g} $ be the Lie
algebra of $G$. Three closely related invariants $WF(\pi)$,
$AS(\pi)$ and $Ass(\pi)$ in ${\mathfrak {g} ^\ast}$ which are
called the {\it wave\ front\ set\ } of $\pi$, the {\it asymptotic\
support\ of\ the\ character\ } of $\pi$ and the {\it associated\
variety\ } of $\pi$ respectively, are attached to a given
admissible representation $\pi$. The subsets $WF(\pi)$, $AS(\pi)$,
and $Ass(\pi)$ are contained in the cone ${\mathcal N }^{\ast}$
consisting of nilpotent elements in ${\mathfrak {g} ^\ast}$. They
are all invariant under the coadjoint action of $G$. Each of them
is a closed subvariety of ${\mathfrak {g} ^\ast}$, and is the
union of finitely many nilpotent orbits. It is known that
$WF(\pi)=AS(\pi)$. W. Schmid and K. Vilonen \cite{ScV2} proved that
$Ass(\pi)$ coincides with $WF(\pi)$ via the Kostant-Sekuguchi
correspondence. The dimensions of all three invariants are the
same, and is always even. We define the {\it Gelfand}-{\it
Kirillov\ dimension } of $\pi$ by
\begin{equation}
{\text{dim}_{G-K}}\,\pi := {1 \over 2 }\, \text{dim}\, WF(\pi)
\end{equation}
\noindent If ${\mathfrak {g} _\mathbb  {C}}$ is simple, there
exist a unique nonzero nilpotent ${G_\mathbb {C}}$-orbit
${{\mathcal {O} }_{min}} \subset \mathfrak {g} ^*_{\mathbb {C}}$
of minimal dimension, which is contained in the closure of any
nonzero nilpotent ${G_\mathbb {C}}$-orbit. In this case, we have
${\mathcal {O} }_{min}={{\mathcal {O} }_{X_\alpha}}$, where
${{\mathcal {O} }_{X_\alpha}}$ is the ${G_\mathbb {C}}$-orbit of a
nonzero highest root vector ${X_\alpha}$.

 \smallskip

\indent A nilpotent $G$-orbit ${\mathcal{O}} \subset {\mathfrak
{g} ^*}$ is said to be {\it minimal} if
\begin{equation}
{\text{dim}_{\mathbb{R}}} {\mathcal {O} } = {\text{dim}_\mathbb
{C}} {{\mathcal {O} }_{min}},
 \end{equation}
 \noindent equivalently, ${\mathcal {O} }$
is nonzero and contained in ${{\mathcal {O} }_{min}} \cap
{\mathfrak {g} ^\ast}$. An irreducible unitary representation
$\pi$ of $G$ is called {\it minimal} if
\begin{equation}
{\text{dim}_{G-K}}\pi = {1 \over 2} {\text{dim}_\mathbb
{C}}{{\mathcal {O} }_{min}}.
\end{equation}
\vspace{0.1in}\\
\noindent {\bf Remark\ 7.1.}\quad (1) If $G$ is not of type
${A_n}$, there are at most finity many minimal representations.
These are the unipotent representations attached to the minimal
orbit ${{\mathcal {O} }_{min}}$.

\smallskip

\indent (2) In many cases, the minimal representations are
isolated in the unitary dual $\hat{G}_u$ of $G$.\\
\indent (3) A minimal representation $\pi$ is almost always {\it
automorphic}, namely, $\pi$ occurs in ${L^2}({\Gamma\backslash
G})$ for some lattice $\Gamma$ in $G$. The theory of minimal
representation is the basis for the construction of large families
of other interesting automorphic representations. For example, it
is known that the end of complementary series of $Sp(n,1)$ and
${F_{4,1}}$ are both automorphic.
\vspace{0.1in}\\
\noindent {\bf Remark\ 7.2.}\quad Let $(\pi, V)$ be an irreducible
admissible representation of $\pi$. We denote by ${V^K}$ the space
of $K$-finite vectors for $\pi$, where $K$ is a maximal compact
subgroup of $G$. Then ${V^K}$ is a $U(\mathfrak {g} _{\mathbb
{C}})$-module. We fix any vector $0\neq x_{\pi}\in V^K$. Let
$U_n(\mathfrak {g} _{\mathbb {C}})$ be the subspace of
$U(\mathfrak {g} _{\mathbb  {C}})$ spanned by products of at most
$n$ elements of ${\mathfrak {g} _\mathbb  {C}}$. Put
$$X_n (\pi) := U_n (\mathfrak {g} _{\mathbb {C}}) {x_\pi}.$$ D. Vogan \cite{Vo1} proved
that $\text{dim}\, {X_n}(\pi)$ is asymptotic to ${{c(\pi)} \over
{d!}} \cdot {n^d}$ as $n \rightarrow \infty$. Here $c(\pi)$ and
$d$ are positive integers independent of the choice of ${x_\pi}$.
In fact, $d$ is the Gelfand-Kirillov dimension of $\pi$. We may
say that $d={\text{dim}_{G-K}}\pi$ is a good measurement of the
size of $\pi$.

\smallskip

\indent Suppose $F=\mathbb{R}$ or $\mathbb {C}$. Let $G$ be a
connected simple Lie group over $F$ and $K$ a maximal compact
subgroup of $G$. Let $\mathfrak {g} =\mathfrak{k}+\mathfrak {p}$
be the corresponding Cartan decomposition of $\mathfrak {g} $. If
$G/K$ is hermitian symmetric, all the minimal representations are
known to be either holomorphic or antiholmorphic. They can be
found in the list of unitary highest weight modules given in
\cite{EHW}. D. Vogan \cite{Vo2} proved the existence and unitarity
of the minimal representations for a family of split simple group
including $E_8$, $F_4$ and all classical groups except for the
${B_n}$-case $(n \geq 4)$, which no minimal representation seems
to exist. The construction of the minimal representations for
${G_2}$ was given by M. Duflo \cite{Du3} in the complex case and
by D. Vogan \cite{Vo6} in the real case. D. Kazhdan and G. Savin
\cite{K-S} constructed the spherical minimal representation for
every simple, split, simply laced group. B. Gross and N. Wallach
\cite{G-W} constructed minimal representations of all exceptional
groups of real rank 4. R. Brylinski and B. Kostant
\cite{BrK1}-\cite{BrK2} gave a construction of minimal
representations for any simple real Lie group $G$ under the
assumption that $G/K$ is not hermitian symmetric and minimal
representations exist.

\smallskip

\indent For a complex group $G$ not of type ${A_n}$, the
Harish-Chandra module of the spherical minimal representation can
be realized on $U(\mathfrak {g} )/J$, where $U(\mathfrak {g} )$ is
the universal enveloping algebra of $\mathfrak {g} $ and $J
\subset U(\mathfrak {g} )$ is the Joseph ideal of $\mathfrak {g} $
\cite{J}. It is known that $Sp(2n,\mathbb {C})$ has a
non-spherical minimal representation, that is, the odd piece of
the Weil representation. The following natural question arises.

\medskip

\noindent {\bf Question.}\quad Are there non-spherical minimal
representations for complex Lie groups other than $Sp(2n,\mathbb
{C})$\,?

\smallskip

\indent Recently P. Torasso \cite{To2} gave a uniform construction of
minimal representations for a simple group over any local field of
characteristic 0 with split rank $\geq 3$. He constructs a minimal
representation for each set of admissible datum associated to the
minimal orbit defined over $F$.\\

\smallskip

\indent Let ${\pi_{min}}$ be a minimal representation of $G$. It
is known that the annihilator of the Harish-Chandra module of
${\pi_{min}}$ in $U(\mathfrak {g} )$ is the Joseph ideal. D. Vogan
\cite{Vo2} proved that the restriction of ${\pi_{min}}$ to $K$  is
given by
$${\pi_{min}} {\vert_K}=\oplus^{\infty}_{n=0} V({\mu_0}+n{\beta}),$$
where $\beta$ is a highest weight for the action of $K$ on
$\mathfrak {p}$, ${\mu_0}$ is a fixed highest weight depending on
${\pi_{min}}$ and $V({\mu_0}+n\beta)$ denotes the highest weight
module with highest weight ${\mu_0}+n{\beta}$. Indeed, there are
two or one possibilities for $\beta$ depending on whether $G/K$ is
hermitian symmetric or not.
\vspace{0.1in}\\
\noindent {\bf Definition\ 7.3.}\quad (1) A {\it reductive dual
pair} in a reductive Lie group $G$ is a pair $(A,B)$ of closed
subgroups of $G$, which are both reductive and are centralizers of
each other.

\smallskip

\noindent (2) A reductive dual pair $(A,B)$ in $G$ is said to be
{\it compact\ } if at least one of $A$ and $B$ is compact.
\vspace{0.1in}\\
\noindent {\bf Duality  Conjecture.}\quad Let $(A,B)$ be a
reductive dual pair in a reductive Lie group $G$. Let
${\pi_{min}}$ be a minimal representation of $G$. Can you find a
Howe type correspondence between suitable subsets of the
admissible duals of $A$ and $B$ by restricting ${\pi_{min}}$ to $A
\times B$\,?

\smallskip

\indent In the 1970s R. Howe\,\cite{Ho3} first formulated the duality
conjecture for the Weil representation (which is a minimal
representation) of the symplectic group $Sp(n,F)$ over any local
field $F$. He \cite{Ho4} proved the duality conjecture for $Sp(n,F)$
when $F$ is archimedean and J. L. Waldspurger \cite{W} proved the
conjecture when $F$ is non-archimedean with odd residue
characteristic.
\vspace{0.1in}\\
\noindent {\bf Example\ 7.4.}\quad Let $G$ be the simply connected
split real group ${E_8}$. Then $K=Spin(16)$ is a maximal compact
subgroup of $G$. We take $A=B=Spin(8)$. It is easy to see that the
pair $(A, B)$ is not only a reductive dual pair in $Spin(16)$, but
also a reductive dual pair in $G$. Let ${\pi_{min}}$ be a minimal
representation of $G$. J.-S. Li \cite{Li1} showed that the
restriction of ${\pi_{min}}$ to $A \times B \subset G$ is
decomposed as follows:
$${\pi_{min}}{\vert_{A \times B}} = {\oplus_\pi}m(\pi)\cdot (\pi
\otimes \pi),$$ where $\pi$ runs over all irreducible
representations of $Spin(8)$ and $m(\pi)$ is the multiplicity with
which $\pi\otimes\pi$ occurs. It turns out that $m(\pi)=+\infty$
for all $\pi$.
\vspace{0.1in}\\

\noindent{\bf Example\ 7.5.}\quad Let $G$ be the simply connected
quaternionic ${E_8}$ with split rank 4. We let $A=Spin(8)$ and
$B=Spin(4,4)$. Then the pair $(A,B)$ is a reductive dual pair in
$G$. Let ${\pi_{min}}$ be the minimal representation of $G$. H. Y.
Loke\,\cite{L} proved that the restriction of ${\pi_{min}}$ to $A
\times B$ is decomposed as follows.
$$\pi_{min}\vert_{A \times B} =\oplus_{\pi} m(\pi)\cdot (\pi
\otimes \pi'),$$ where $\pi$ runs over all irreducible
representations of $A$ and $\pi'$ is the discrete series
representation of $B$ which is uniquely determined by $\pi$. All
the multiplicities are {\it finite}. But they are {\it unbounded}.
\pagestyle{myheadings}
 \markboth{\headd Jae-Hyun Yang $~~~~~~~~~~~~~~~~~~~~~~~~~~~~~~~~~~~~~~~~~~~~~~~~~~~~~~~~~~$}
 {\headd $~~~~~~~~~~~~~~~~~~~~~~~~~~~~~~~~~~~$The Method of Orbits for Real Lie Groups}

\indent The following interesting problem is proposed by Li.
\vspace{0.1in}\\
\noindent {\bf Problem\ 7.6.} \quad Let $(A,B)$ be a {\it compact
} reductive dual pair in $G$. Describe the explicit decomposition
of the restriction ${\pi_{min}} {\vert_{A \times B}}$ of a minimal
representation ${\pi_{min}}$ of $G$ to $A \times B$.
\vspace{0.1in}\\
\noindent {\bf Remark\ 7.7.}\quad In \cite{HPS}, J. Huang, P.
Paudzic and G. Savin dealt with the family of dual pairs $(A, B)$,
where $A$ is the split exceptional group of type ${G_2}$ and $B$
is compact. \pagestyle{myheadings}
 \markboth{\headd Jae-Hyun Yang $~~~~~~~~~~~~~~~~~~~~~~~~~~~~~~~~~~~~~~~~~~~~~~~~~~~~~~~~~~$}
 {\headd $~~~~~~~~~~~~~~~~~~~~~~~~~~~~~~~~~~~$The Method of Orbits for Real Lie Groups}
\vspace{0.1in}\\
\noindent{\bf 8.\ The Heisenberg Group $H_{\mathbb{R}}^{(g,h)}$ }
\setcounter{equation}{0}
\renewcommand{\theequation}{8.\arabic{equation}}
\def\gh{\Cal G^J}
\def\de{\delta}
\def\lambdamk{(\lambdaambda,\mu,\kappa)}
\def\Om{\Omega}
\def\bA{\bold A}
\def\bH{\bold H}
\def\J{J\in {\Bbb Z}^{(h,g)}_{\geq 0}}
\def\N{N\in {\Bbb Z}^{(h,g)}}
\def\Dm{\lambdaeft[\matrix -A\\ -B\endmatrix\right]}
\def\dt{{{d}\over {dt}}\bigg|_{t=0}}
\def\lambdat{\lambdaim_{t\to 0}}
\def\zhg{\BZ^{(h,g)}}
\def\bhg{\BR^{(h,g)}}
\def\ex{\par\smallpagebreak\noindent}
\def\pis{\pi i \sigma}
\def\sd{\,\,{\vartriangleright}\kern -1.0ex{<}\,}
\def\sc{\bf}
\def\wt{\widetilde}
\vspace{0.1in}\\
\indent For any positive integers $g$ and $h$, we consider the
Heisenberg group
$$H_{\mathbb{R}}^{(g,h)}=\left\{\,(\lambda,\mu,\kappa)\,\vert\ \lambda,\mu\in
\mathbb{R}^{(h,g)},\ \kappa\in \mathbb{R}^{(h,h)},\
\kappa+\mu\,^t\!\lambda\ \text{symmetric}\ \right\}$$ with the
multiplication law
$$(\lambda,\mu,\kappa)\circ
(\lambda',\mu',\kappa')=(\lambda+\lambda',\mu+\mu',\kappa+
\kappa'+\lambda\,^t\!\mu'-\mu\,^t\!\lambda'). $$ Here
$\mathbb{R}^{(h,g)}\,(\,\text{resp.}\ \mathbb{R}^{(h,h)})$ denotes
the all $h\times g\,(\,\text{resp.}\ h\times h)$ real matrices.

\smallskip

\indent The Heisenberg group $H_{\mathbb{R}}^{(g,h)}$ is embedded
to the symplectic group $Sp(g+h,\mathbb{R})$ via the mapping
$$H_{\mathbb{R}}^{(g,h)}\ni (\lambda,\mu,\kappa)\longmapsto \begin{pmatrix} E_g & 0
& 0 & ^t\mu \\ \lambda & E_h & \mu & \kappa \\ 0 & 0 & E_g & -^t\lambda \\
0 & 0 & 0 & E_h \end{pmatrix} \in Sp(g+h,\mathbb{R}).$$ This
Heisenberg group is a 2-step nilpotent Lie group and is important
in the study of toroidal compactifications of Siegel moduli
spaces. In fact, $H_{\mathbb{R}}^{(g,h)}$ is obtained as the
unipotent radical of the parabolic subgroup of
$Sp(g+h,\mathbb{R})$ associated with the rational boundary
component $F_g$\,(\,cf.\,\cite{F-C}\,p.\,123 or \cite{N}
p.\,21\,). For the motivation of the study of this Heisenberg
group we refer to \cite{Y7}-\cite{Y11} and \cite{Zi}. We refer to
\cite{Y2}-\cite{Y6} for more results on $H_{\mathbb{R}}^{(g,h)}$.
\pagestyle{myheadings}
 \markboth{\headd Jae-Hyun Yang $~~~~~~~~~~~~~~~~~~~~~~~~~~~~~~~~~~~~~~~~~~~~~~~~~~~~~~~~~~$}
 {\headd $~~~~~~~~~~~~~~~~~~~~~~~~~~~~~~~~~~~$The Method of Orbits for Real Lie Groups}

\indent In this section, we describe the Schr{\"o}dinger
representations of $H_{\mathbb{R}}^{(g,h)}$ and the coadjoint
orbits of $H_{\mathbb{R}}^{(g,h)}$. The results in this section
are based on the article \cite{YS} with some corrections.
\vspace{0.2in}\\
\noindent{\bf 8.1. Schr\"{o}dinger Representations}
\vspace{0.1in}\\
\indent First of all, we observe that $H_{\mathbb{R}}^{(g,h)}$ is
a 2-step nilpotent Lie group. It is easy to see that the inverse
of an element $(\lambda,\mu,\kappa)\in H_{\mathbb {R}}^{(g,h)}$ is
given by
$$(\lambda,\mu,\kappa)^{-1}=(-\lambda,-\mu,-\kappa+\lambda\,^t\!\mu-\mu\,^t\!\lambda).$$
Now we set
\begin{equation}
[\lambda,\mu,\kappa]:=(0,\mu,\kappa)\circ
(\lambda,0,0)=(\lambda,\mu,\kappa-\mu\,^t\! \lambda).
\end{equation}
\noindent Then $H_{\mathbb {R}}^{(g,h)}$ may be regarded as a
group equipped with the following multiplication
\begin{equation}
[\lambda,\mu,\kappa]\diamond
[\lambda_0,\mu_0,\kappa_0]=[\lambda+\lambda_0,\mu+\mu_0,
\kappa+\kappa_0+\lambda\,^t\!\mu_0+\mu_0\,^t\!\lambda].
\end{equation}
\noindent The inverse of $[\lambda,\mu,\kappa]\in H_{\mathbb
{R}}^{(g,h)}$ is given by
$$[\lambda,\mu,\kappa]^{-1}=[-\lambda,-\mu,-\kappa+\lambda\,^t\!\mu+\mu\,^t\!\lambda].$$
We set
\begin{equation}
 K=\left\{\,[0,\mu,\kappa]\in H_{\mathbb
{R}}^{(g,h)}\,\Big| \, \mu\in \mathbb {R}^{(h,g)},\
\kappa=\,^t\!\kappa\in \mathbb {R}^{(h,h)}\ \right\}.
 \end{equation}
\noindent Then $K$ is a commutative normal subgroup of $H_{\mathbb
{R}}^{(g,h)}$. Let ${\hat {K}}$ be the Pontrajagin dual of $K$,
i.e., the commutative group consisting of all unitary characters
of $K$. Then ${\hat {K}}$ is isomorphic to the additive group
$\mathbb {R}^{(h,g)}\times \text{Symm}\,(h,\mathbb {R})$ via
\begin{equation}
<a,{\hat a}>=e^{2 \pi i\sigma({\hat {\mu}}\,^t\!\mu+{\hat
{\kappa}}\kappa)}, \ \ \ a=[0,\mu,\kappa]\in K,\ {\hat {a}}=({\hat
{\mu}},{\hat {\kappa}})\in \hat{K}.
\end{equation}
\noindent We put
\begin{equation}
S=\left\{\,[\lambda,0,0]\in H_{\mathbb {R}}^{(g,h)}\,\Big|\
\lambda\in \mathbb {R}^{(h,g)}\, \right\}\cong \mathbb
{R}^{(h,g)}.
\end{equation}
\noindent Then $S$ acts on $K$ as follows:
\begin{equation}
\alpha_{\lambda}([0,\mu,\kappa])=[0,\mu,\kappa+\lambda\,^t\!\mu+\mu\,^t\!\lambda],
\ \ \ [\lambda,0,0]\in S.
\end{equation}
\noindent It is easy to see that the Heisenberg group $\left(
H_{\mathbb {R}}^{(g,h)}, \diamond\right)$ is isomorphic to the
semi-direct product $S\ltimes K$ of $S$ and $K$ whose
multiplication is given by
$$(\lambda,a)\cdot
(\lambda_0,a_0)=(\lambda+\lambda_0,a+\alpha_{\lambda}(a_0)),\ \
\lambda,\lambda_0\in S,\ a,a_0\in K.$$ On the other hand, $S$ acts
on ${\hat {K}}$ by
\begin{equation}
\alpha_{\lambda}^{*}({\hat {a}})=({\hat {\mu}}+2{\hat
{\kappa}}\lambda, {\hat {\kappa}}),\ \ [\lambda,0,0]\in S,\ \
a=({\hat {\mu}},{\hat {\kappa}})\in {\hat {K}}.
\end{equation}
\noindent Then, we have the relation $<\alpha_{\lambda}(a),{\hat
{a}}>=<a,\alpha_{\lambda}^{*} ({\hat {a}})>$ for all $a\in K$ and
${\hat {a}}\in {\hat {K}}.$

\smallskip

 \indent We have three types of $S$-orbits in ${\hat
{K}}.$

\smallskip

\noindent {\scshape Type I.} Let ${\hat{\kappa}} \in
\text{Sym}\,(h,\mathbb {R})$ be nondegenerate. The $S$-orbit of
${\hat {a}}({\hat {\kappa}}):=(0,{\hat {\kappa}}) \in {\hat {K}}$
is given by

\begin{equation}
\hat{\mathcal{O}}_{\hat{\kappa}}= \left\{(2{\hat {\kappa}}\lambda,{\hat
{\kappa}}) \in {\hat {K}}\ \Big|\ \lambda \in \mathbb{R}^{(h,g)}
\right\} \cong \mathbb {R}^{(h,g)}.
\end{equation}
\noindent {\scshape Type II.}\ \ Let $({\hat{\mu}},{\hat{\kappa}})\in
\mathbb{R}^{(h,g)}\times\text{Sym}\,(h,\mathbb {R})$ with
degenerate  ${\hat{\kappa}}\neq 0.$ Then
\begin{equation}
\hat{\mathcal{O}}_{(\hat{\mu}, \hat{\kappa})} = \left\{ \hat  \mu +
2\hat{\kappa}\lambda, \hat{\kappa}) \Big|\ \lambda \in \mathbb
{R}^{(h,g)} \right\} \subsetneqq \mathbb {R}^{(h,g)}\times \{
\hat{\kappa} \} .
\end{equation}
\noindent {\scshape Type III.} Let $\hat{y} \in \mathbb{R}^{(h,g)}$. The
$S$-orbit ${\hat {\mathcal {O} }}_{\hat {y}}$ of $\hat{a}(\hat{y})
:= (\hat{y} ,0)$ is given by
\begin{equation}
{\hat {\mathcal {O} }}_{\hat {y}}=\left\{\,({\hat
{y}},0)\,\right\}={\hat {a}} ({\hat {y}}).
\end{equation}
\noindent We have
$$
{\hat{K}}= \left( \bigcup_{\begin{subarray}{c} \hat{\kappa} \in
\text{Sym} (h,\mathbb{R}) \\ {\hat{\kappa}} \,\ \text{nondegenerate} \end{subarray}}
 \hat {\mathcal{O}}_ {\hat\kappa }
 \right)
  \bigcup
\left( \bigcup_{{\hat {y}}\in \mathbb{R}^{(h,g)}}{\hat {\mathcal {O} }}_{\hat
{y}}\right)  \bigcup  \left( \bigcup_{\begin{subarray}{c}({\hat{\mu}},{\hat
{\kappa}}) \in \mathbb{R}^{(h,g)} \times \text{Sym}(h,\mathbb{R}) \\ \hat{\kappa} \neq 0 \,\ \text{degenerate} \end{subarray}}
 {\hat {\mathcal {O}}}
_({\hat \mu},{\hat \kappa }) \right)
$$
\noindent as a set. The
stabilizer $S_{\hat {\kappa}}$ of $S$ at ${\hat {a}}({\hat
{\kappa}})=(0,{\hat {\kappa}})$ is given by
\begin{equation}
S_{\hat {\kappa}}=\{0\}.
 \end{equation}
\noindent And the stabilizer $S_{\hat {y}}$ of $S$ at ${\hat
{a}}({\hat {y}})= ({\hat {y}},0)$ is given by
\begin{equation}
S_{\hat {y}}=\left\{\,[\lambda,0,0]\,\Big|\ \lambda\in \mathbb
{R}^{(h,g)}\,\right\}=S \,\cong\,\mathbb {R}^{(h,g)}.
\end{equation}
\indent From now on, we set $G=H_{\mathbb {R}}^{(g,h)}$ for
brevity. It is known that $K$ is a closed, commutative normal
subgroup of $G$. Since $(\lambda,\mu,\kappa)=(0,\mu,
\kappa+\mu\,^t\!\lambda)\circ (\lambda,0,0)$ for
$(\lambda,\mu,\kappa)\in G,$ the homogeneous space $X=K\backslash G$
can be identified with $\mathbb{R}^{(h,g)}$ via
$$
Kg=K\circ (\lambda,0,0)\longmapsto \lambda,\ \ \ g=(\lambda,\mu,\kappa)\in G.
$$
We observe that $G$ acts on $X$ by
\begin{equation}
(Kg)\cdot g_0=K\,(\lambda+\lambda_0,0,0)=\lambda+\lambda_0,
\end{equation}
\noindent where $g=(\lambda,\mu,\kappa)\in G$ and
$g_0=(\lambda_0,\mu_0,\kappa_0)\in G.$

\medskip

\indent If $g=(\lambda,\mu,\kappa)\in G$, we have
\begin{equation}
k_g=(0, \mu, \kappa+\mu\,^t\!\lambda),\ \ \ s_g=(\lambda,0,0)
\end{equation}
\noindent in the Mackey decomposition of $g=k_g \circ s_g$
(cf.\cite{M}). Thus if $g_0=(\lambda_0,\mu_0,\kappa_0)\in G,$ then we
have
\begin{equation}
s_g\circ g_0=(\lambda,0,0)\circ
(\lambda_0,\mu_0,\kappa_0)=(\lambda+\lambda_0,\mu_0,
\kappa_0+\lambda\,^t\!\mu_0)
\end{equation}
\noindent and so
\begin{equation}
k_{s_g\circ
g_0}=(0,\mu_0,\kappa_0+\mu_0\,^t\!\lambda_0+\lambda\,^t\!\mu_0+
\mu_0\,^t\!\lambda).
\end{equation}
 \indent For a real symmetric matrix
$c=\,^tc\in \mathbb{R}^{(h,h)}$ with $c\neq 0$, we consider the
one-dimensional unitary representation $\sigma_c$ of $K$ defined
by
\begin{equation}
\sigma_c\left((0,\mu,\kappa)\right)=e^{2\pi
i \sigma(c\kappa)}\,I,\ \ \ (0,\mu,\kappa)\in K,
\end{equation}
\noindent where $I$ denotes the identity mapping. Then the induced representation
$U(\sigma_c):=\text{Ind}_K^G\,\sigma_c$ of $G$ induced from $\sigma_c$ is
realized in the Hilbert space ${\mathcal H}_{\sigma_c}=L^2(X,d{\dot
{g}},\mathbb {C}) \cong L^2\left(\mathbb{R}^{(h,g)}, d\xi\right)$ as
follows. If $g_0=(\lambda_0,\mu_0,\kappa_0)\in G$ and $x=Kg\in X$
with $g=(\lambda,\mu,\kappa)\in G,$ we have
\begin{equation}
\left(U_{g_0}(\sigma_c)f\right)(x)=\sigma_c\left(k_{s_g\circ
g_0}\right)\left( f(xg_0)\right),\ \ f\in {\mathcal
H}_{\sigma_c}.
\end{equation}
\noindent It follows from (8.16) that
\begin{equation}
\left(U_{g_0}(\sigma_c)f\right)(\lambda)=e^{2\pi
i\sigma\{c(\kappa_0+\mu_0\,^t\!\lambda_0+
2\lambda\,^t\!\mu_0)\}}\,f(\lambda+\lambda_0).
\end{equation}
\noindent Here, we identified $x=Kg$\,(resp.\,$xg_0=Kg g_0$) with $\lambda$\,(resp.\,
$\lambda+\lambda_0$). The induced representation $U(\sigma_c)$ is
called the {\it Schr{\" {o}}dinger\ representation} of $G$
associated with $\sigma_c.$ Thus $U(\sigma_c)$ is a monomial
representation.

\medskip

\indent Now, we denote by ${\,mathcal
H}^{\sigma_c}$ the Hilbert space consisting of all functions
$\phi:G\longrightarrow \mathbb {C}$ which satisfy the
following conditions:

\smallskip

\indent (1) $\phi(g)$ is measurable measurable with respect to
$dg$,\\
\indent (2) $\phi\left((0,\mu,\kappa)\circ g)\right)= e^{2\pi i\sigma(c\kappa)}\phi(g)$ \ \
for\ all $g\in G,$\\
\indent (3) $\parallel\phi\parallel^2:=\int_X\,\vert \phi(g)\vert^2\, d{\dot
{g}} < \infty,\ \ \ {\dot {g}}=Kg,$\\
\noindent where $dg$\,(resp.\,$d{\dot {g}}$) is a $G$-invariant measure on
$G$ (resp.\,$X=K\backslash G$). The inner product $(\,,\,)$ on
${\mathcal H}^{\sigma_c}$ is given by
$$
(\phi_1,\phi_2)=\int_G\,\phi_1(g)\,{\overline {\phi_2(g)}}\,dg\
\ \ \ \text{for} \ \phi_1,\, \phi_2\in {\mathcal H}^{\sigma_c}.
$$
We observe that the mapping $\Phi_c:{\mathcal
H}_{\sigma_c}\longrightarrow {\mathcal H}^{\sigma_c}$ defined by
\begin{equation}
\left( \Phi_c(f)\right)(g)=e^{2\pi i\sigma
\{c(\kappa+\mu\,^t\!\lambda)\}}\,f(\lambda),\ \ \ f\in {\mathcal
H}_{\sigma_c},\ g=(\lambda,\mu, \kappa)\in G
\end{equation}
\noindent is an isomorphism of Hilbert spaces. The inverse $\Psi_c:{\mathcal H}^{\sigma_c}
\longrightarrow {\mathcal H}_{\sigma_c}$ of $\Phi_c$ is given by
\begin{equation}
\left( \Psi_c(\phi)\right)(\lambda)=\phi((\lambda,0,0)),\ \ \ \phi\in
{\mathcal H}^{\sigma_c},\ \lambda\in \mathbb{R}^{(h,g)}.
\end{equation}
\noindent The Schr{\"{o}}dinger representation $U(\sigma_c)$ of $G$ on ${\mathcal H}^{\sigma_c}$ is
given by
\begin{equation}
\left( U_{g_0}(\sigma_c)\phi\right)(g)= e^{2\pi i\sigma\{
c(\kappa_0+\mu_0\,^t\!\lambda_0+\lambda\,^t\!\mu_0-\lambda_0\,^t\!\mu)\}}\,
\phi\left( (\lambda_0,0,0)\circ g\right),
\end{equation}
\noindent where $g_0=(\lambda_0,\mu_0,\kappa_0),\ g=(\lambda,\mu,\kappa)\in G$ and
$\phi\in {\mathcal H}^{\sigma_c}.$ (8.22) can be expressed as follows.
\begin{equation}
\left(\,U_{g_0}(\sigma_c)\phi\,\right)(g)=e^{2\pi
i\sigma\{c(\kappa_0+\kappa+
\mu_0\,^t\!\lambda_0+\mu\,^t\!\lambda+2\lambda\,^t\mu_0)\}}\,\phi((\lambda_0+\lambda,0,0)).
\end{equation}
\vspace{0.1in}\\
\noindent {\bf Theorem\ 8.1.}\quad {\it Let $c$ be a positive
symmetric half-integral matrix of degree $h$. Then the Schr{\"
{o}}dinger representation $U(\sigma_c)$ of $G$ is irreducible.}
\vspace{0.05in}\\
\noindent{\it Proof.}\quad The proof
can be found in \cite{Y3},\ Theorem 3.
\vspace{0.2in}\\
\noindent{\bf 8.2.  The Coadjoint Orbits of Picture}
\vspace{0.1in}\\
\indent In this subsection, we find the coadjoint orbits of the
Heisenberg group $H^{(g,h)}_\mathbb {R}$ and describe the
connection between the coadjoint orbits and the unitary dual of
$H^{(g,h)}_\mathbb{R}$ explicitly.\\
\indent For brevity, we let $G := H^{(g,h)}_\mathbb {R}$ as
before. Let $\frak g$ be the Lie algebra of $G$ and let ${\frak
g}^*$ be the dual space of $\frak g.$ We observe that $\frak g$
can be regarded as the subalgebra consisting of all $(g+h)\times
(g+h)$ real matrices of the form
$$X(\alpha,\beta,\gamma)= \begin{pmatrix} 0&0&0&{}^t\!\beta\\
                                    \alpha&0&\beta&\gamma\\
                    0&0&0&-{}^t\!\alpha\\
            0&0&0&0 \end{pmatrix}, \;\alpha,\beta \in \mathbb{R}^{(h,g)},\;
                    \gamma={}^t\!\gamma \in \mathbb {R}^{(h,h)} $$
of the lie algebra  $\mathfrak {sp} (g+h,\mathbb {R})$ of the
symplectic group $Sp(g+h,\mathbb {R}).$ An easy computation yields
$$[X(\alpha,\beta,\gamma), \,X(\delta,\epsilon,\xi)] =
X(0,0,\alpha^t\epsilon+
\epsilon^t\alpha-\beta^t\delta-\delta^t\beta).$$ The dual space
${\frak g}^*$ of $\frak g$ can be identified with the vector space
consisting of all $(g+h)\times (g+h)$ real matrices of the form $$
F(a,b,c) = \begin{pmatrix} 0& {}^t\!a&0&0\\
                  0&0&0&0\\
          0&{}^t\!b&0& 0\\
          b&c&-a&0\end{pmatrix},\; a,b \in \mathbb{R}^{(h,g)},\;c ={}^t\!c \in
          \mathbb{R}^{(h,h)} $$
so that
\begin{equation}
<F(a,b,c), X(\alpha,\beta,\gamma)> =
\sigma(F(a,b,c)\,X(\alpha,\beta,\gamma)) =2 \sigma({}^t\alpha a
+{}^t\!b\beta)+\sigma(c\gamma).
\end{equation}

\smallskip

\noindent The adjoint representation $Ad$ of $G$ is given by $Ad_G(g)X = g X g^{-1}$
for $g\in G$ and $X\in \frak g.$ For $g \in G$ and $F\in {\frak
g}^*, \; gFg^{-1}$ is not of the form $F(a,b,c).$ We denote by
$(gFg^{-1})_*$ the $$\begin{pmatrix} 0&*&0&0 \\
           0&0&0&0 \\
       0&*&0&0\\
       *&*&*&0 \end{pmatrix}-\text{part}$$
of the matrix $gFg^{-1}.$ Then it is easy to see that the
coadjoint representation $Ad_G^* : G \longrightarrow GL({\frak
g}^*)$ is given by $Ad_G^*(g)F = (gFg^{-1})_*,$ where $g \in G$
and $F\in {\frak g}^*.$ More precisely,
\begin{equation}
Ad_G^*(g)F(a,b,c) =
F(a+c\mu, b-c\lambda,c),
\end{equation}
\noindent where $g=(\lambda,\mu,\kappa)\in G.$ Thus the coadjoint orbit
$\Omega_{a,b}$ of $G$ at $F(a,b,0) \in {\frak g}^*$ is given by
\begin{equation}
\Omega_{a,b} = Ad_G^*(G)\,F(a,b,0)=\{F(a,b,0)\},\;\text{a single
point}
\end{equation}
\noindent and the coadjoint orbit $\Omega_c$ of $G$ at
$F(0,0,c)\in {\frak g}^*$ with $c \ne 0$ is given by
\begin{equation}
\Omega_c =
Ad_G^*(G)\,F(0,0,c) = \{ F(a,b,c) \vert a,b \in \mathbb{R}^{(h,g)}\}\cong
 \mathbb {R}^{(h,g)}\times \mathbb {R}^{(h,g)}.
\end{equation}
\noindent Therefore the coadjoint orbits of $G$ in $\mathfrak{g}^{*}$ fall into two classes :\\
\indent ({\bf I}) \ \ The single point $\{ \Omega_{a,b} \vert a,b \in
\mathbb {R}^{(h,g)} \}$
 located in the plane $c=0.$ \\
\indent ({\bf II}) \ \ The affine planes $\{ \Omega_c \vert c={}^t\!c \in
\mathbb{R}^{(h,h)}, \;c\ne o \}$ parallel to the homogeneous plane $c=0.$\\
\noindent In other words, the orbit space $\mathcal {O} (G)$ of
coadjoint orbits is parametrized by
$$\begin{cases}
 c-\text{axis},\;c\ne
0, \; c={}^t\!c \in \mathbb{R}^{(h,h)} ; \\ (a,b)-\text{plane}
\approx \mathbb {R}^{(h,g)} \times \mathbb{R}^{(h,g)}.
\end{cases}$$
\noindent The single point coadjoint orbits of the type
$\Omega_{a,b}$ are said to be the {\it degenerate} orbits of $G$
in $\frak g^*.$ On the other hand, the flat coadjoint orbits of
the type $\Omega_c$ are said to be the {\it non-degenerate} orbits
of $G$ in $\frak g^*.$ Since $G$ is connected and simply connected
2-step nilpotent Lie group, according to A. Kirillov (cf.
\cite{Ki1} or \cite{Ki2} p.249, Theorem 1), the unitary dual
$\widehat{G}$ of $G$ is given by
\begin{equation}
\widehat{G} = \left( \mathbb {R}^{(h,g)}\times \mathbb {R}^{(h,g)} \right)
\coprod \left\{ z \in \mathbb {R}^{(h,h)} \;\vert\;  z={}^t\!z ,\;
z\ne 0 \right\},
\end{equation}
\noindent where $\coprod$ denotes the disjoint union. The topology
of $\widehat{G}$ may be described as follows. The topology on
$\{c-\text{axis}-(0)\}$ is the usual topology of the Euclidean
space and the topology on $\{F(a,b,0) \vert a,b \in \mathbb
{R}^{(h,g)} \}$ is the usual Euclidean topology. But a sequence on
the $c$-axis which converges to $0$ in the usual topology
converges to the whole Euclidean space $\mathbb {R}^{(h,g)} \times
\mathbb {R}^{(h,g)}$ in the topology of $\hat{G}$. This is just
the quotient topology on ${\frak g}^*/G$ so that algebraically and
topologically $\widehat{G} = {\frak g}^* \slash G.$ \\
\indent It is well known that each coadjoint orbit is a symplectic manifold. We will
 state this fact in detail. For the present time being, we fix an element $F$
of $\frak g^*$ once and for all. We consider the alternating
$\mathbb{R}$-bilinear form ${\bold B}_F$ on $\frak g$ defined by
\begin{equation}
{\bold B}_F(X,Y) = <F,[X,Y]>=< ad_{\frak g}^*(Y)F,X>,\; X,Y \in {\frak
g},
\end{equation}
where $ad^*_{\frak g} : {\frak g}
\longrightarrow \text{End}({\frak g}^*)$ denotes the
differential of the coadjoint representation $Ad^*_G : G
\longrightarrow GL({\frak g}^*).$ More precisely, if
$F=F(a,b,c),\; X=X(\alpha,\beta,\gamma),\;\text{and}\;
Y=X(\delta,\epsilon,\xi),$ then
\begin{equation}
{\bold B}_F(X,Y) = \sigma \{
c(\alpha^t\epsilon + \epsilon^t\alpha -\beta^t\delta
-\delta^t\beta)\}.
\end{equation}
 For $F\in {\frak g}^*,$ we let $$G_F
 = \{ g  \subset G \vert Ad_G^*(g)F=F \} $$ be the stabilizer of the
coadjoint action $Ad^*$ of $G$ on ${\frak g}^*$ at $F.$ Since
$G_F$ is a closed subgroup of $G,$ $G_F$ is a Lie subgroup of $G.$
We denote by ${\frak g}_F$ the Lie subalgebra of $\frak g$
corresponding to $G_F.$ Then it is easy to show that
\begin{equation}
{\frak g}_F = rad\,{\bold B}_F = \{ X \in {\frak g} \vert ad^*_{\frak g}(X)F=
0\}.
 \end{equation}
 Here $rad\,{\bold B}_F$ denotes the radical of
${\bold B}_F$ in $\frak g.$ We let $\dot{{\bold B}_F}$ be the
non-degenerate alternating $\mathbb{R}$-bilinear form on the quotient
vector space $\frak g \slash rad\;{\bold B}_F$ induced from
${\bold B}_F.$ Since we may identify the tangent space of the
coadjoint orbit $\Omega_F \cong G \slash G_F $ with $\frak g
\slash {\frak g}_F = \frak g \slash rad{\bold B}_F,$ we see that
the tangent space of $\Omega_F$ at $F$ is a symplectic vector
space with respect to the symplectic form $\dot{{\bold B}_F}.$\\
\indent Now we are ready to prove that the coadjoint orbit $\Omega_F =
Ad^*_G(G)F$ is a symplectic manifold. We denote by $\widetilde X$
the smooth vector field on $\frak g^*$ associated to $X\in \frak
g.$ That means that for each $\ell \in {\frak g}^*,$ we have
\begin{equation}
{\widetilde X}(\ell) = ad^*_{\frak g}(X)\; \ell.
\end{equation}
 We define the differential 2-form $B_{\Omega_F}$ on $\Omega_F$ by
 \begin{equation}
B_{\Omega_F} (\widetilde X, \widetilde Y) = B_{\Omega_F}(
ad^*_{\frak g}(X)F, ad^*_{\frak g}(Y)F) :={\bold B}_F(X,Y),
\end{equation}
where $X,Y \in {\frak g}.$
\vspace{0.1in}\\
\noindent {\bf Lemma 8.2.}\quad {\it $B_{\Omega_F}$ is
non-degenerate.}
\vspace{0.05in}\\
\noindent{\it Proof.}\quad Let $\widetilde X$ be the smooth vector field on
$\frak g^*$  associated to $X \in \frak g$ such that $B_{\Omega_F}(\widetilde X, \widetilde
Y)= 0$ for all $\widetilde Y$ with $Y \in \frak g.$ Since
$B_{\Omega_F}(\widetilde X, \widetilde Y) = {\bold B}_F(X,Y) =0$
for all $Y \in {\frak g},\; X \in {\frak g}_F.$ Thus
$\widetilde{X} =0.$ Hence $B_{\Omega_F}$ is non-degenerate.
\vspace{0.1in}\\
\noindent {\bf Lemma 8.3.}\quad {\it $B_{\Omega_F}$ is closed.}
\vspace{0.05in}\\
\noindent{\it Proof.}\quad If $\widetilde{X_1},\,
\widetilde{X_2}, \text{and} \widetilde{X_3}$ are three smooth vector fields on
$\frak g^*$ associated to $X_1,\,X_2,\, X_3 \in \frak g,$ then
$$
\begin{array}{lllll}
dB_{\Omega_F}(\widetilde{X_1},\widetilde{X_2},\widetilde{X_3})
&=\widetilde{X_1}(B_{\Omega_F}(\widetilde{X_2},\widetilde{X_3}))-\widetilde{X_2}(B_{\Omega_F}(\widetilde{X_1},\widetilde{X_3}))+\widetilde{X_3}(B_{\Omega_F}(\widetilde{X_1},\widetilde{X_2}))\\
&-B_{\Omega_F}([\widetilde{X_1},\widetilde{X_2}],\widetilde{X_3})+B_{\Omega_F}([\widetilde{X_1},\widetilde{X_3}],\widetilde{X_2})-B_{\Omega_F}([\widetilde{X_2},\widetilde{X_3}],\widetilde{X_1})\\
&=-<F,[[X_1,X_2],X_3]+[[X_2,X_3],X_1]+[[X_3,X_1],X_2]>\\
&=0
\qquad (\text {by the Jacobi identity}).
\end{array}
$$
Therefore $B_{\Omega_F}$ is closed. \hfill$\square$

\smallskip

\indent In summary, $(\Omega_F, B_{\Omega_F})$ is a symplectic manifold of
dimension $2hg$ or $0.$

\smallskip

\indent In order to describe the
irreducible unitary representations of $G$ corresponding to the
coadjoint orbits under the Kirillov correspondence, we have to
determine the polarizations of $\frak g$ for the linear forms
$F\in {\frak g}^*.$

\smallskip

\noindent {\slshape Case} {\bf I.}\,\ $F=F(a,b,0);$ the degenerate
case.\par According to (8.26), $\Omega_F = \Omega_{a,b}
=\{F(a,b,0)\}$ is a single point. It follows from (8.30) that
${\bold B}_F(X,Y) =0$ for all $X,Y \in \frak g.$ Thus $\frak g$ is
the unique polarization of $\frak g$ for $F.$ The Kirillov
correspondence says that the irreducible unitary representation
$\pi_{a,b}$ of $G$ corresponding to the coadjoint orbit
$\Omega_{a,b}$ is given by
\begin{equation}
\pi_{a,b}( \text{exp}\,
X(\alpha,\beta,\gamma)) = e^{2\pi i <F,X(\alpha,\beta,\gamma)>} =
e^{4 \pi i \sigma({}^ta\alpha +{}^tb\beta)}.
\end{equation}
That is, $\pi_{a,b}$ is a one-dimensional degenerate representation of
$G.$

\smallskip

\noindent {\slshape Case} {\bf II.} \,\ $F=F(0,0,c),\; 0\ne c={}^t
c \in \mathbb{R}^{(h,h)}:$ the non-degenerate case. \\
\indent According to (8.27), $\Omega_F =\Omega_c =\{ F(a,b,c)
\vert a,b \subset \mathbb{R}^{(h,g)} \}.$ By (8.30), we see that
\begin{equation}
\frak k =
\{\;X(0,\beta,\gamma) \vert \beta \in \mathbb {R}^{(h,g)},\;
 \gamma={}^t\gamma \in \mathbb {R}^{(h,h)} \}
 \end{equation}
 is a polarization of $\frak g$ for $F,$ i.e.,$\frak k$ is a Lie subalgebra of
 $\frak g$ subordinate to $F\in \frak g^*$ which is maximal among the totally
 isotropic vector subspaces of $\frak g$ relative to the alternating $\mathbb{R}$-bilinear
 form ${\bold B}_F.$ Let $K$ be the simply connected Lie subgroup of $G$
 corresponding to the Lie subalgebra $\frak k$ of $\frak g.$ We let
 $$ \chi_{c,\frak k} : K \longrightarrow \mathbb  {C}^{\times}_1 $$
 be the unitary character of $K$ defined by
\begin{equation}
\chi_{c,\frak k}(\text{exp}\, X(0,\beta,\gamma)) =
 e^{2\pi i <F, X(0,\beta,\gamma)>}= e^{2 \pi i \sigma(c\gamma)} .
 \end{equation}
 The Kirillov correspondence says that the irreducible unitary representation
 $\pi_{c,\frak k}$ of $G$ corresponding to the coadjoint orbit
 $\Omega_F =\Omega_c$ is given by
 \begin{equation}
 \pi_{c,\frak k} = \text{Ind}_K^G\,\chi_{c,\frak k}.
 \end{equation}
 According to Kirillov's Theorem (cf. \cite{Ki1}), we know that the induced
 representation $\pi_{c,\frak k}$ is, up to equivalence, independent of the
choice of a polarization of $\frak g$ for $F.$ Thus we denote the
equivalence class of $\pi_{c,\frak k}$ by $\pi_c.$ $\pi_c$ is
realized on the representation space $L^2(\mathbb{R}^{(h,g)}, d\xi)$ as
follows:
\begin{equation}
(\pi_c(g)f)(\xi) = e^{2\pi i
\sigma\{c(\kappa+\mu^t\lambda + 2 \xi^t\mu)\}}f(\xi+\lambda),
\end{equation}
 where $g=(\lambda, \mu, \kappa) \in G$ and $\xi \in
\mathbb{R}^{(h,g)}.$ Using the fact that $$ \text{exp}\,
X(\alpha,\beta,\gamma) =(\alpha, \beta,
\gamma+\frac12(\alpha^t\beta- \beta^t\alpha)),$$ we see that
$\pi_c$ is nothing but the Schr\"{o}dinger representation
$U(\sigma_c)$ of $G$ induced from the one-dimensional unitary
representation $\sigma_c$ of $K$ given by
$\sigma_c((0,\mu,\kappa)) = e^{2 \pi i \sigma(c\kappa)}I$. We note
that $\pi_c$ is the non-degenerate representation of $G$ with
central character $\chi_c : Z \longrightarrow \mathbb
{C}^{\times}_1$ given by $\chi_c((0,0,\kappa))=e^{2\pi i
\sigma(c\kappa)}.$ Here $Z=\{(0,0,\kappa) \vert \kappa={}^t\kappa
\in \mathbb{R}^{(h,h)} \}$ denotes the center of $G.$ \\
\indent It is well
known that the monomial representation
$(\pi_c,L^2(\mathbb{R}^{(h,g)},d\xi))$ of $G$ extends to an operator of
operator of trace class
\begin{equation}
\pi_c(\phi) : L^2(\mathbb{R}^{(h,g)},d\xi)\longrightarrow
L^2(\mathbb{R}^{(h,g)},d\xi)
\end{equation}
for all $\phi \in C_c^{\infty}(G).$ Here $C_c^{\infty}(G)$ is the
vector space of all smooth functions on $G$ with compact support. We let
$C_c^{\infty}(\frak g)$ and $C(\frak g^*)$ the vector space of all
smooth functions on $\frak g$ with compact support and the vector
space of all continuous functions on $\frak g^*$ respectively. If $f \in \mathcal{C}_{c}^{\infty}(\mathfrak{g})$,
we define the Fourier cotransform
$$
{\mathcal C}F_{\frak g} : C_c^\infty(\frak g)\longrightarrow C({\frak g}^*)
$$
by
\begin{equation}
\left({\mathcal C}F_{\frak g}(f)\right)(F') := \int_{\frak
g}\, f(X)\, e^{2\pi i<F',X>}dX,
\end{equation}
where $F'\in {\frak
g}^*$ and $dX$ denotes the usual Lebesgue measure on $\frak g.$
According to A. Kirillov (cf. \cite{Ki1}), there exists a measure
$\beta$ on the coadjoint orbit $\Omega_c \approx \mathbb{R}^{(h,g)} \times \mathbb{R}^{(h,g)}$ which is invariant
under the coadjoint action og $G$ such that
\begin{equation}
\text{tr}\,\pi_c^1(\phi) = \int_{\Omega_c} {\mathcal C} F_{\frak g} (\phi\circ
\text{exp}) (F') d\beta(F')
\end{equation}
holds for all test
functions $\phi \in C_c^\infty(G),$ where $\text{exp}$ denotes the
exponentional mapping of $\frak g$ onto $G.$ We recall that
$$\pi_c^1(\phi)(f) = \int_G \phi(x) \left(\pi_c (x) f\right) dx,$$
where $\phi \in C_c^\infty(G)$ and $f\in L^2(\mathbb{R}^{(h,g)},d\xi).$
By the Plancherel theorem, the mapping
$$S(G\slash Z) \ni \varphi \longmapsto \pi_c^1(\varphi) \in
TC(L^2(\mathbb{R}^{(h,g)}, d\xi))$$ extends to a unitary isometry
\begin{equation}
\pi_c^2 : L^2(G\slash Z, \chi_c) \longrightarrow
HS(L^2(\mathbb{R}^{(h,g)},d\xi))
\end{equation}
 of the representation space
$L^2(G\slash Z), \chi_c)$ of $\text{Ind}^G_Z\,\chi_c$ onto the
complex Hilbert space\\
\noindent $HS(L^2(\mathbb{R}^{(h,g)},d\xi))$ consisting of all
Hilbert-Schmidt operators on $L^2(\mathbb{R}^{(h,g)},d\xi),$ where
$S(G\slash z)$ is the Schwartz space of all infinitely
differentiable complex-valued functions on $G\slash Z \cong
\mathbb{R}^{(h,g)} \times \mathbb{R}^{(h,g)}$ that are rapidly
decreasing at infinity and $TC(L^2(\mathbb{R}^{(h,g)},d\xi))$
denotes the complex vector space of all continuous $\mathbb
{C}$-linear mappings of $L^2(\mathbb{R}^{(h,g)},d\xi)$ into itself
which are of trace
class. \\
\indent In summary, we have the following result.
\vspace{0.1in}\\
\noindent {\bf Theorem 8.4.}\quad {\it For $F=F(a,b,0) \in {\frak
g}^*,$ the irreducible unitary representation $\pi_{a,b}$ of $G$
corresponding to the coadjoint orbit $\Omega_F = \Omega_c$ under
the Kirillov correspondence is degenerate representation of $G$
given by
$$\pi_{a,b} (\text{exp}\,X(\alpha,\beta,\gamma)) =
e^{4\pi i \sigma({}^ta\alpha -{}^tb\beta)}.$$
 On the other hand,
for $F=F(0,0,c)\in {\frak g}^*$ with $0\ne c ={}^t c \in
\mathbb{R}^{(h,h)},$ the irreducible unitary representation
$(\pi_c,L^2(\mathbb{R}^{(h,g)},d\xi))$ of $G$ corresponding to the
coadjoint orbit $\Omega_c$ under the Kirillov correspondence is unitary equivalent to the Schr\"{o}dinger representation
$U(\sigma_c), L^2(\mathbb{R}^{(h,g)},d\xi))$ and this non-degenerate
representation $\pi_{c}$ is square integrable module its center $Z$.
For all test functions $\phi \in C_c^\infty(G),$ the character
formula
$$\text{tr}\, \pi_c^2(\phi) = {\mathcal C}(\phi, c)\,
 \int_{\mathbb{R}^{(h,g)}}\,\phi(0,0,\kappa)\,e^{2\pi i \sigma(c\kappa)} d\kappa$$
holds for some constant ${\mathcal C}(\phi,c)$ depending on $\phi$ and
$c,$ where $d\kappa$ is the Lebesgue measure on the Euclidean
space $\mathbb{R}^{(h,h)}.$}

\smallskip

\indent   Now we consider the subgroup $K$ of $G$ given by
  $$K = \{ (0,0,\kappa) \in G \,\vert \, \mu \in \mathbb{R}^{(h,g)}, \; \kappa ={}^t\kappa
  \in \mathbb{R}^{(h,h)} \}.$$
The Lie algebra $\frak k$ of $K$ is given by (8.35). The dual
space $\frak k^*$ of $\frak k$ may be identified with the space
$$\{ F(0,b,c) \,\vert \, b \in \mathbb{R}^{(h,g)},\, c={}^tc \in
\mathbb{R}^{(h,h)} \}.$$
We let $\text{Ad}_K^* : K \longrightarrow
GL({\frak k}^*)$ be the coadjoint representation of $K$ on $\frak
k^*.$ The coadjoint orbit $\omega_{b,c}$ of $K$ at $F(0,b,c)\in
{\frak k}^*$ is given by
\begin{equation}
\omega_{b,c} = Ad^*_K (K) \,F(0,b,c) =
\{ F(0,b,c)\}, \;\;\text{a single point}.
\end{equation}
Since $K$ is a commutative group, $[\frak k, \frak k]=0$ and so
the alternating $\mathbb{R}$-bilinear form ${\bold B}_f$ on $\frak
k$ associated to $F:=F(0,b,c)$ identically vanishes on ${\frak k}
\times {\frak k}$(cf. (8.29)). $\frak k$ is the unique
polarization of $\frak k$ for $F=F(0,b,c).$ The Kirillov
correspondence says that the irreducible unitary representation
$\chi_{b,c}$ of $K$ corresponding to the coadjoint orbit
$\omega_{b,c}$ is given by
\begin{equation}
\chi_{b,c}( \text{exp}\,X(0,\beta,\gamma)) = e^{2\pi
i<F(0,b,c),X(0,\beta,\gamma)>} =e^{2\pi i \sigma(2{}^tb\beta +
c\gamma)}
\end{equation}
or
\begin{equation}
\chi_{b,c}((0,\mu,\kappa))= e^{2 \pi i
\sigma(2{}^tb\mu +c\kappa)},\;\; (0,\mu,\kappa) \in K.
\end{equation}
 For $0\ne c={}^t c \in \mathbb{R}^{(h,h)},$ we let $\pi_c$ be the
Schr\"o{}dinger representation of $G$ given by (8.38). We know
that the irreducible unitary representation of $G$ corresponding
to the coadjoint orbit
$$\Omega_c =\text{Ad}^*_G(G)\,F(0,0,c) =\{
F(a,b,c) \,\vert\, a,b \in \mathbb{R}^{(h,g)}\}.$$
 Let $p:{\frak g}^*
\longrightarrow {\frak k}^*$ be the natural projection
defined by $p(F(a,b,c))=F(o,b,c).$ Obviously we have
$$p(\Omega_c)
= \left \{ F(0,b,c) \,\vert\, b \in \mathbb{R}^{(h,g)} \right\}=
\cup_{b\in \mathbb{R}^{(h,g)}}\,\omega_{b,c}.$$
 According to Kirillov Theorem
(cf. \cite{Ki2} p.249, Theorem1), The restriction $\pi_{c}|_{K}$ of $\pi_{c}$ to $K$ is the direct
integral of all one-dimensional representations $\chi_{b,c}$ of
$K\; (b\in \mathbb{R}^{(h,g)}).$ Conversely, we let
$\chi_{b,c}$ be the element of $\hat{K}$ corresponding to the coadjoint
orbit $\omega_{b,c}$ of $K.$ The induced representation
$\text{Ind}_K^G\,\chi_{b,c}$ is nothing but the Schr\"{o}dinger
representation $\pi_c.$ The coadjoint orbit $\Omega_c$ of $G$ is
the only coadjoint orbit such that $\Omega_c \cap
p^{-1}(\omega_{b,c})$ is nonempty.
\vspace{0.2in}\\
\noindent{\bf 9.\ The Jacobi Group} \setcounter{equation}{0}
\renewcommand{\theequation}{9.\arabic{equation}}
\vspace{0.1in}\\
\indent In this section, we study the unitary representations of
the Jacobi group which is a semi-product of a a symplectic group
and a Heisenberg group, and their related topics. In the
subsection 9.1, we present basic ingredients of the Jacobi group
and the Iwasawa decomposition of the Jacobi group. In the
subsection 9.2, we find the Lie algebra of the Jacobi group in
some detail. In the subsection 9.3, we give a definition of Jacobi
forms. In the subsection 9.4, we characterize Jacobi forms as
functions on the Jacobi group satisfying certain conditions. In
the subsection 9.5, we review some results on the unitary
representations, in particular, the Weil representation of the
Jacobi group. Most of the materials here are contained in
\cite{T1}-\cite{T3}. In the subsection 9.6, we describe the
duality theorem for the Jacobi group. In the final subsection, we
study the coadjoint orbits for the Jacobi group and relate these
orbits to the unitary representations of the Jacobi group.
\vspace{0.2in}\\
\noindent{\bf 9.1 The Jacobi Group $G^J$ }
\vspace{0.1in}\\
\indent In this section, we give the standard coordinates of the
Jacobi group $G^J$ and an Iwasawa decomposition of
$G^J.$

\medskip

\noindent {\bf {\scshape 9.1.1. The Standard\
Coordinates\ of\ the\ Jacobi\ Group}} $G^J$

\medskip

\indent Let $m$ and $n$ be two fixed positive integers.
Let
$$Sp(n,\mathbb{R})=\{ M\in \mathbb{R}^{(2n,2n)}\,\ | \,\ ^tM{J_n}M= J_n \}$$
be the symplectic group of degree $n$, where
$$J_n=\left(\begin{matrix} 0&E_n\\
                   -E_n&0\end{matrix}\right)$$
is the symplectic matrix of degree $n$. We let
$$H_n=\{\,Z\in \mathbb{C}^{(n,n)}\,\vert \ Z=\,^tZ,\ \ \text{Im}\,Z\,>\,0\ \}$$
be the Siegel upper half plane of degree $n$. Then
it is easy to see that $Sp(n,\mathbb{R})$ acts on $H_n$ transitively
by
\begin{equation}
M<Z>=(AZ+B)(CZ+D)^{-1},
\end{equation}
where $M=\left( \begin{matrix} A&B \\ C&D\end{matrix}\right)\in Sp(n,\mathbb{R})$ and $Z\in H_n.$
\smallskip

\indent We consider the Heisenberg group
$$
H_{\mathbb{R}}^{(n,m)}=\left\{\,(\lambda,\mu,\kappa) \biggl{|}
 \lambda,\mu\in \mathbb{R}^{(m,n)},\
\kappa\in \mathbb{R}^{(m,m)},\ \kappa+\mu\,^t\!\lambda\ \text{symmetric}\ \right\}
$$
endowed with the following multiplication law
\begin{equation}
(\lambda,\mu,\kappa)\circ (\lambda',\mu',\kappa')=(\lambda+\lambda',\mu+\mu',\kappa+
\kappa'+\lambda\,^t\!\mu'-\mu\,^t\!\lambda').
\end{equation}
We already studied this Heisenberg group in the previous section.

\smallskip

\indent Now we let
$$G^J_{n,m}:=Sp(n,\mathbb{R})\ltimes H_{\mathbb{R}}^{(n,m)}$$
the semidirect product of the symplectic group $Sp(n,\mathbb{R})$ and
the Heisenberg group $H_{\mathbb{R}}^{(n,m)}$
endowed with the following multiplication law
\begin{equation}
(M,(\lambda,\mu,\kappa))\cdot(M',(\lambda',\mu',\kappa'))
=\, (MM',(\tilde{\lambda}+\lambda',\tilde{\mu}+
\mu', \kappa+\kappa'+\tilde{\lambda}\,^t\!\mu'
-\tilde{\mu}\,^t\!\lambda'))
\end{equation}
with $M,M'\in Sp(n,\mathbb{R}), (\lambda,\mu,\kappa),\,(\lambda',\mu',\kappa')
\in H_{\mathbb{R}}^{(n,m)}$
and $(\tilde{\lambda},\tilde{\mu}):=(\lambda,\mu)M'$. We call $G^J_{n,m}$
the {\it Jacobi\ group} of {\it degree} $(n,m).$ If there is no confusion
about the degree $(n,m),$ we write $G^J$ briefly instead of $G^J_{n,m}.$
It is easy to see that $G^J$
acts on $H_{n,m}:=H_n\times \mathbb{C}^{(m,n)}$ transitively by
\begin{equation}
(M,(\lambda,\mu,\kappa))\cdot (Z,W)=(M<Z>,(W+\lambda Z+\mu)
(CZ+D)^{-1}),
\end{equation}
where $M=\left(\begin{matrix} A&B\\ C&D\end{matrix}\right)
\in Sp(n,\mathbb{R}),\ (\lambda,\mu,
\kappa)\in H_{\mathbb{R}}^{(n,m)}$ and $(Z,W)\in H_{n,m}.$

\medskip

\indent Now we define the linear mapping
$$Q\,:\,({\mathbb{R}}^{(m,n)}\times \mathbb{R}^{(m,n)})\times
({\mathbb{R}}^{(m,n)}\times \mathbb{R}^{(m,n)})   \longrightarrow \mathbb{R}^{(m,m)}$$
by
$$
Q((\lambda,\mu),(\lambda',\mu'))=\lambda{^t\!\mu'}-\mu{^t\!\lambda'},\ \ \lambda,\mu,\lambda',\mu'
\in \mathbb{R}^{(m,n)}.
$$
Clearly we have\\
\begin{equation}
^t\!Q(\xi,\eta)=-Q(\eta,\xi) \,\ \text{for all} \xi,\eta\in
\mathbb{R}^{(m,n)}\times \mathbb{R}^{(m,n)},
\end{equation}
\begin{equation}
Q(\xi M,\eta M)=Q(\xi,\eta) \,\ \text{for all} M\in Sp(n,\mathbb{R}).
\end{equation}

\smallskip

\indent For a reason of the convenience, we write an element of $G^J$ as
$$
g=[M,(\lambda,\mu,\kappa)]:=(E_{2n},(\lambda,\mu,\kappa))\cdot (M,(0,0,0)).
$$
Then the multiplication becomes
$$
[M,(\xi,\kappa)]\circ [M',(\xi',\kappa')]=[MM',(\xi+\xi'M^{-1},\kappa+\kappa'+
Q(\xi,\xi'M^{-1}))].
$$
For brevity, we set $G=Sp(n,\mathbb{R}).$
We note that the stabilizer $K$ of $G$ at $iE_n$
under the symplectic action is
given by
$$K=\left\{\begin{pmatrix} A&-B\\ B&A \end{pmatrix}
\in GL(2n,\mathbb{R})\,\bigg|\ \, ^t\!AB=\,^tBA,\ \
^t\!AA+\,^tBB=E_n\ \right\}$$
and is a maximal compact subgroup of $G$. We also recall that the Jacobi
group $G^J$ acts on $H_{n,m}$ transitively via (9.4). Then it is easy
to see that  the stabilizer
$K^J$ of $G^J$ at $(iE_n,0)$ under this action is given by
\begin{align*}
K^J & =\left\{ [k,(0,0,\kappa)] \biggl{|} \ k \in K, \,\   k = ^tk
\in \mathbb{R}^{(m,m)}
\right\}\\
& \cong K\times \{ \,\ (0,0,\kappa)| \,\ \kappa=^t{\kappa} \in
\mathbb{R}^{(m,m)} \,\ \} \cong K \times
\text{Symm}^2(\mathbb{R}^m).
\end{align*}

\noindent Thus on $G^J/K^J\cong H_{n,m},$ we have the coordinate
$$
g\cdot (iE_n,0):=(Z,W):=(X+iY,\lambda Z+\mu),\ \ \ g\in G^J.
$$
In fact, if $g=[M,(\lambda,\mu,\kappa)]\in G^J$ with $M=\begin{pmatrix} A&B\\ C&D
\end{pmatrix}
\in G,$
\begin{align*}
\ \ \ Z&=M<iE_n>=(iA+B)(iC+D)^{-1}=X+iY,\\
\ \ \ W&=\{i(\lambda A+\mu C)+\lambda B+\lambda D\}(iC+D)^{-1}\\
&=\{\lambda(iA+B)+\mu(iC+D)\}(iC+D)^{-1}\\
&=\lambda Z+\mu.
\end{align*}
We set
$$
dX=\begin{pmatrix} dX_{11}&\ldots&dX_{1n}\\
\vdots&\ddots&\vdots\\
dX_{n1}&\ldots&dX_{nn}\end{pmatrix}, \ \ dW=
\begin{pmatrix} dW_{11}&\ldots&dW_{1n}\\
\vdots&\ddots&\vdots\\
dW_{m1}&\ldots&dW_{mn}\end{pmatrix}
$$
and
$$
{\partial\over {\partial X}}=\begin{pmatrix} {{\partial}\over {\partial X_{11}}}&\ldots
&{{\partial}\over {\partial X_{1n}}}\\
\vdots&\ddots&\vdots\\
{{\partial}\over {\partial X_{n1}}}&\ldots&{{\partial}\over {\partial X_{nn}}}\end{pmatrix},\ \ \
{{\partial}\over {\partial W}}=\begin{pmatrix} {{\partial}\over {\partial W_{11}}}&\ldots &
{{\partial}\over {\partial W_{m1}}}\\
\vdots&\ddots&\vdots\\
{{\partial}\over {\partial W_{1n}}}&\ldots &{{\partial}\over {\partial W_{mn}}}
\end{pmatrix}.
$$
Similarly we set $dY=(dY_{ij}),\ d\lambda=(d\lambda_{pq}),\ d\mu=(d\mu_{pq}),\cdots
$etc.\ \, By an easy calculation, we have
\begin{align*}
\displaystyle {{\partial}\over {\partial W}}&=\displaystyle{{1\over {2i}}\,Y^{-1}\left( {{\partial}\over {\partial \lambda}}-
\bar{Z}{{\partial}\over {\partial \mu}} \right)},\\
\displaystyle{{\partial}\over {\partial \overline{W}}}&=
\displaystyle{{i\over 2}\,Y^{-1}\left( {{\partial}\over {\partial \lambda}}-
Z{{\partial}\over {\partial \mu}}\right)},\\
\displaystyle{{\partial}\over {\partial X}}&=\displaystyle{{\partial}\over {\partial Z}}+{{\partial}\over {\partial \bar{Z}}}+
 {{\partial}\over {\partial W}}\lambda+{{\partial}\over {\partial \overline{W}}}\lambda,\\
\displaystyle{{\partial}\over {\partial Y}}&=\displaystyle i{{\partial}\over {\partial Z}}-i{{\partial}\over {\partial \bar{Z}}}
+i{{\partial}\over {\partial W}}\lambda-i{{\partial}\over {\partial \overline{W}}}\lambda.
\end{align*}
We set
$$P_+={\frac 12}\left({{\partial}\over {\partial
X}}-i{{\partial}\over {\partial Y}}\right)= {{\partial}\over
{\partial Z}}+{{\partial}\over {\partial
W}}\,(\text{Im}\,W)\,Y^{-1}$$ and
$$P_{-}={\frac 12}\left( {{\partial}\over {\partial X}}+i{{\partial}\over {\partial Y}}\right)=
{{\partial}\over {\partial \bar{Z}}}+{{\partial}\over {\partial \overline{W}}}\,
(\text{Im}\,W)\,Y^{-1}.$$
Let ${\frak g}$ be the Lie algebra of $G$ and ${\frak g}_{\mathbb{C}}$ its
complexification. Then
$${\frak g}_{\mathbb{C}}=\left\{\begin{pmatrix} A & B\\ C & -^t{A} \end{pmatrix}
\in \mathbb{C}^{(2n,2n)}\,
\bigg|\
B=\,^tB,\ \ C=\,^tC\ \right\}.$$
We let $\hat{J}:=iJ_n$ with $J_n=\begin{pmatrix} 0 & E_n\\ -E_n & 0 \end{pmatrix}.$
We define an involution $\sigma$ of $G$ by
\begin{equation}
\sigma(g):=\hat{J}g\hat{J}^{-1},\ \ \ g\in G.
\end{equation}
The differential map $d\sigma=\text{Ad}\,(\hat{J})$ of $\sigma$ extends complex linearly
to the complexification ${\frak g}_{\mathbb{C}}$ of ${\frak g}.
\ \text{Ad}\,(\hat{J})$ has
1 and -1 as eigenvalues. The $(+1)$-eigenspace of $\text{Ad}\,(\hat{J})$
is given by
\begin{equation}
{\frak k}_{\mathbb{C}}=\left\{
\begin{pmatrix} A&-B\\ B & A\end{pmatrix} \in \mathbb{C}^{(2n,2n)}\ \bigg|\
^t\!A+A=0,\ B=\,^tB\ \right\}.
\end{equation}
We note that ${\frak k}_{\mathbb{C}}$ is the complexification of the Lie algebra
${\frak k}$ of a maximal compact subgroup $K=G\cap SO(2n,\mathbb{R})\cong U(n)$ of
$G$. The $(-1)$-eigenspace of $\text{Ad}\,(\hat{J})$ is given by
\begin{equation}
{\frak p}_{\mathbb{C}}=\left\{ \begin{pmatrix} A & B\\ B & -A\end{pmatrix}
\in \mathbb{C}^{(2n,2n)}\
\bigg|\
A=\,^t\!A,\ B=\,^tB\ \right\}.
\end{equation}
We observe that ${\frak p}_{\mathbb{C}}$ is not a Lie algebra.
But ${\frak p}_{\mathbb{C}}$
has the following decomposition
$${\frak p}_{\mathbb{C}}={\frak p}_{+}\oplus {\frak p}_{-},$$
where
\begin{equation}
{\frak p}_{+}=\left\{\begin{pmatrix} X & iX \\ iX & -X
\end{pmatrix} \in \mathbb{C}^{(2n,2n)}\
\bigg|\
X=\,^tX\ \right\}
\end{equation}
and
\begin{equation}
{\frak p}_{-}=\left\{\begin{pmatrix} Y & -iY\\ -iY & -Y
\end{pmatrix} \in \mathbb{C}^{(2n,2n)}\
\bigg|\
Y=\,^tY\ \right\}.
\end{equation}
We observe that ${\frak p}_+$ and ${\frak p}_{-}$ are abelian subalgebras of
${\frak g}_{\mathbb{C}}.$
Since $\text{Ad}\,(\hat{J})[X,Y]=[\text{Ad}\,(\hat{J})X,\,
\text{Ad}\,
(\hat{J})Y]$ for
all $X,Y\in {\frak g}_{\mathbb{C}},$ we have
\begin{equation}
[{\frak k}_{\mathbb{C}},{\frak k}_{\mathbb{C}}]\subset {\frak k}_{\mathbb{C}},\ \
[{\frak k}_{\mathbb{C}},{\frak p}_{\mathbb{C}}]\subset {\frak p}_{\mathbb{C}},\ \
[{\frak p}_{\mathbb{C}},{\frak p}_{\mathbb{C}}]\subset {\frak k}_{\mathbb{C}}.
\end{equation}
Since $\text{Ad}(k)X=kXk^{-1}\,(\,k\in K,\ X\in {\frak g}_{\mathbb{C}}\,),$ we obtain
\begin{equation}
\text{Ad}(k){\frak p}_+\subset {\frak p}_+,\ \ \
\text{Ad}(k){\frak p}_{-}\subset
{\frak p}_{-}.
\end{equation}
For instance, if $k=\begin{pmatrix} A & B\\ -B & A\end{pmatrix} \in K,$ then
\begin{equation}
\text{Ad}(k)\begin{pmatrix} X & \pm iX\\ \pm iX & -X
\end{pmatrix} =\begin{pmatrix} X'&\pm iX'\\
\pm iX' & -X'\end{pmatrix},\ \ \ X=\,^tX,
\end{equation}
where
$$X'=(A+iB)X\,^t(A+iB).$$
If we identify ${\frak p}_{-}$ with $\text{Symm}^2(\mathbb{C}^n)$ and $K$
with $U(n)$ as a
subgroup of $GL(n,\mathbb{C})$ via the mapping $K\ni \begin{pmatrix} A & -B\\
B&A \end{pmatrix} \longrightarrow A+iB\in U(n),$
then the action of $K$ on ${\frak p}_{-}$ is
compatible with the natural representation $\rho^{[1]}$ of $GL(n,\mathbb{C})$ on
$\text{Symm}^2(\mathbb{C}^n)$ given by
$$\rho^{[1]}(g)X=gX\,^tg,\ \ g\in GL(n,\mathbb{C}),\ \ X\in \text{Symm}^2(\mathbb{C}^n).$$
The Lie algebra ${\frak g}$ of $G$ has a Cartan decomposition
\begin{equation}
{\frak g}={\frak k}\oplus {\frak p},
\end{equation}
where
\begin{align*}
{\frak k}&=\left\{ \begin{pmatrix} A & B\\ -B & A
\end{pmatrix} \in \mathbb{R}^{(2n,2n)}\ \bigg|\
A+{^t\!A}=0,\ \ B=\,^tB\ \right\},\\
{\frak p}&=\left\{ \begin{pmatrix} A & B\\ B &-A
\end{pmatrix} \in \mathbb{R}^{(2n,2n)}\ \bigg|\
A={^t\!A},\ \ B=\,^tB\ \right\}.
\end{align*}
Then $\theta:=\text{Ad}\,(\hat{J})$ is a Cartan involution because
$$-B(W,\theta(W))=-B(X,X)+B(Y,Y) > 0$$
for all $\,W=X+Y,\ X\in {\frak k},\
Y\in {\frak p}.$
Here $B$ denotes the Cartan-Killing form for ${\frak g}.$ \\
\noindent Indeed,
\begin{equation}
B(X,Y)=2(n+1)\,\sigma(XY),\ \ \ X,Y \in {\frak g}.
\end{equation}
The vector space ${\frak p}$ is identified with the tangent space of $H_n$
at $iE_n.$ The correspondence
\begin{equation}
{1\over 2}\begin{pmatrix} B & A\\ A &-B
\end{pmatrix} \longmapsto A+iB
\end{equation}
yields an isomorphism of ${\frak p}$ onto $\text{Symm}^2(\mathbb{C}^n).$
The Lie algebra ${\frak g}^J$ of the Jacobi group $G^J$ has a decomposition
\begin{equation}
{\frak g}^J={\frak k}^J + {\frak p}^J,
\end{equation}
where
\begin{align*}
{\frak k}^J&=\left\{(X,(0,0,\kappa )\,| \,\ X \in {\frak k}, \,\ \kappa=
\,^t\kappa \in \mathbb{R}^{(m,m)}
\right\},\\
{\frak p}^J&=\left\{(Y,(P,Q,0)\ | \,\ Y\in {\frak p}, \,\ P,Q\in \mathbb{R}^{(m,n)}
\right\}.
\end{align*}
Thus the tangent
space of the homogeneous space $H_{n,m}\cong G^J/K^J$ at $(iE_n,0)$
is given by
$${\frak p}^J\cong{\frak p}\oplus (\mathbb{R}^{(n,m)}\times \mathbb{R}^{(n,m)})\cong
{\frak p}\oplus \mathbb{C}^{(n,m)}.$$
We define a complex structure $I^J$ on the tangent space ${\frak p}^J$ of
$H_{n,m}$ at $iE_n$ by
\begin{equation}
I^J \left( \begin{pmatrix} Y & X\\ X &-Y
\end{pmatrix},(P,Q)\right):=
\left(\begin{pmatrix} X & -Y\\ -Y & -X
\end{pmatrix},(Q,-P)\right).
\end{equation}
Identifying $\mathbb{R}^{(m,n)}\times \mathbb{R}^{(m,n)}$ with $\mathbb{C}^{(m,n)}$ via
\begin{equation}
(P,Q) \longmapsto iP+Q,\ \ P,Q\in \mathbb{R}^{(m,n)},
\end{equation}
we may regard the complex structure $I^J$ as a real linear map
\begin{equation}
I^J(X+iY,Q+iP)=(-Y+iX,-P+iQ),
\end{equation}
where $X+iY\in \text{Symm}^2(\mathbb{C}^n),\ Q+iP\in \mathbb{C}^{(m,n)}.\ I^J$
extends complex
linearly on the complexification ${\frak p}_{\mathbb{C}}^J={\frak p}\otimes_{\mathbb{R}}
\mathbb{C}$ of ${\frak p}.\ {\frak p}_{\mathbb{C}}$ has a decomposition
\begin{equation}
{\frak p}_{\mathbb{C}}={\frak p}_{+}^J\oplus {\frak p}_{-}^J,
\end{equation}
where ${\frak p}_+^J\,(\text{resp.}\ {\frak p}_{-}^J\,)$
denotes the $(+i)$-eigenspace
\,$(\text{resp.}\ (-i)$-eigenspace) of $I^J$.
Precisely, both ${\frak p}_+^J$ and
${\frak p}_{-}^J$ are given by
$${\frak p}_+^J=\left\{\left( \begin{pmatrix} X & iX \\ iX & -X
\end{pmatrix},(P,iP)
\right)\ \bigg|\
X\in \text{Symm}^2(\mathbb{C}^n),\ \ P\in \mathbb{C}^{(m,n)}\right\}$$
and
$${\frak p}_{-}^J=\left\{\left( \begin{pmatrix} X & -iX \\ -iX & -X
\end{pmatrix},(P,-iP)
\right)\ \bigg|
\ X\in \text{Symm}^2(\mathbb{C}^n),\ \ P\in \mathbb{C}^{(m,n)}\right\}.$$
With respect to this complex structure $I^J,$ we may say that $f$ is
{\it holomorphic} if and only if $\xi f=0$ for all
$\xi\in {\frak p}_{-}^J.$

\medskip

We fix $g=[M,(\l,\mu;\kappa)]\in G^J$ with $M= \begin{pmatrix} A & B\\
C & D\end{pmatrix} \in G.$ Let $T_g:H_n \longrightarrow H_n$ be
the mapping defined by (9.4). We consider the behavior of the
differential map $dT_g$ of $T_g$ at $(iE_n,0)$
$$dT_g:T_{(iE_n,0)}(H_{n,m})\longrightarrow T_{(Z,W)}(H_{n,m}),\ \ (Z,W):=g\cdot
(iE_n,0).$$
Now we let $\alpha(t)=(Z(t),\xi(t))$ be a smooth curve in $H_{n,m}$ passing
through $(iE_n,0)$ with $\alpha'(0)=(V,iP+Q)\in T_{(iE_n,0)}(H_{n,m}).$ Then
\begin{align*}
\gamma(t):&=g\cdot \alpha(t)=(Z(g;t),\xi(g;t))\\
            &=(M<Z(t)>,(\xi(t)+{\tilde {\lambda}}Z(t)+{\tilde {\mu}})
               (CZ(t)+D)^{-1})
\end{align*}
is a curve in $H_{n,m}$ passing through $\gamma(0)=(Z,W)$ with $({\tilde {\lambda}},
{\tilde {\mu}})=(\lambda,\mu)M.$ Using the relation
$${{\partial}\over {\partial t}}\Bigg|_{t=0}(CZ(t)+D)^{-1}=
-(iC+D)^{-1}CZ'(0)(iC+D)^{-1},$$
we have
\begin{align*}
\displaystyle{{\partial}\over {\partial t}}\Bigg|_{t=0}Z(g;t)&=AZ'(0)(iC+D)^{-1}+(iA+B)
\,\displaystyle{{\partial}\over {\partial t}}\Bigg|_{t=0}(CZ(t)+D)^{-1}\\
&=AZ'(0)(iC+D)^{-1}-(iA+B)(iC+D)^{-1}CZ'(0)(iC+D)^{-1}\\
&=\{ A\,^t(iC+D)-(iA+B)\,^tC\}\,^t(iC+D)^{-1}Z'(0)(iC+D)^{-1}\\
&=\,^t(iC+D)^{-1}Z'(0)(iC+D)^{-1}.
\end{align*}
and
\begin{align*}
\displaystyle{{\partial}\over {\partial t}}\Bigg|_{t=0}\xi(g;t)
&=(\xi'(0)+{\tilde {\lambda}}Z'(0))
(iC+D)^{-1}\\
&\ \ \ +(\xi(0)+i{\tilde {\lambda}}+{\tilde {\mu}})
\,\displaystyle{{\partial}\over {\partial t}}\Bigg|_{t=0}
(CZ(t)+D)^{-1}\\
&=(iP+Q+{\tilde {\lambda}}Z'(0))(iC+D)^{-1}\\
& \ \ \ -(i{\tilde {\lambda}}+{\tilde {\mu}})(iC+D)^{-1}CZ'(0)(iC+D)^{-1}\\
&=(iP+Q)(iC+D)^{-1}\\
& \ \ \
+\left\{ {\tilde \lambda}\,^t(iC+D)-(i{\tilde \lambda}+{\tilde \mu})\,^tC\right\}
\,^t(iC+D)^{-1}Z'(0)(iC+D)^{-1}\\
&=(iP+Q)(iC+D)^{-1}+\lambda\,^t(iC+D)^{-1}Z'(0)(iC+D)^{-1}.
\end{align*}
Here we used the fact that $(iC+D)^{-1}C$ is symmetric and the relation
$${\tilde \lambda}=\lambda A+\mu C,\ \ \ {\tilde \mu}=\lambda B+\mu D.$$
Therefore we obtain
\begin{align*}
\ \ Z'(g;0)&=\,^t(iC+D)^{-1}Z'(0)(iC+D)^{-1},\\
\ \
\xi'(g;0)&=\xi'(0)(iC+D)^{-1}+\lambda\,^t(iC+D)^{-1}Z'(0)(iC+D)^{-1}.
\end{align*}
\indent In summary, we have
\vspace{0.1in}\\
\noindent {\bf Proposition\ 9.1.}\quad {\it Let $g=[M,(\lambda,\mu;\kappa)]\in G^J$ with
$M=\begin{pmatrix} A & B\\ C&D\end{pmatrix} \in G$ and let $(Z,W)=g\cdot (iE_n,0).$
Then the differential map $dT_g:T_{(iE_n,0)}(H_{n,m}) \longrightarrow
T_{(Z,W)}(H_{n,m})$ is
given by}
\begin{equation}
(v,w)\longmapsto (v(g),w(g)),\ \ \ v\in \text{Symm}^2(\mathbb{C}^n),\ \ w\in \mathbb{C}^{(m,n)}
\end{equation}
{\it with}
\begin{align*}
\ \ \ v(g)&=\,^t(iC+D)^{-1}v(iC+D)^{-1},\\
\ \ \ w(g)&=w(iC+D)^{-1}+\lambda\,^t(iC+D)^{-1}v(iC+D)^{-1}.
\end{align*}

\medskip

\noindent {\bf {\scshape 9.1.2. An Iwasawa\ Decomposition\ of\ the\ Jacobi\ Group} $G^J$}

\medskip

\indent First of all, we give the Iwasawa decomposition of $G=Sp(n,\mathbb{R}).$
For a positive diagonal matrix $H$ of degree $n$, we put
$$t(H)=\begin{pmatrix} H & 0\\ 0 & H^{-1}\end{pmatrix}$$
and for an upper triangular matrix $A$ with $1$ in every diagonal entry and
$B\in \mathbb{R}^{(m,n)},$ we write
$$n(A,B)=\begin{pmatrix} A & B\\ 0 & ^tA^{-1}\end{pmatrix}.$$
We let $A$ be the set of such all $t(H)$ and let $N$ be the set of
such all $n(A,B)$ such that $n(A,B)\in G$, namely,
$A\,^tB=B\,^tA.$ It is clear that $A$ is an abelian subgroup of
$G$ and $N$ is a nilpotent subgroup of $G.$ Then we have the
so-called {\it Iwasawa\ decomposition}
\begin{equation}
G=NAK=KAN.
\end{equation}
\indent Now we define the subgroups $A^J,\,N^J$ and ${\tilde N}^J$ of $G^J$ by
\begin{align*}
A^J&=\left\{\,t(H,\lambda):=[t(H),(\lambda,0,0)]\,\bigg|\ t(H)\in A,\
 \lambda \in \mathbb{R}^{(m,n)}\,\right\},\\
N^J&=\left\{\,n(A,B;\mu):=[n(A,B),(0,\mu,0)]\,\bigg|\
n(A,B)\in N,\ \mu\in \mathbb{R}^{(m,n)}\,
\right\}
\end{align*}
and
$${\tilde N}^J:=\left\{\,{\tilde n}(A,B;\mu,\kappa)=[n(A,B),(0,\mu,\kappa)]\,\bigg|\
n(A,B)\in N,\ \mu \in \mathbb{R}^{(m,n)},\ \kappa \in \mathbb{R}^{(m,m)}\,
\right\}.$$
For $t(H,\lambda),\ t(H',\lambda')\in A^J,$ we have
$$t(H,\lambda)\circ t(H',\lambda')=t(HH',\lambda+ \lambda'H^{-1}).$$
Thus $A^J$ is the semidirect product of $\mathbb{R}^{(m,n)}$ and
${\mathbb{ D}}^+$, where ${\mathbb{D}}^+$ denotes the subgroup of
$GL(n,\mathbb{R})$ consisting of positive diagonal matrices of
degree $n$. Furthermore we have for $t(H,\lambda)\in A^J$ and
${\tilde n}(A,B;\mu,\kappa) \in {\tilde N}^J$
$${\tilde n}(A,B;\mu,\kappa)\circ t(H,\lambda)=[n(A,B)t(H),(\lambda A^{-1},\,\mu-\lambda A^{-1}B\,^tA,\,-\mu\,^tA^{-1}\,^t\lambda)]$$
and
\begin{align*}
t(H,\lambda)\circ {\tilde n}(A,B;\mu,\kappa)&=[t(H)n(A,B),(\lambda,\mu H,\,\kappa+\lambda H\,^t\mu)]\\
&=\left[ \begin{pmatrix} HA & HB\\ 0& H^{-1}\,^tA^{-1}\end{pmatrix},(\lambda,\mu H,\kappa+\lambda H\,^t\mu)\right].\end{align*}
Therefore we have
\begin{align*}
& t(H,\lambda)\circ {\tilde n}(A,B;\mu,\kappa)\circ t(H,\lambda)^{-1}\\
=& [t(H)n(X)t(H^{-1}),(\lambda-\lambda HA^{-1}H^{-1},\,\mu H+\lambda HA^{-1}B\,^tA H,\\
 &\ \ \kappa+\lambda HA\,^tB
\,^tA^{-1}H\,^t\lambda-\mu\,^tA^{-1}H\,^t\lambda)]\\
=&[n(HAH^{-1},HBH),\,(\lambda-\lambda HA^{-1}H^{-1},\mu H+\lambda HA^{-1}B\,^tA H,\\
& \ \ \kappa+\lambda HA\,^tB
\,^tA^{-1}H\,^t\lambda-\mu\,^tA^{-1}H\,^t\lambda)].
\end{align*}
Thus there is a decomposition
\begin{equation}
G^J={\tilde N}^JA^JK.
\end{equation}
For $g\in G^J$, one has
\begin{align*}
g & = [n(A,B)t(H)\kappa,(\lambda,\mu,\kappa)], \quad  \kappa \in K \\
& = {\tilde{n}}(A,B;\mu^{\ast},\kappa^{\ast})\circ
t(H,\lambda^{\ast})\circ k
\end{align*}
with
$$\lambda^{\ast}=\lambda H, \quad  \mu^{\ast}=\mu+\lambda B^t A \,\,\ \text{and}  \,\,\ \kappa^{\ast}=\kappa+\mu^t\lambda+
\lambda B\,^t(\lambda A).$$
Recalling the subgroup $K^J$ of $G^J$ defined by
$$K^J=\{\,[k,(0,0,\kappa)]\,\vert \ k\in K,\ \kappa=\,^t\kappa\in \mathbb{R}^{(m,m)}\,\},$$
we also have a decomposition
\begin{equation}
G^J=N^JA^JK^J.
\end{equation}
For $g\in G^J$, one has
\begin{align*}
g&=[n(A,B)t(H)k,(\lambda,\mu,\kappa)]\\
&=n(A,B;{\tilde{\mu}})\circ t(H,{\tilde{\lambda}})\circ [k,(0,0,{\tilde{\kappa}})]
\end{align*}
with
$${\tilde{\lambda}}=\lambda A,\quad    {\tilde{\mu}}=\mu+\lambda A^{-1}B\,^tA,\quad
{\tilde{\kappa}}=\kappa+(\mu+\lambda A^{-1}B\,^tA)\,^t\lambda.$$
We call the decomposition (9.25) or (9.26) an {\it Iwasawa\ decomposition}
of $G^J.$
Finally we note that the decomposition (9.25) or (9.26) may be understood
as the product of the usual Iwasawa decomposition (9.24)
of $G$ with a decomposition
\begin{equation}
H_{\mathbb{R}}^{(n,m)}={\tilde N}_0A_0
\end{equation}
of the Heisenberg group $H_{\mathbb{R}}^{(n,m)}$ into the group $
A_0=\{\,(\lambda,0,0)\,\vert\ \lambda\in \mathbb{R}^{(m,n)}\,\} $
which normalizes the maximal abelian subgroup ${\tilde N}_0 = \{
\,\ (0, \mu, \kappa)\,\ | \,\ \mu \in \mathbb{R}^{(m,n)} , \kappa=
\,\ ^t\kappa \in \mathbb{R}^{(m,n)}\,\}. $
\vspace{0.2in}\\
\noindent{\bf 9.2. The Lie Algebra of the Jacobi Group $G^J$ }
\vspace{0.1in}\\
\indent In this section, we describe the Lie algebra ${\frak g}^J$ of the
Jacobi group $G^J$ explicitly.

\smallskip
\indent First of all, we observe that ${\frak g}$ of $G$ may be
regarded as a subalgebra of ${\frak g}^J$ by identifying ${\frak
g}$ with ${\frak g}\times \{ 0\}$ and the Lie algebra ${\frak h}$
of the Heisenberg group $H_{\mathbb{R}}^{(n,m)}$ may be regarded
as an ideal of ${\frak g}^J$ by identifying ${\frak h}$ with $\{
0\}\times {\frak h}.$  We denote by $E_{ij}$ the matrix with entry
$1$ where the $i$-th row and the $j$-th column meet, all other
entries 0.

\smallskip

\indent For $1\leq a,b,p\leq m,\ 1\leq i,j,q\leq n,$ we set
\begin{align*}
\ \ \ \ \ A_{ij}:&=\begin{pmatrix} E_{ij}+E_{ji} & 0 & 0 & 0\\
0 & 0 & 0 & 0\\ 0 & 0 & -(E_{ij}+E_{ji}) & 0 \\
0 & 0 & 0 & 0 \end{pmatrix}, \\
\ \ \ \ \ B_{ij}:&=\begin{pmatrix} 0 & 0 & E_{ij}+E_{ji} & 0 \\
0 & 0 & 0 & 0 \\ E_{ij}+E_{ji} & 0 & 0 & 0 \\ 0 & 0 & 0 & 0
\end{pmatrix},\\
\ \ \ \ \ S_{ij}:&=\begin{pmatrix} E_{ij}-E_{ji} & 0 & 0 & 0 \\
0 & 0 & 0 & 0 \\ 0 & 0 & E_{ij}-E_{ji} & 0 \\ 0 & 0 & 0 & 0
\end{pmatrix},\\
\ \ \ \ \ T_{ij}:&=\begin{pmatrix} 0 & 0 & E_{ij}+E_{ji} & 0 \\
0 & 0 & 0 & 0 \\ -(E_{ij}+E_{ji}) & 0 & 0 & 0 \\ 0 & 0 & 0 & 0
\end{pmatrix},\\
\ \ \ \ \ D_{ab}^0:&=\begin{pmatrix} 0 & 0 & 0 & 0 \\
0 & 0 & 0 & {\frac 12}(E_{ab}+E_{ba})\\ 0 & 0 & 0 & 0 \\
0 & 0 & 0 & 0\end{pmatrix},\\
\ \ \ \ \ D_{pq}:&=\begin{pmatrix} 0 & 0 & 0 & 0 \\ E_{pq} & 0 & 0 & 0 \\
0 & 0 & 0 & -E_{qp} \\ 0 & 0 & 0 & 0
\end{pmatrix},\\
\ \ \ \ \ {\hat D}_{pq}:&=\begin{pmatrix} 0 & 0 & 0 & E_{qp} \\
0 & 0 & E_{pq} & 0 \\ 0 & 0 & 0 & 0 \\ 0 & 0 & 0 & 0 \end{pmatrix}.
\end{align*}
We observe that the set
$$\left\{ S_{ij},\,T_{kl},\,D_{ab}^0 \bigg|\ 1\leq i < j
\leq n,\ 1\leq k\leq l\leq n,\ 1\leq a\leq b\leq m\ \right\}$$
form a basis of ${\frak k}^J$ and the set
$$\left\{ A_{ij},\,B_{ij},\,D_{pq},\,{\hat D}_{rs}\,\bigg| \ 1\leq i\leq j
\leq n,\ 1\leq p,r\leq m,\ 1\leq q,s\leq n\ \right\}$$
form a basis of
${\frak p}^J$\,(cf. (9.15)). We note that
\begin{align*}
\ \ \ A_{ij}&=A_{ji},\ \ B_{ij}=B_{ji},\ \ S_{ij}=-S_{ji},\ \ T_{ij}=T_{ji},
\\ \ \ \ D_{ab}^0&=D_{ba}^0,\ \ \ D_{pq}^2={\hat D}_{pq}^2=0.
\end{align*}
\vspace{0.1in}\\
\noindent {\bf Lemma\ 9.2.} \quad {\it We have the following
commutation relation :}
\def\d{\delta}
\begin{align*}
 [A_{ij},A_{kl}]&=\d_{ik}S_{jl}+\d_{il}S_{jk}+\d_{jk}S_{il}+
\d_{jl}S_{ik},\\
 [A_{ij},B_{kl}]&=\d_{ik}T_{jl}+\d_{il}T_{jk}+\d_{jk}T_{il}+
\d_{jl}T_{ik},\\
 [A_{ij},S_{kl}]&=\d_{ik}A_{jl}-\d_{il}A_{jk}+\d_{jk}A_{il}-
\d_{jl}A_{ik},\\
 [A_{ij},T_{kl}]&=\d_{ik}B_{jl}+\d_{il}B_{jk}+\d_{jk}B_{il}+
\d_{jl}B_{ik},\\
 [B_{ij},B_{kl}]&=\d_{ik}S_{jl}+\d_{il}S_{jk}+\d_{jk}S_{il}+
\d_{jl}S_{ik},\\
 [B_{ij},S_{kl}]&=\d_{ik}B_{jl}-\d_{il}B_{jk}+\d_{jk}B_{il}-
\d_{jl}B_{ik},\\
 [B_{ij},T_{kl}]&=-\d_{ik}A_{jl}-\d_{il}A_{jk}-\d_{jk}A_{il}-
\d_{jl}A_{ik},\\
 [S_{ij},S_{kl}]&=-\d_{ik}S_{jl}+\d_{il}S_{jk}+\d_{jk}S_{il}-
\d_{jl}S_{ik},\\
 [S_{ij},T_{kl}]&=-\d_{ik}T_{jl}-\d_{il}T_{jk}+\d_{jk}T_{il}+
\d_{jl}T_{ik},\\
 [T_{ij},T_{kl}]&=-\d_{ik}S_{jl}-\d_{il}S_{jk}-\d_{jk}S_{il}-
\d_{jl}S_{ik},\\
 [D_{ab}^0,A_{ij}]&=[D_{ab}^0,B_{ij}]=[D^0_{ab},S_{ij}]=
[D_{ab}^0,T_{ij}]=0,\\
 [D_{ab}^0,D_{cd}^0]&=[D_{ab}^0,D_{pq}]=[D_{ab}^0,{\hat D}_{pq}]=0,\\
 [D_{pq},A_{ij}]&=\d_{qi}D_{pj}+\d_{qj}D_{pi},\\
 [D_{pq},B_{ij}]&=[D_{pq},T_{ij}]=\d_{qi}{\hat D}_{pj}+\d_{qj}
{\hat D}_{pi},\\
 [D_{pq},S_{ij}]&=\d_{qi}D_{pj}-\d_{qj}D_{pi},\\
 [D_{pq},D_{rs}]&=0,\ \ \ [D_{pq},{\hat D}_{rs}]=2\d_{qs}D_{pr}^0,\\
 [{\hat D}_{pq},A_{ij}]&=-\d_{qi}{\hat D}_{pj}-\d_{qj}{\hat D}_{pi},\\
 [{\hat D}_{pq},B_{ij}]&=\d_{qi}D_{pj}+\d_{qj}D_{pi},\\
 [{\hat D}_{pq},S_{ij}]&=\d_{qi}D_{pj}-\d_{qj}{\hat D}_{pi},\\
 [{\hat D}_{pq},T_{ij}]&=-\d_{qi}D_{pj}-\d_{qj}D_{pi},\\
 [{\hat D}_{pq},{\hat D}_{rs}]&=0.
\end{align*}
\noindent Here $1\leq a,b,c,d,p,r\leq m,\ 1\leq i,j,k,l,q,s\leq n$
and $\d_{ij}$ denotes the Kronecker delta symbol.
\vspace{0.05in}\\
\noindent {\it Proof.}\quad The proof follows from a
straightforward calculation. \hfill $\square$
\vspace{0.1in}\\
\noindent {\bf Corollary\ 9.3.}\quad {\it We have the following
relation :}
\begin{align*}
 [{\frak k}^J,{\frak k}^J]&\subset {\frak k}^J,\ \ \ [{\frak k}^J,
{\frak p}^J]\subset {\frak p}^J,\\
 [{\frak p},{\frak h}]&\subset {\frak h},\ \ \ [{\frak h},{\frak h}]
\subset {\frak h}.
\end{align*}
\vspace{0.05in}\\
\noindent {\it Proof.}\quad It follows immediately from Lemma 9.2.
\hfill $\square$
\vspace{0.1in}\\
\noindent {\bf Remark\ 9.4.}\quad We remark that the relation
$$[{\frak p}^J,{\frak p}^J]\subset {\frak k}^J$$
does not hold.
\smallskip

Now we set
\begin{align*}
 Z_{ab}^0:&=-\sqrt{-1} D_{ab}^0,\\
 Y_{pq}^{\pm}:&={\frac 12}(D_{pq}\pm \sqrt{-1}{\hat D}_{pq}),\\
 Z_{ij}^{+}:&=-S_{ij},\\
 Z_{ij}^{-}:&=-\sqrt{-1}T_{ij},\\
 X_{ij}^{\pm}:&={\frac 12}(A_{ij}\pm \sqrt{-1} B_{ij}).
\end{align*}
\vspace{0.1in}\\
\noindent {\bf Lemma 9.5.} \quad {\it We have the following
commutation relation :}
\begin{align*}
 [Z_{ab}^0,Z_{cd}^0]&=[Z_{ab}^0,Y^{\pm}_{pq}]=
 [Z_{ab}^0,Z_{ij}^{\pm}]=[Z_{ab}^0,X_{ij}^{\pm}] = 0,\\
 [Y_{pq}^+,Y_{rs}^+]&=0,\ \ \ [Y_{pq}^+,Y_{rs}^-]=\d_{qs}Z_{pr}^0,\\
 [Y_{pq}^+,Z_{ij}^+]&=-\d_{qi}Y_{pj}^+ + \d_{qj}Y_{pi}^+,\\
 [Y_{pq}^+,Z_{ij}^-]&=-\d_{qi}Y_{pj}^+ - \d_{qj}Y_{pi}^+,\\
 [Y_{pq}^+,X_{ij}^+]&=0,\\
 [Y_{pq}^+,X_{ij}^-]&=\d_{qi}Y_{pj}^- + \d_{qj}Y_{pi}^-,\\
 [Y_{pq}^-,Y_{rs}^-]&=0,\\
 [Y_{pq}^-,Z_{ij}^+]&=-\d_{qi}Y_{pj}^- +\d_{qj}Y_{pi}^-,\\
 [Y_{pq}^-,Z_{ij}^-]&=\d_{qi}Y_{pj}^- + \d_{qj}Y_{pi}^-,\\
 [Y_{pq}^-,X_{ij}^+]&=\d_{qi}Y_{pj}^+ + \d_{qj}Y_{pi}^+,\\
 [Y_{pq}^-,X_{ij}^-]&= 0,\\
 \end{align*}

\begin{align*}
 [Z_{ij}^+,Z_{kl}^+]&= \d_{ik}Z_{jl}^+ -\d_{il}Z_{jk}^+ -\d_{jk}
Z_{il}^+ + \d_{jl}Z_{ik}^+,\\
 [Z_{ij}^+,Z_{kl}^-]&= \d_{ik}Z_{jl}^- -\d_{il}Z_{jk}^- +\d_{jk}
Z_{il}^- -\d_{jl}Z_{ik}^-,\\
 [Z_{ij}^+,X_{kl}^{\pm}]&=\d_{ik}X_{jl}^{\pm}-\d_{jk}X_{il}^{\pm}
+\d_{il}X_{jk}^{\pm}-\d_{jl}X_{ik}^{\pm},\\
 [Z_{ij}^-,Z_{kl}^-]&=-\d_{ik}Z_{jl}^+ -\d_{il}Z_{jk}^+ -\d_{jk}
Z_{il}^+ -\d_{jl}Z_{ik}^+,\\
 [Z_{ij}^-,X_{ij}^+]&=\d_{ik}X_{jl}^+ +\d_{il}X_{jk}^+ +\d_{jk}
X_{il}^+ +\d_{jl}X_{ik}^+,\\
 [Z_{ij}^-,X_{ij}^-]&=-\d_{ik}X_{jl}^- -\d_{il}X_{jk}^- -\d_{jk}
X_{il}^- -\d_{jl}X_{ik}^-,\\
 [X_{ij}^+,X_{kl}^+]&=[X_{ij}^-,X_{kl}^-]=0,\\
 [X_{ij}^+,X_{kl}^-]&=-{\frac 12}(\d_{ik}Z_{jl}^+
+\d_{il}Z_{jk}^+
+\d_{jk}Z_{il}^+ +\d_{jl}Z_{ik}^+) \\
&  +{{\sqrt{-1}}\over 2}(\d_{ik}Z_{jl}^- +\d_{il}Z_{jk}^- +\d_{jk}Z_{il}^-
+\d_{jl}Z_{ik}^-).
\end{align*}
\vspace{0.05in}\\
\noindent {\it Proof.}\quad It follows from Lemma 9.2. \hfill
$\square$
\vspace{0.1in}\\
\noindent {\bf Corollary\ 9.6.}\quad {\it The set
$$\left\{ Z_{ab}^0,\,Z_{ij}^{+},
\, Z_{kl}^{-}\,\bigg|
\ 1\leq a\leq b\leq m,\ 1\leq i< j\leq n,\ 1\leq k < l\leq n\ \right\}$$
form a basis of the
complexification ${\frak k}_{\mathbb{C}}^J$ of ${\frak k}^J$ and the set
$$\left\{ X_{ij}^{\pm},\,Y_{pq}^{\pm}\,\bigg| \ 1\leq i\leq j\leq n,\
1\leq p\leq m,\ 1\leq q\leq n\ \right\}$$
form a basis of ${\frak p}^J_{\mathbb{C}}.$
And $\left\{ X_{ij}^+,\, Y_{pq}^+\, \bigg| \ 1\leq i\leq j\leq n,\ 1\leq p
\leq m,\ 1\leq q\leq n\ \right\}$ form a basis of ${\frak p}^J_+$ and
$\left\{ X_{ij}^-,\,Y_{pq}^-\,\bigg| \ 1\leq i\leq j\leq n,\ 1\leq p\leq m,
\ 1\leq q\leq n\ \right\}$ form a basis of ${\frak p}_{-}^J.$ Both
${\frak p}_+^J$ and ${\frak p}_{-}^J$ are all abelian subalgebras of
${\frak g}_{\mathbb{C}}^J.$ We have the relation
$$[{\frak k}_{\mathbb{C}}^J,\,{\frak p}_+^J]\subset {\frak p}_+^J,\ \ \ \ \
[{\frak k}_{\mathbb{C}}^J,\,{\frak p}_{-}^J]\subset {\frak p}_{-}^J.$$
${\frak g}_{\mathbb{C}}$ is a subalgebra of ${\frak g}_{\mathbb{C}}^J$ and
${\frak h}_{\mathbb{C}}$,
the complexification of ${\frak h},$ is an ideal of}
${\frak g}_{\mathbb{C}}^J.$
\vspace{0.05in}\\
\noindent {\it Proof.}\quad It follows immediately from Lemma 9.5.
\hfill $\square$
\vspace{0.2in}\\
\noindent{\bf 9.3. Jacobi Forms }
\vspace{0.1in}\\
\indent Let $\rho$ be a rational representation of $GL(n,\mathbb{C})$ on a finite dimensional
complex vector space $V_{\rho}.$ Let ${\mathcal M}\in \mathbb R^{(m,m)}$ be a symmetric
half-integral semi-positive definite matrix of degree $m$.
Let $C^{\infty}(H_{n,m},V_{\rho})$ be the algebra of all
$C^{\infty}$ functions on $H_{n,m}$
with values in $V_{\rho}.$ For $f\in C^{\infty}(H_{n,m},
V_{\rho}),$
we define
\begin{align}
  & (f|_{\rho,{\mathcal M}}[(M,(\lambda,\mu,\kappa))])(Z,W) \notag \\
:=
\,& e^{-2\pi i\sigma({\mathcal M}[W+\lambda Z+\mu](CZ+D)^{-1}C)}
\times
e^{2\pi i\sigma({\mathcal M}(\lambda Z^t\!\lambda+2\lambda^t\!W+(\kappa+
\mu^t\!\lambda)))} \\
&\times\rho(CZ+D)^{-1}f(M<Z>,(W+\lambda Z+\mu)(CZ+D)^{-1}),\notag
\end{align}
where $M=\left(\begin{matrix} A&B\\ C&D\end{matrix}\right)\in Sp(n,\mathbb R),\
(\lambda,\mu,\kappa)\in H_{\mathbb R}^{(n,m)}$
and $(Z,W)\in H_{n,m}.$
\vspace{0.1in}\\
\noindent {\bf Definition\ 9.7.}\quad Let $\rho$ and $\mathcal M$ be as above. Let
$$H_{\mathbb Z}^{(n,m)}:= \{ (\lambda,\mu,\kappa)\in H_{\mathbb R}^{(n,m)}\, \vert
\, \lambda,\mu\in \mathbb Z^{(m,n)},\ \kappa\in \mathbb Z^{(m,m)}\,\ \}.$$
A {\it Jacobi\ form} of index $\mathcal M$ with respect to $\rho$
on $\Gamma_n$ is a holomorphic
function $f\in C^{\infty}(H_{n,m},V_{\rho})$ satisfying the
following conditions (A) and (B):

\smallskip

\noindent (A) \,\ $f|_{\rho,{\mathcal M}}[\tilde{\gamma}] = f$ for
all $\tilde{\gamma}\in \Gamma^J_n := \Gamma_n \ltimes H_{\mathbb
Z}^{(n,m)}$.

\smallskip

\noindent (B) \,\ $f$ has a Fourier expansion of the following
form :
$$f(Z,W) = \sum\limits_{T\ge0\atop \text {half-integral}}
\sum\limits_{R\in \mathbb Z^{(n,m)}}
c(T,R)\cdot e^{{2\pi i}\,\sigma(TZ)}\cdot
e^{2\pi i\sigma(RW)}$$
with $c(T,R)\ne 0$ only if $\left(\begin{matrix}
T & \frac 12R\\ \frac 12\,^t\!R&{\mathcal M}\end{matrix}\right)
\ge 0$.

\medskip

\indent If $n\geq 2,$ the condition (B) is superfluous by
K{\" o}cher principle\,(\,cf.\,
\cite{Zi} Lemma 1.6). We denote by $J_{\rho,\mathcal M}(\Gamma_n)$
the vector space of all
Jacobi forms of index $\mathcal{M}$
with respect to $\rho$ on $\Gamma_n$.
Ziegler\,(\,cf.\,\cite{Zi} Theorem 1.8 or \cite{E-Z} Theorem 1.1\,)
proves that the vector space $J_{\rho,\mathcal {M}}(\Gamma_n)$ is finite dimensional.
For more results on Jacobi forms with $n>1$ and $m>1$, we refer to
\cite{Kr},\,\cite{Y7}-\cite{Y11} and \cite{Zi}.
\vspace{0.1in}\\
 \noindent {\bf Definition\ 9.8.}\quad A Jacobi form $f\in J_{\rho,\mathcal {M}}(\Gamma)$ is said to be
a {\it cusp}\,(\,or {\it cuspidal}\,) form if
$\begin{pmatrix} T & {\frac 12}R\\ {\frac 12}\,^t\!R & \mathcal {M}\end{pmatrix} > 0$ for
any $T,\,R$ with $c(T,R)\ne 0.$ A Jacobi form $f\in J_{\rho,\mathcal{M}}(\Gamma)$ is
said to be {\it singular} if it admits a Fourier expansion such that a
Fourier coefficient $c(T,R)$ vanishes unless $\text{det}\begin{pmatrix}
T &{\frac 12}R\\ {\frac 12}\,^t\!R & \mathcal {M}\end{pmatrix}=0.$

\medskip
\noindent {\bf Example\ 9.9.}\quad Let $S\in \mathbb Z^{(2k,2k)}$
be a symmetric, positive definite unimodular even integral matrix
and $c\in \mathbb Z^{(2k,m)}.$ We define the theta series
$$\vartheta_{S,c}^{(g)}(Z,W):=\sum_{\lambda \in \mathbb Z^{(2k,n)}}
e^{\pi i\{ \sigma (S \lambda Z\,^t \lambda)+2\sigma (^t\!cS \lambda,^t\!W)\} },\ \ \
Z\in H_n,\ W\in \mathbb C^{(m,n)}.$$
\noindent We put $\mathcal{M}:={\frac 12}\,^t\!cSc.$ We assume that $2k<g+\text{rank}\,(\mathcal{M}).$
Then it is easy to see that $\vartheta_{S,c}^{(g)}$ is a singular form in
$J_{k,\mathcal{M}}(\Gamma_g)$\,(cf.\,\cite{Zi} p.\,212).
\vspace{0.2in}\\
\noindent{\bf 9.4. Characterization of Jacobi Forms as Functions
on the Jacobi Group  $G^J$}
\vspace{0.1in}\\
\indent In this section, we lift a Jacobi form $f\in J_{\rho,\mathcal M}(\Gamma_n)$ to a
smooth function $\Phi_f$ on the Jacobi group $G^J$ and characterize the
lifted function $\Phi_f$ on $G^J.$

\smallskip

\indent We recall that for given $\rho$ and $\mathcal M$, the canonical automorphic factor
$J_{\mathcal M,\rho}:G^J\times H_{n,m} \longrightarrow GL(V_{\rho})$ is given by
\begin{align*}
J_{\mathcal M,\rho}(g,(Z,W))&=e^{-2\pi i\sigma (\mathcal M[W+\lambda Z+\mu](CZ+D)^{-1}C)}\\
& \times e^{2\pi i\sigma(\mathcal M(\lambda Z\,^t\lambda +2\lambda \,^tW+\kappa+
\mu\,^t \lambda))}\rho(CZ+D)^{-1},
\end{align*}
where $g=(M,(\lambda,\mu;\kappa))\in G^J$ with
$M=\begin{pmatrix} A & B\\ C & D\end{pmatrix}
\in G.$ It is easy to see that the automorphic factor
$J_{\mathcal M,\rho}$ satisfies the cocycle condition:
\begin{equation}
J_{\mathcal M,\rho}(g_1g_2,(Z,W))=J_{\mathcal M,\rho}(g_2,(Z,W))\,J_{\mathcal M,\rho}(g_1,
g_2\cdot (Z,W))
\end{equation}
for all $g_1,\,g_2\in G^J$ and $(Z,W)\in H_{n,m}.$

\smallskip

\indent Since the space $H_{n,m}$ is diffeomorphic to the homogeneous space
$G^J/K^J$, we may lift a function $f$ on $H_{n,m}$ with values in $V_{\rho}$
to a function $\Phi_f$ on $G^J$ with values in $V_{\rho}$ in the following
way. We define the lifting
\begin{equation}
\varphi_{\rho,\mathcal M}\,:\,{\mathcal F}(H_{n,m},V_{\rho}) \longrightarrow
{\mathcal F}(G^J,V_{\rho}),\ \ \ \varphi_{\rho,\mathcal M}(f):=\Phi_f
\end{equation}
by
\begin{align*}
\Phi_f(g):&=(f|_{\rho,\mathcal M}[g])(iE_n,0)\\
            &=J_{\mathcal M,\rho}(g,(iE_n,0))\,f(g\cdot (iE_n,0)),
\end{align*}
where $g\in G^J$ and ${\mathcal F}(H_{n,m},V_{\rho})\,(\text{resp.}
\ {\mathcal F}(G^J,
V_{\rho}))$ denotes the vector space consisting of functions on $H_{n,m}\,
(\text{resp.}\ G^J)$ with values in $V_{\rho}.$

\smallskip

\indent For brevity, we set $\Gamma:=\Gamma_n=Sp(n,\mathbb Z)$
and $\Gamma^J=\Gamma \ltimes H_{\mathbb Z}^{(n,m)}.$
We let ${\mathcal F}_{\rho,\mathcal M}^{\Gamma}$ be the space of all functions $f$ on
$H_{n,m}$ with values in $V_{\rho}$ satisfying the transformation formula
\begin{equation}
f|_{\rho,\mathcal M}[\gamma]=f\ \ \ \ \text{for\ all}\ \gamma\in \Gamma^J.
\end{equation}
And we let ${\mathcal F}_{\rho,\mathcal M}^{\Gamma}(G^J)$ be the space of functions
$\Phi:G^J \longrightarrow V_{\rho}$ on $G^J$ with values in $V_{\rho}$ satisfying the
following conditions (9.32) and (9.33):

\smallskip

\begin{equation}
\Phi(\gamma g)= \Phi(g) \qquad  \text{for all} \quad \gamma \in \Gamma^J
\quad  \text{and} \quad g\in G^J.
\end{equation}
\begin{equation}
\Phi(g \,\ r(k,\kappa))=e^{2\pi i\sigma(\mathcal M\kappa)}\rho(k)^{-1}\Phi(g),
\qquad \forall\ \ r(k,\kappa):=[k,(0,0;\kappa)]\in K^J.
\end{equation}
\vspace{0.1in}\\
\noindent {\bf Lemma\ 9.10.}\quad {\it The space ${\mathcal F}_{\rho,\mathcal M}^{\Gamma}$ is isomorphic to
the space ${\mathcal F}_{\rho,\mathcal M}^{\Gamma}(G^J)$ via the lifting} $\varphi_{\rho,\mathcal M}.$
\vspace{0.05in}\\
\noindent {\it Proof.}\quad Let $f\in {\mathcal F}_{\rho,\mathcal M}^{\Gamma}.$ If $\gamma \in \Gamma^J,\ g\in
G^J$ and $r(k,\kappa)=[k,(0,0;\kappa)]\in K^J,$ then we have
\begin{align*}
\Phi_f(\gamma g)&=(f|_{\rho,\mathcal M}[\gamma g])(iE_n,0)\\
   &=((f|_{\rho,\mathcal M}[\gamma])|_{\rho,\mathcal M}[g])(iE_n,0)\\
   &= (f|_{\rho,\mathcal M}[g])(iE_n,0)\qquad \qquad\qquad\qquad(\,\text{since}\ f\in {\mathcal F}_{\rho,\mathcal M}^{\Gamma}\,)
   \\
   &= \Phi_f(g)
\end{align*}
and
\begin{align*}
\Phi_f(g\,r(k,\kappa))&= J_{\mathcal M,\rho}(g\,r(k,\kappa),(iE_n,0))\,
          f(g\,r(k,\kappa)\cdot (iE_n,0))\\
          &=J_{\mathcal M,\rho}(r(k,\kappa),(iE_n,0))\,J_{\mathcal M,\rho}(g,(iE_n,0))\,
             f(g\cdot (iE_n,0))\\
          &= e^{2\pi i\sigma (\mathcal M \kappa)}\rho(k)^{-1}\Phi_f(g).
\end{align*}
Here we identified $k=\begin{pmatrix} A & -B\\ B & A\end{pmatrix} \in K$ with
$A+iB\in U(n).$

\medskip

\indent Conversely, if $\Phi\in {\mathcal F}_{\rho,\mathcal M}^{\Gamma}(G^J),\ G^J$ acting on
$H_{n,m}$ transitively, we may define a function $f_{\Phi}$ on $H_{n,m}$ by
\begin{equation}
f_{\Phi}(g\cdot (iE_n,0)):=J_{\mathcal M,\rho}(g,(iE_n,0))^{-1}\Phi(g).
\end{equation}
Let $\gamma\in \Gamma^J$ and $(Z,W)=g\cdot (iE_n,0)$ for some $g\in G^J.$
Then using the cocycle condition (9.29), we have
\begin{align*}
\ \ \ (f_{\Phi}|_{\rho,\mathcal M}[\gamma])(Z,W)&= J_{\mathcal M,\rho}(\gamma,(Z,W))f_{\Phi}
(\gamma\cdot (Z,W))\\
&= J_{\mathcal M,\rho}(\gamma,g\cdot (iE_n,0))\,f_{\Phi}(\gamma g\cdot (iE_n,0))\\
&= J_{\mathcal M,\rho}(\gamma,g\cdot (iE_n,0))J_{\mathcal M,\rho}(\gamma g,(iE_n,0))^{-1}
\Phi(\gamma g)\\
&=J_{\mathcal M,\rho}(\gamma,g\cdot (iE_n,0))\,J_{\mathcal M,\rho}
(\gamma,g\cdot (iE_n,0))^{-1}\\
& J_{\mathcal M,\rho}(g,(iE_n,0))^{-1}\Phi(g)\\
&=J_{\mathcal M,\rho}(g,(iE_n,0))^{-1}\Phi(g)\\
&=f_{\Phi}(g\cdot (iE_n,0))=f_{\Phi}(Z,W).
\end{align*}
This completes the proof. \hfill $\square$

\medskip

Now we have the following two algebraic representations
$T_{\rho,\mathcal M}$ and
$\dot{T}_{\rho,\mathcal M}$ of $G^J$ defined by
\begin{equation}
T_{\rho,\mathcal M}(g)f:=f|_{\rho,\mathcal M}[g^{-1}],\ \ \ g\in G^J,\ f\in
{\mathcal F}_{\rho,\mathcal M}^{\Gamma}
\end{equation}
and
\begin{equation}
\dot{T}_{\rho,\mathcal M}(g)\Phi(g'):=\Phi(g^{-1}g'),\ \ \ g,\,g'\in G^J,\ \
\Phi\in {\mathcal F}_{\rho,\mathcal M}^{\Gamma}(G^J).
\end{equation}
Then it is easy to see that these two models $T_{\rho,\mathcal M}$ and
$\dot{T}_{\rho,\mathcal M}$ are intertwined by the lifting $\varphi_{\rho,\mathcal M}.$
\vspace{0.1in}\\
\noindent {\bf Proposition\ 9.11.} \quad The vector space $J_{\rho,\mathcal M}(\Gamma_n)$ is
isomorphic to the space $A_{\rho,\mathcal M}(\Gamma^J)$ of smooth functions $\Phi$ on
$G^J$ with values in $V_{\rho}$ satisfying the following conditions:

\smallskip

\indent \quad (1a)\,\ $\Phi(\gamma g)=\Phi(g)$\ \ for\ all\ $\gamma \in \Gamma^J.$\\
\indent \quad (1b)\,\ $\Phi(g  r(k,\kappa))=e^{2\pi i\sigma
(\mathcal M\kappa)}\rho(k)^{-1}
\Phi(g)$ \,\ for\ all\ $g \in G^J,\ r(k,\kappa)\in K^J.$\\
\indent \quad (2)\,\ $X_{ij}^{-}\Phi=Y_{ij}^{-}\Phi=0,\ \ \ 1\leq i,j\leq n.$\\
\indent \quad (3)\,\ For all $M\in Sp(n,\mathbb R),$ the function
$\psi:G^J \longrightarrow V_{\rho}$ defined by
$$\psi(g):=\rho(Y^{-{1\over 2}})\,\Phi(Mg),\ \ \ g\in G^J$$
is bounded in the domain $Y \geq Y_0.$ Here $g\cdot (iE_n,0)=(Z,W)$ with
$Z=X+iY,\ Y>0.$
\vspace{0.1in}\\
\noindent {\bf Corollary\ 9.12.}\quad $J_{\rho,\mathcal M}^{\text{cusp}}(\Gamma_n)$
is isomorphic to the
subspace $A_{\rho,\mathcal M}^0(\Gamma^J)$ of $A_{\rho,\mathcal M}(\Gamma^J)$ with the condition
(3') the function $g \longmapsto \Phi(g)$ is bounded.
\vspace{0.2in}\\
\noindent{\bf 9.5.\ Unitary Representations of the Jacobi Group $G^J$}
\vspace{0.1in}\\
\indent In this section, we review some results of Takase\,(cf. \cite{T1}-\cite{T3}) on the
unitary representations of the Jacobi group $G^J$. We follow the notations
in the previous sections.

\smallskip

First we observe that $G^J$ is not reductive because the center of $G^J$
is given by
$${\mathcal Z}=\left\{ [E_{2n},(0,0;\kappa)]\in G^J\ \bigg|
\kappa=\,^t\kappa \in \mathbb R^{(m,m)}\ \right\}
\cong \text{Sym}^2(\mathbb R^m).$$
Let $d_K(k)$ be a normalized Haar measure on $K$ so that
$\int_K d_K(k)=1$ and $d_{\mathcal Z}
(\kappa)=\prod\limits_{i\leq j}d \kappa_{ij}$ a Haar measure on
${\mathcal Z}$. We let $d_{K^J}=d_K\times d_{\mathcal Z}$ be the product measure on
$K^J=K\times {\mathcal Z}.$
The Haar measure $d_{G^J}$ on $G^J$ is normalized so that
$$\int_{G^J}f(g)d_{G^J}(g)=\int_{G^J/K^J}\left( \int_{K^J}f(gh)d_{K^J}(h)
\right) d_{G^J/K^J}({\dot {g}})$$
for all $f\in C_c(G^J).$

\smallskip

From now on, we will fix a real positive definite symmetric matrix
$S\in \text{Sym}^2(\mathbb R^m)$ of degree $m$. For any fixed $Z=X+iY\in H_n,$ we
define a measure $\nu_{S,Z}$ on $\mathbb C^{(m,n)}$ by
\begin{equation}
d\nu_{S,Z}(W)=\left(\text{det}\, 2S\right)^n \left( \text{det}
\,Y\right)^{-m}
\kappa_S(Z,W)\,dUdV,
\end{equation}
where $W=U+iV\in \mathbb C^{(m,n)}$ with $U,V\in \mathbb R^{(m,n)}$ and
\begin{equation}
 \kappa_S(Z,W)=e^{-4\pi \sigma(\,^tVSVY^{-1})}.
 \end{equation}
\indent Let $H_{S,Z}$ be the complex Hilbert space consisting of all $\mathbb C$-valued
holomorphic functions $\varphi$ on $\mathbb C^{(m,n)}$ such that
$\int_{\mathbb C^{(m,n)}}\,|\varphi(W)|^2d\nu_{S,Z} < +\infty.$ The inner product on
$H_{S,Z}$ is given by
$$(\varphi,\psi)=\int_{\mathbb C^{(m,n)}}\varphi(W)\overline{\psi(W)}
d\nu_{S,Z}(W),\ \ \ \ \varphi,\,\psi\in H_{S,Z}.$$
We put
$$\eta_S=J_{S,\delta}^{-1}\ \ \ (\text{see\ subsection}\ 9.4),$$
where $\delta$ denotes the trivial representation of $GL(n,\mathbb C)$. Now we
define a unitary representation $\Xi_{S,Z}$ of $H_{\mathbb R}^{(n,m)}$ by
\begin{equation}
\left( \Xi_{S,Z}(h)\varphi\right)(W)=\eta_S(h^{-1},(Z,W))^{-1}\cdot
\varphi(W-\lambda Z-\mu),
\end{equation}
where $h=(\lambda,\mu,\kappa)\in H_{\mathbb R}^{(n,m)}$ and $\varphi\in H_{S,Z}.$
It is easy to see that $(\Xi_{S,Z},H_{S,Z})$ is irreducible and
$\Xi_{S,Z}(0,0,\kappa)=e^{-2\pi i \sigma(S\kappa)}.$

\smallskip
Let $${\mathcal X}=\left\{ T\in \mathbb C^{(n,n)}\,\bigg| \ T=\,^tT,\ \ \text{Re}\,T>0\
\right\}$$
be a connected simply connected open subset of $\mathbb C^{(n,n)}$. Then there
exists uniquely a holomorphic function $\text{det}^{\frac 12}$ on
${\mathcal X}$
such that
\begin{align*}
&(1)\ \ \ \ \left( \text{det}^{\frac 12}\,
T\right)^2=\text{det}\,T\ \ \ \ \text{for\ all}\ T\in
T\in {\mathcal X},\\
&(2)\ \ \ \ \ \text{det}^{\frac 12}\,T=(\text{det}\,T)^{\frac 12}\ \ \
\ \text{for\ all}\ T\in {\mathcal X}\cap
\mathbb R^{(n,n)}.
\end{align*}
For any integer $k\in \mathbb Z,$ we set
\begin{equation}
\text{det}^{\frac k2}\,T=\left
( \text{det}^{\frac 12}\,T\right)^k,\ \ \ T\in {\mathcal X}.
\end{equation}
For any $g=(\sigma,h)\in G^J$ with $\sigma\in G$, we define an integral operator
$T_{S,Z}(g)$ from $H_{S,Z}$ to $H_{S,\sigma<\!Z\!>}$ by
\begin{equation}
\left(T_{S,Z}(g)\varphi\right)(W)=\eta_S(g^{-1},(Z,W))^{-1}
\varphi(W'),
\end{equation}
where $(Z',W')=g^{-1}\cdot (Z,W)$ and $\varphi\in H_{S,Z}.$ And for any
fixed $Z$ and $Z'$ in $H_n,$ we define a unitary mapping
\begin{equation}
U_{Z',Z}^S:\,H_{S,Z}\longrightarrow H_{S,Z'}
\end{equation}
by
$$\left( U_{Z',Z}^S\varphi\right)(W')=\gamma(Z',Z)^m\cdot
\int_{\mathbb C^{(m,n)}}\,\kappa_S((Z',W'),(Z,W))^{-1}\varphi(W)\,d\nu_{S,Z}(W),$$
where
$$\gamma(Z',Z)=\text{det}^{-{\frac 12}}
\left( {{Z'-{\bar Z}}\over {2i}}\right)
\cdot \text{det\,(Im}\,Z')^{\frac 14}\cdot \text{det\,(Im}\,Z)^{\frac 14}$$
and
$$\kappa_S((Z',W'),(Z,W))=e^{2\pi i \sigma(S[W'-{\overline W}]\cdot
(Z'-{\overline Z})^{-1}
)}.$$
For any $g=(\sigma,h)\in G^J,$ we define a unitary operator $T_S(g)$ of
$H_{S,iE_n}$ by
\begin{equation}
T_S(g)=T_{S,\sigma^{-1}<iE_n>}(g)\circ U_{\sigma^{-1}<iE_n>,iE_n}^S.
\end{equation}
We put, for any $\sigma_1,\sigma_2\in G$,
\begin{equation}
\beta(\sigma_1,\sigma_2)={{\gamma(\sigma_1^{-1}<iE_n>,iE_n)}\over {\gamma(\sigma_2^{-1}
\sigma_1^{-1}<iE_n>,\sigma_1^{-1}<iE_n>)}}\ \ .
\end{equation}
Then the function $\beta(\sigma_1,\sigma_2)$ satisfies the cocycle condition
$$\beta(\sigma_2,\sigma_3)\beta(\sigma_1\sigma_2,\sigma_3)^{-1}\beta(\sigma_1,\sigma_2\sigma_3)
\beta(\sigma_1,\sigma_2)^{-1}=1$$
for all $\sigma_1,\sigma_2,\sigma_3\in G.$ Thus $\beta(\sigma_1,\sigma_2)$ defines a group
extension $G\ltimes \mathbb C_1$ by $\mathbb C_1=\{ z\in \mathbb C\,|\ |z|=1\ \}.$ Precisely,
$G\ltimes \mathbb C_1$ is a topological group with multiplication
$$(\sigma_1,\epsilon_1)\cdot (\sigma_2,\epsilon_2)=(\sigma_1\sigma_2,\beta(\sigma_1,\sigma_2)\epsilon_1\epsilon_2)$$
for all $(\sigma_i,\epsilon_i)\in G\times \mathbb C_1\,(i=1,2).$ If we put
\begin{equation}
\epsilon(\sigma)={{\text{det}\,J(\sigma^{-1},iE_n)}\over
{|\text{det}\,J(\sigma^{-1},iE_n)|}}\ \ ,\ \
\ \
\sigma\in G,
\end{equation}
then we have the relation
$$\beta(\sigma_1,\sigma_2)^2=\epsilon(\sigma_1)\cdot \epsilon(\sigma_1\sigma_2)^{-1}\cdot \epsilon(\sigma_2),\ \ \
\sigma_1,\,\sigma_2\in G.$$
Therefore we have a closed normal subgroup
\begin{equation}
G_2=\left\{ (\sigma,\epsilon)\in G\ltimes \mathbb C_1\,\bigg| \ \epsilon^2=\epsilon(\sigma)^{-1}
\right\}
\end{equation}
of $G_2\ltimes \mathbb C_1$ which is a connected two-fold covering group of $G$.
Since $G_2$ acts on the Heisenberg group $H_{\mathbb R}^{(n,m)}$ via the
projection $p:G_2\longrightarrow G,$ we may put
$$G_2^J=G_2\ltimes H_{\mathbb R}^{(n,m)}.$$
Now we define the unitary representation $\omega_S$ of $G_2^J$ by
\begin{equation}
\omega_S(g)=\epsilon ^m\cdot T_S(\sigma,h),\ \ \ g=((\sigma,\epsilon),h)\in G_2^J.
\end{equation}
It is easy to see that $(\omega_S,H_{S,iE_n})$ is irreducible and the
restriction of $\omega_S$ to $G_2$ is the $m$-fold tensor product of the
Weil representation. $\omega_S$ is called the {\it Weil\ representation}
of the Jacobi group $G^J.$

\smallskip

We set
\begin{align*}
\ \ \ &p:G_2 \longrightarrow G,\ \ \ p(\sigma,\epsilon)=\sigma,\\
&p^J:G_2^J \longrightarrow G^J,\ \ \ p^J((\sigma,\epsilon),h)=(\sigma,h),\\
&q:G^J \longrightarrow G,\ \ \ q(\sigma,h)=\sigma,\\
&q^J:G_2^J\longrightarrow G_2,\ \ \ q^J((\sigma,\epsilon),h)=(\sigma,\epsilon).
\end{align*}
\vspace{0.1in}\\
\noindent {\bf Proposition\ 9.13.}\quad {\it Let $\chi_S$ be the
character of ${\mathcal Z}\cong \text{Sym}^2(\mathbb R^m)$ defined
by $\chi_S(\kappa)=e^{2\pi i\sigma(S\kappa)},\, \kappa\in
{\mathcal Z}.$ We denote by ${\hat {G}}_2^J({\bar {\chi}}_S)$ the
set of all equivalence classes of irreducible unitary
representations $\tau$ of $G_2^J$ such that
$\tau(\kappa)=\chi_S(\kappa)^{-1}$ for all $\kappa\in {\mathcal
Z}.$ We put ${\tilde {\pi}}=\pi\circ q^J\in {\hat {G}}_2^J$ for
any $\pi\in {\hat {G}}_2.$ The correspondence
$$\pi \longmapsto {\tilde {\pi}}\otimes \omega_S$$
is a bijection from ${\hat {G}}_2$ to ${\hat {G}}_2^J({\bar{\chi}}_S).$
And ${\tilde {\pi}}\otimes \omega_S$
is square-integrable modulo ${\mathcal Z}$
if and only if $\pi$ is square integrable.}
\vspace{0.05in}\\
\noindent {\it Proof.} \quad See \cite{T1}, Proposition 11.8.\hfill $\Box$
\vspace{0.1in}\\
\noindent {\bf Proposition\ 9.14.} {\it Let $m$ be even. We put ${\check{\pi}}
=\pi\circ q\in
{\hat {G}}^J$ for any $\pi\in {\hat G}.$ Then the correspondence
$$\pi \longmapsto {\check{\pi}}\otimes \omega_S$$
is a bijection of ${\hat {G}}$ to ${\hat {G}}^J.$ And ${\check{\pi}}\otimes
\omega_S$ is square integrable modulo $A$ if and only if $\pi$ is square
integrable.}
\vspace{0.05in}\\
\noindent {\it Proof.}\quad See \cite{T1}.\hfill$\square$

\smallskip

The above proposition was proved by Satake \cite{S2} or by Berndt \cite{Be1} in the
case $m=1.$

\smallskip

\indent Let $(\rho,V_{\rho})$ be an irreducible representation of $K=U(n)$ with
highest weight $l=(l_1,l_2,\cdots,l_n)\in \mathbb Z^n,\ l_1\geq\cdots\geq l_n
\geq 0.$ Then $\rho$ is extended to a rational representation of $GL(n,\mathbb C)$
which is also denoted by $\rho.$ The representation space $V_{\rho}$ of
$\rho$ has an hermitian inner product ( , ) such that
$(\rho(g)u,v)=(u,\rho(g^*)v)$ for all $g\in GL(n,\mathbb C),\ u,v\in V_{\rho},$
where $g^*=\,^t{\bar g}$. We let the mapping $J:G\times H_n \longrightarrow GL(n,\mathbb C)$
be the automorphic factor defined by
$$J(\sigma,Z)=CZ+D,\ \ \ \sigma=\begin{pmatrix} A & B\\ C & D
\end{pmatrix} \in G.$$
We define a unitary representation $\tau_l$ of $K$ by
\begin{equation}
\tau_l(k)=\rho(J(k,iE_n)),\ \ \ k\in K.
\end{equation}
We set $J_{\rho,S}=J_{S,\rho}^{-1}$\,(\,cf.\,subsection 9.4\,).
According to the definition, we have
$$J_{\rho,S}(g,(Z,W))=\eta_S(g,(Z,W))\,\rho(J(\sigma,Z))$$
for all $g=(\sigma,h)\in G^J$ and $(Z,W)\in H_{n,m}.$ For any $g=(\sigma,h)\in
G^J$ and $(Z,W)\in H_{n,m},$ we set
\begin{align*}
\overline{J_{\rho,S}(g,(Z,W))}&=\overline{\eta_S(g,(Z,W))}\,
\rho(\overline{J(\sigma,Z)}),\\
^tJ_{\rho,S}(g,(Z,W))&=\eta_S(g,(Z,W))\,\rho(\,^tJ(\sigma,Z)),\\
J_{\rho,S}(g,(Z,W))^*&=\,^t\overline{J_{\rho,S}(g,(Z,W))}.
\end{align*}
Then for all $g\in G^J,\ (Z,W)\in H_{n,m}$ and $u,v\in V_{\rho},$ we have
$$(J_{\rho,S}(g,(Z,W))u,v)=(u,J_{\rho,S}(g,(Z,W))^*v)$$
We denote by $E(\rho,S)$ the Hilbert space consisting of $V_{\rho}$-valued
measurable functions $\varphi$ on $H_{n,m}$ such that
$$|\varphi|^2=\int_{H_{n,m}}\left(\rho(Im\,Z)\,\varphi(Z,W),\varphi(Z,W)
\right)\,\kappa_S(Z,W)\,d(Z,W),$$
where
$$d(Z,W)=(\text{det}\,Y)^{-(m+n+1)}dXdYdUdV,\ \ Z=X+iY,\ W=U+iV$$
denotes a $G^J$-invariant volume element on $H_{n,m}.$
The induced representation $\text{Ind}_{K^J}^{G^J}
(\rho\otimes {\bar {\chi}}_S)$
is realized on $E(\rho,S)$ as follows: For any $g\in G^J$ and $\varphi\in
E(\rho,S),$ we have
$$\left( \text{Ind}_{K^J}^{G^J}(\rho\otimes {\bar {\chi}}_S)(g)\varphi
\right)(Z,W)=J_{\rho,S}(g^{-1},(Z,W))^{-1}\varphi(g^{-1}\cdot (Z,W)).$$
We recall that $\chi_S$ is the unitary character of $A$ defined by
$\chi_S(\kappa)=e^{2\pi i\sigma(S\kappa)},\ \kappa\in {\mathcal Z}.$ Let $H(\rho,S)$ be
the subspace of $E(\rho,S)$ consisting of $\varphi\in E(\rho,S)$ which is
holomorphic on $H_{n,m}.$ Then $H(\rho,S)$ is a closed $G^J$-invariant
subspace of $E(\rho,S).$ Let $\pi^{\rho,S}$ be the restriction of the
induced representation $\text{Ind}_{K^J}^{G^J}
(\rho\otimes {\bar {\chi}}_S)$ to
$H(\rho,S).$

\smallskip

\indent Takase\,(cf. \cite{T2},\,Theorem  1.1) proved the following
\vspace{0.1in}\\
\noindent {\bf Theorem\ 9.15.}\quad {\it Suppose $l_n>n+{\frac
m2}.$ Then $H(\rho,S)\neq 0$ and $\pi^{\rho,S}$ is an irreducible
unitary representation of $G^J$ which is square integrable modulo
${\mathcal Z}$. The multiplicity of $\rho_l$ in $\pi^{\rho,S}|_K$
is equal to one.}

\medskip

We put
$$K_2=p^{-1}(K)=\left\{ (k,\epsilon)\in K\times \mathbb C_1\,\bigg|\ \epsilon^2=
\text{det}\,J(k,iE_n)\,\right\}.$$
The Lie algebra ${\frak k}$ of $K_2$ and its Cartan algebra are given by
$${\frak k}=\left\{\begin{pmatrix} A&-B\\ B&A
\end{pmatrix} \in \mathbb R^{(2n,2n)}\,\bigg|\
A+\,^tA=0,\ B=\,^tB\right\}$$
and
$${\frak h}=\left\{\begin{pmatrix} 0 & -C\\ C & 0
\end{pmatrix} \in \mathbb R^{(2n,2n)}\,\bigg|\
C=\text{diag}\,(c_1,c_2,\cdots,c_n)\ \right\}.$$
Here $\text{diag}\,(c_1,c_2,\cdots,c_n)$
denotes the diagonal matrix of degree $n$.
We define $\lambda_j\in {\frak h}_{\mathbb C}^*$ by $\lambda_j\begin{pmatrix}
0 & -C\\ C & 0\end{pmatrix}:
=\sqrt{-1}c_j.$ We put
$$M^+=\left\{\sum_{j=1}^n m_j\lambda_j\bigg|\ m_j\in {\frac 12}\mathbb Z,\ m_1\geq
\cdots\geq m_n,\ m_i-m_j\in \mathbb Z\
\text{for\ all}\ i,j\right\}.$$
We take an element $\lambda=\sum_{j=1}^n m_j\lambda_j\in M^+.$ Let $\rho$ be an
irreducible representation of $K$ with highest weight $l=(l_1,\cdots,l_n)
\in \mathbb Z^n,$ where $l_j=m_j-m_n\,(1\leq j\leq n-1).$ Let $\rho_{[\lambda]}$ be the
irreducible representation of $K_2$ defined by
\begin{equation}
\rho_{[\lambda]}(k,\epsilon)=\epsilon^{2m_n}\cdot \rho(J(k,iE_n)),\ \ (k,\epsilon)\in K_2.
\end{equation}
Then $\rho_{[\lambda]}$ is the irreducible representation of $K_2$ with highest
weight $\lambda=(m_1,\cdots,m_n)$ and $\lambda \longmapsto\rho_{[\lambda]}$ is a bijection from
$M^+$ to ${\hat K}_2,$ the unitary dual of $K_2.$

\smallskip

The following proposition is a special case of \cite{KaV1}, Theorem 7.2.
\vspace{0.1in}\\
\noindent {\bf Proposition\ 9.16.}\quad {\it We have an
irreducible decomposition
$$\omega_S\bigg|_{K_2}=\oplus_{\lambda}\,m(\lambda)\rho_{[\lambda]},$$
where $\lambda$ runs over
\begin{align*}
\ \ \ &\lambda=\sum_{j=1}^{\nu}l_j\lambda_j+{\frac m2}\sum_{j=1}^n\lambda_j\in M^+\
\ (\nu=\text{min}\,\{m,n\}),\\
\ \ \ &\lambda_j\in \mathbb Z\ \,\text{such\ that}\ l_1\geq\l_2\cdots\geq l_{\nu}\geq 0
\end{align*}
and the multiplicity $m(\lambda)$ is given by
$$m(\lambda)=\prod_{1\leq i< j \leq m}\left( 1+{{l_i-l_j}\over {j-i}}\right),$$
where $l_j=0$ if $j>\nu.$
Let ${\hat G}_{2,d}$ be the set of all the unitary equivalence classes of
square integrable irreducible unitary representations of $G_2.$ The
correspondence
$$\pi\longmapsto \text{Harish}\!-\!\text{Chandra\ parameter\ of}\ \pi$$
is a bijection from ${\hat G}_{2,d}$ to $\Lambda^+,$ where
$$\Lambda^+=\left\{\sum_{j=1}^n m_j\lambda_j\in M^+\ \bigg|\ m_1>\cdots>m_n,\
m_i-m_j\neq 0\ \text{for\ all}\ i,j,\ i\neq j\right\}.$$}
See \cite{Wa}, Theorem 10.2.4.1 for the details.
\smallskip

We take an element $\lambda=\sum_{j=1}^n m_j\lambda_j\in M^+.$ Let $\pi^{\lambda}\in
{\hat G}_{2,d}$ be the representation corresponding to the Harish-Chandra
parameter
$$\sum_{j=1}^n(m_j-j)\lambda_j\in \Lambda^+.$$
The representation $\pi^{\lambda}$ is realized  as follows\,(see \cite{Kn}, Theorem
6.6)\,: Let $(\rho,V_{\rho})$ be the irreducible representation of $K$ with
highest weight $l=(l_1,\cdots,l_n),\
l_i=m_i-m_n\,(\,1\leq j\leq n\,).$
Let $H^{\lambda}$ be a complex Hilbert space consisting of the $V_{\rho}$-valued
holomorphic functions $\varphi$ on $H_n$ such that
$$|\varphi|^2=\int_{H_n}\left(\rho(\text{Im}\,Z)\,
\varphi(Z),\varphi(Z)\right)\cdot
(\text{det\,Im}\,Z)^{m_n}\,dZ < +\infty,$$
where $dZ$ is the usual $G_2$-invariant measure on $H_n$. Then $\pi^{\lambda}$
is defined by
$$\left(\pi^{\lambda}(g)\varphi\right)(Z)=J_{\lambda}(g^{-1},Z)^{-1}\varphi(g^{-1}
<Z>)$$
for all $g=(\sigma,\epsilon)\in G_2$ and $\varphi\in H^{\lambda}.$ Here
$$J_{\lambda}(g,Z)=\rho(J(\sigma,Z))\cdot J_{{\frac 12}}(g,Z)^{m_n},$$
where
$$J_{{\frac 12}}(g,Z)={{\gamma(\sigma<Z>,\sigma<iE_n>)}\over {\gamma(Z,iE_n)}}\cdot
\beta(\sigma,\sigma^{-1})\cdot \epsilon \cdot |det\,J(\sigma,Z)|^{{\frac 12}}\ .$$
\vspace{0.1in}\\
\noindent {\bf Proposition\ 9.17.} {\it Suppose $l_n>n+{\frac m2}.$ We put
$\lambda=\sum_{j=1}^n(l_j-{\frac m2})\lambda_j\in M^+.$ Then $\pi^{\rho,S}$ is an
irreducible unitary representation of $G^J$ and we have a unitary
equivalence
$$(\pi^{\lambda}\circ q^J)\otimes \omega_S\longrightarrow \pi^{\rho,S}\circ p^J$$
via the intertwining operator $\Lambda_{\rho,S}:H^{\lambda}\otimes H_{S,iE_n}
\longrightarrow H(\rho,S)$ defined by
$$\left(\Lambda_{\rho,S}(\varphi\otimes \psi)\right)(Z,W)=
(\,\text{det}\,2S\,)^n(\,\text{det\,Im}\,Z\,)^{-{\frac m4}}
\varphi(Z)(U_{Z,iE_n}^S\psi)
(W)$$
for all $\varphi\in H^{\lambda}$ and $\psi\in H_{S,iE_n}.$}
\vspace{0.2in}\\
\noindent{\bf 9.6.\ Duality Theorem for $G^J$}
\vspace{0.1in}\\
\indent In this subsection, we state the duality theorem for the Jacobi
group $G^J.$

\smallskip
\indent Let $E_{ij}$ denote a square matrix of degree $2n$ with
entry 1 where the $i$-th row and the $j$-th column meet, all other
entries being $0$. We put
$$H_i=E_{ii}-E_{n+i,n+i}\,(1\leq i\leq n),\quad
\mathfrak {h}=\sum_{i=1}^n\mathbb{C} H_i.$$ Then $\mathfrak {h}$
is a Cartan subalgebra of $\mathfrak {g} .$ Let $e_j:\mathfrak
{h}\longrightarrow \mathbb {C}\,(1\leq j\leq n)$
be the linear form on $\mathfrak {h}$ defined by
$$e_j(H_i)=\delta_{ij},$$
where $\delta_{ij}$ denotes the Kronecker delta symbol. The roots
of $\mathfrak {g} $ with respect to $\mathfrak {h}$ are given by
$$\pm 2e_i\ (1\leq i\leq n),\ \pm e_k\pm e_l\,(1\leq k<l
\leq n).$$ The set $\Phi^+$ of positive roots is given by
$$\Phi^+=\{ 2e_i\,(1\leq i\leq n),\ e_k+e_l\,(\,1\leq k<l\leq
n\,)\}.$$
\def\fn{{\frak n}}
Let $$\mathfrak {g} _{\alpha}=\left\{\,X\in\mathfrak {g}
\,\ \vert \,\ [H,X]=\alpha (H)X\ \text{for\ all}\ H\in \mathfrak
{h}\,\right\}$$ be the root space corresponding to a root $\alpha$ of
$\mathfrak {g} $ with respect to $\mathfrak {h}$. We put
$\fn=\sum_{\Phi^+}\mathfrak {g} _{\alpha}.$ We define
$$N^J=\left\{ [\text{exp}\,X,(0,\mu,0)]\in G^J\Big|\
X\in\fn\,\right\},$$ where $\text{exp}:\mathfrak {g}
\longrightarrow G$ denotes the exponential mapping from
$\mathfrak {g} $ to $G$. A subgroup $N^g$ of $G^J$ is said to be
{\it horosherical} if it is conjugate to $N^J$, that is,
$N^g=gN^Jg^{-1}$ for some $g\in G$. A horospherical subgroup $N^g$
is said to be {\it cuspidal} for $\Gamma^J=\Gamma_n\ltimes
H^{(n,m)}_{\mathbb Z}$ in $G^J$ if $(N^g \cap \Gamma^J)\backslash N^g$ is compact.
Let $L^2(\Gamma^J\backslash G^J,\rho)$ be the complex Hilbert space
consisting of all $\Gamma^J$-invariant $V_{\rho}$-valued measurable
functions $\Phi$ on $G^J$ such that $||\Phi||<\infty,$ where $||\
\, ||$ is the norm induced from the norm $|\ \,|$ on $E(\rho,{\mathcal
M})$ by the lifting from $H_{n,m}$ to $G^J$. We denote by
$L_0^2(\Gamma^J \backslash G^J,\rho)$ the subspace of $L^2(\Gamma^J\backslash G^J,\rho)$
consisting of functions $\varphi$ on $G^J$ such that $\varphi\in
L^2(\Gamma^J\backslash G^J,\rho)$  and $$\int_{N^g\cap \Gamma^J\backslash
N^g}\varphi(ng_0)dn=0$$ for any cuspidal subgroup $N^g$ of $G^J$
and almost all $g_0\in G^J$. Let $R$ be the right regular
representation of $G^J$ on $L_0^2(\Gamma^J\backslash G^J,\rho)$.

\smallskip

Now we state the duality theorem for the Jacobi group $G^J.$

\smallskip

\noindent {\bf Duality Theorem.}\quad {\it Let $\rho$ be an
irreducible representation of $K$ with highest weight
$l=(l_1,\cdots,l_n)\in \mathbb Z^n,\ l_1\geq l_2\geq\cdots\geq
l_n.$ Suppose $l_n>n+{\frac 12}$ and let ${\mathcal M}$ be a half
integrable positive definite symmetric matrix of degree $m$. Then
the multiplicity $m_{\rho,{\mathcal M}}$ of $\pi^{\rho,{\mathcal
M}}$ in the right regular representation $R$ of $G^J$ in
$L^2_0(\Gamma^J\backslash G^J,\rho)$ is equal to the dimension of
$J_{\rho,{\mathcal M}}^{\text{cusp}}(\Gamma_n)$, that is,}
$$m_{\rho,{\mathcal M}}=\text{dim}_{\mathbb{C}} \,
J_{\rho,{\mathcal M}}^{\text{cusp}}(\Gamma_n).$$ We may prove the
above theorem following the argument of \cite{B-B} in the case
$m=n=1.$ So we omit the detail of the proof.
\vspace{0.2in}\\
\noindent{\bf 9.7. Coadjoint Orbits for the Jacobi Group $G^J$}
\vspace{0.1in}\\
\indent  We observe that the Jacobi group $G^J$ is embedded in
$Sp(n+m,\mathbb R)$ via
\begin{equation}
(M,(\lambda,\mu,\kappa))\mapsto \begin{pmatrix} A & 0& B
&A{}^t\!\mu-B{}^t\!\lambda\\
                                    \lambda& E_m&\mu&\kappa\\
                    C&0&D&C{}^t\!\mu-D{}^t\!\lambda \\
            0&0&0&E_m \end{pmatrix},
            \end{equation}
where $(M,(\lambda,\mu,\kappa))\in G^J$ with $M=\begin{pmatrix} A
& B\\ C & D\end{pmatrix} \in Sp(n,\mathbb{R}).$ The Lie algebra
$\mathfrak {g} ^J$ of $G^J$ is given by
\begin{equation}
\mathfrak {g}
^J=\left\{ (X,(P,Q,R))\,|\ X\in \mathfrak {g} ,\ P,Q\in
\mathbb{R}^{(m,n)},\ R={}^t\!R\in \mathbb{R}^{(m,m)}\,\right\}
\end{equation}
 with the bracket
\begin{equation}
[(X_1,(P_1,Q_1,R_1)),\,(X_2,(P_2,Q_2,R_2))]=({\tilde X},({\tilde
P},{\tilde Q},{\tilde R})),
\end{equation}
where
$$X_1=\begin{pmatrix} a_1 & b_1\\ c_1 & -{}^t\!a_1\end{pmatrix},\quad
X_2=\begin{pmatrix} a_2 & b_2\\ c_2 & -{}^t\!a_2\end{pmatrix} \in \mathfrak
{g} $$ and
\begin{align*} \quad {\tilde X}&=X_1X_2-X_2X_1,\\ \quad {\tilde
P}&=P_1a_2+Q_1c_2-P_2a_1-Q_2c_1,\\ \quad {\tilde
Q}&=P_1b_2-Q_1\,^ta_2-P_2b_1+Q_2\,^ta_1,\\ {\tilde
R}&=P_1\,^tQ_2-Q_1\,^tP_2-P_2\,^tQ_1+Q_2\,^tP_1.
\end{align*}
Indeed, an element $(X,(P,Q,R))$ in $\mathfrak {g} ^J$ with
$X=\begin{pmatrix} a & b\\ c & -\,^ta\end{pmatrix} \in \mathfrak {g} $ may
be identified with the matrix
\begin{equation}
\begin{pmatrix} a & 0 & b & {}^tQ \\ P &
0 & Q & R \\ c & 0 & -{}^ta & -{}^tP \\ 0 & 0 & 0 & 0
\end{pmatrix},\quad b={}^tb,\ c={}^tc,\ R={}^tR
\end{equation}
in $\frak{sp}(n+m,\mathbb R).$
\smallskip

 Let us identify $\mathfrak {g}
_{n+m}:=\frak{sp}(n+m,\mathbb R)$ with its dual $\mathfrak {g}
_{n+m}^{\ast}$\,(see Proposition 6.1.3. (6.5)). In fact, there
exists a $G$-equivariant linear isomorphism
$$\mathfrak {g} _{n+m}^{\ast}\longrightarrow
 \mathfrak {g} _{n+m},\quad \lambda\mapsto
X_{\lambda}$$ characterized by
\begin{equation}
\lambda(Y)={\text{tr}}(X_{\lambda}Y),\quad Y\in\mathfrak {g} _{n+m}.
\end{equation}
Then the dual $(\mathfrak {g} ^J)^{\ast}$ of $\mathfrak {g}^J$ consists
of matrices of the form
\begin{equation}
\begin{pmatrix} x & p & y & 0 \\ 0 & 0 & 0 & 0 \\
z & q & -\,^tx & 0 \\ ^tq & r & -{}^tp & 0 \end{pmatrix},\quad
y={}^ty,\ z={}^tz,\ r={}^tr.
\end{equation}
There is a family of coadjoint orbits $\Omega_{\delta}$ which have the minimal
dimension $2n$, depending on a nonsingular $m\times m$ real
symmetric matrix parameter $\delta$ and are defined by the
equation
\begin{equation}
\delta=r,\quad
XJ_n=\begin{pmatrix} p\\ q\end{pmatrix} \delta^{-1}\,{}^t\!\begin{pmatrix} p \\
q\end{pmatrix},
\end{equation}
 where $X=\begin{pmatrix} x & y\\ z &
-\,^t\!x\end{pmatrix}$ with $y=\,^ty$ and $z=\,^tz$ in (9.55). Let
us denote by $\mathfrak {h}_{n,m}$ the Lie algebra of the
Heisenberg group $H_{\mathbb R}^{(n,m)}.$ Then the family
$\Omega_{\delta}\,(\,\delta={}^t\delta,\ \delta\in GL(m,\mathbb R))$
have the following properties $(\Omega 1)$-$(\Omega 2)$:
\smallskip

$(\Omega 1)$ Under the natural projection on $\mathfrak
{h}_{n,m}^{\ast}$, the orbit $\Omega_{\delta}$ goes to the orbit
which corresponds to the irreducible unitary representation
$U(\delta)$ of the Heisenberg group $H_{\mathbb R}^{(n,m)},$
namely, the Schr{\"o}dinger representation of
$H_{\mathbb{R}}^{(n,m)}$\,(cf.\,(8.19)).

\smallskip

$(\Omega 2)$ Under the projection on $\mathfrak
{g}^{\ast}=\frak{sp}(n,\mathbb{R})^{\ast},$ the orbit $\Omega_{\delta}$ goes to
$\Omega_{\text{sign}(\text{det}(\delta))}.$

\smallskip

In fact, there is an irreducible unitary representation
$\pi_{\delta}\,(\,\delta=\,^t\delta,\ \delta\in GL(m,\mathbb R))$ of
$G^{J}$ (or its universal cover) with properties
\begin{equation}
\text{Res}^{G^J}_{H_{\mathbb R}^{(n,m)}}\pi_{\delta}\cong
U(\delta),\quad \text{Res}^{G^J}_G \pi_{\delta}\cong
\pi_{\text{sign}(\text{det}(\delta))},
\end{equation}
 where $\pi_{\pm}$
are some representations of $G$ (or its universal cover)
corresponding to the minimal orbits $\Omega_{\pm}\subset \mathfrak
{g}^{\ast}.$ Indeed, $\pi_{\pm}$ are two irreducible components of
the Weil representation of $G$ and $\pi_{\delta}$ is one of the
irreducible components of the Weil representation of
$G^J$\,(cf.\,(9.47)). These are special cases of the so-called
{\it unipotent} representations of $G^J$. We refer to
\cite{Vo3}-\cite{Vo4},\,\cite{Vo6} for a more detail on unipotent
representations of a reductive Lie group.
\smallskip

Now we consider the case $m=n=1.$ If
$$g^{-1}=\begin{pmatrix} a & 0 & b & a\mu-b\lambda
\\ \lambda & 1 & \mu & \kappa
\\ c & 0 & d & c\mu-d\lambda \\ 0 & 0 & 0 & 1 \end{pmatrix}$$
is an element of the Jacobi group $G^J$, then its inverse is given by
$$g=\begin{pmatrix} d & 0 & -b & -\mu \\ c\mu -d\lambda  & 1 &  \lambda b-\mu a &
-\kappa \\ -c & 0 & a & \lambda \\ 0 & 0 & 0 & 1 \end{pmatrix}.$$
We put
$$X=\begin{pmatrix} 1 & 0 & 0 & 0\\ 0 & 0 & 0 & 0 \\ 0 & 0 & -1 & 0 \\ 0
& 0 & 0 & 0 \end{pmatrix},\quad Y=\begin{pmatrix} 0 & 0 & 1 & 0\\ 0 & 0 & 0
& 0 \\ 1 & 0 & 0 & 0 \\ 0 & 0 & 0 & 0 \end{pmatrix},\quad Z=\begin{pmatrix}
0 & 0 & 1 & 0\\ 0 & 0 & 0 & 0 \\ -1 & 0 & 0 & 0 \\ 0 & 0 & 0 & 0
\end{pmatrix},   $$ $$P=\begin{pmatrix} 0 & 1 & 0 & 0\\ 0 & 0 & 0 & 0 \\ 0
& 0 & 0 & 0 \\ 0 & 0 & -1 & 0 \end{pmatrix},\quad Q=\begin{pmatrix} 0 & 0 &
0 & 0\\ 0 & 0 & 0 & 0 \\ 0 & 1 & 0 & 0 \\ 1 & 0 & 0 & 0
\end{pmatrix},\quad  R=\begin{pmatrix} 0 & 0 & 0 & 0\\ 0 & 0 & 0 & 0 \\ 0 &
0 & 0 & 0 \\ 0 & 1 & 0 & 0 \end{pmatrix}.   $$ Then according to
(9.55), $X,Y,Z,P,Q,R$ form a basis for $(\mathfrak {g}^J)^{\ast}$.
By an easy computation, we see that the coadjoint orbits
$\Omega_X,\,\Omega_Y,\,\Omega_Z,\,\Omega_P,\,\Omega_Q,\, \Omega_R$
of $X,Y,Z,P,Q,R$ respectively are given by $$\Omega_X=\left\{
\begin{pmatrix} ad+bc & 0 & -2ab & 0\\ 0 & 0 & 0 & 0 \\ 2cd & 0 &
-(ad+bc) & 0 \\ 0 & 0 & 0 & 0 \end{pmatrix}\,\bigg|\ ad-bc=1,\
a,b,c,d\in \mathbb R\,\right\},$$ $$\Omega_Y=\left\{ \begin{pmatrix} bd-ac
& 0 & a^2-b^2 & 0\\ 0 & 0 & 0 & 0 \\ d^2-c^2 & 0 & ac-bd & 0 \\ 0
& 0 & 0 & 0 \end{pmatrix}\,\bigg|\ ad-bc=1,\ a,b,c,d\in
\mathbb R\,\right\},$$
$$\Omega_Z=\left\{ \begin{pmatrix} -(ac+bd) & 0 & a^2+b^2 & 0\\ 0 & 0 & 0
& 0 \\ -(c^2+d^2) & 0 & ac+bd & 0 \\ 0 & 0 & 0 & 0
\end{pmatrix}\,\bigg|\ ad-bc=1,\ a,b,c,d\in \mathbb R\,\right\},$$
$$\Omega_P=\left\{ \begin{pmatrix}  (2ad-1)\lambda-2ac\mu & a & 2ab\lambda -2a^2
\mu & 0\\ 0 & 0 & 0 & 0 \\ 2c^2\mu-2cd\lambda & c &
(1-2ad)\lambda+2ac\mu &
0 \\ c & 0 & -a & 0 \end{pmatrix}\,\bigg|\ \gathered ad-bc=1,\\
a,b,c,d,\lambda,\mu\in \mathbb R\endgathered \,\right\},$$
$$\Omega_Q=\left\{ \begin{pmatrix} (2ad-1)\mu-2bd\lambda & b & 2b^2\lambda-2ab\mu &
0\\ 0 & 0 & 0 & 0
\\ 2cd\mu-2d^2\lambda & d &  (1-2ad)\mu+2bd\lambda & 0 \\ d & 0 & -b & 0
\end{pmatrix}\,\bigg|\gathered ad-bc=1,\\ a,b,c,d,\lambda,\mu\in
\mathbb R\endgathered\,\right\}$$ and $$ \Omega_R= \left\{ \begin{pmatrix}
(a\mu-b\lambda)(c\mu-d\lambda) & a\mu-b\lambda &
-(a\mu-b\lambda)^2 & 0\\ 0 & 0 & 0 &
0 \\ (c\mu-d\lambda)^2 & c\mu-d\lambda & -(a\mu-b\lambda)(c\mu-d\lambda) & 0 \\
c\mu-d\lambda & 1 & b\lambda-a\mu & 0 \end{pmatrix}\,\bigg|\
\gathered ad-bc=1,\\ a,b,c,d,\lambda,\mu\in \mathbb R\endgathered
\,\right\}.$$ Moreover we put $$S= \begin{pmatrix} 0 & 0 & 1 & 0\\ 0 & 0
& 0 & 0 \\ 0 & 0 & 0 & 0 \\ 0 & 0 & 0 & 0 \end{pmatrix}\quad
\text{and}\quad T= \begin{pmatrix} 0 & 0 & 0 & 0\\ 0 & 0 & 0 & 0 \\ 1 & 0
& 0 & 0 \\ 0 & 0 & 0 & 0 \end{pmatrix}.$$ Then the coadjoint orbits
$\Omega_S$ and $\Omega_T$ of $S$ and $T$ are given by
$$\Omega_S=\left\{ \begin{pmatrix} -ab & 0 & a^2 & 0\\ 0 & 0 & 0 & 0 \\
-b^2 & 0 & ab & 0 \\ 0 & 0 & 0 & 0 \end{pmatrix}\,\bigg|\
a,b\in\mathbb R\,\right\}$$ and
$$\Omega_T=\left\{ \begin{pmatrix} ab & 0 & -a^2 & 0\\ 0 & 0 & 0 & 0 \\
b^2 & 0 & -ab & 0 \\ 0 & 0 & 0 & 0 \end{pmatrix}\,\bigg|\
a,b\in\mathbb R\,\right\}.$$ For an element of $(\mathfrak
{g}^J)^{\ast}$, we write
\begin{equation}
\begin{pmatrix} x & p & y+z & 0\\ 0 & 0 & 0 &
0
\\ y-z & q & -x & 0 \\ q & r & -p & 0
\end{pmatrix}=xX+yY+zZ+pP+qQ+rR.
\end{equation}
 The coadjoint orbit
$\Omega_X$ is represented by the one-sheeted hyperboloid
\begin{equation}
x^2+y^2-z^2=1 > 0,\quad p=q=r=0.
\end{equation}
 The coadjoint orbit
$\Omega_Y$ is also represented by the one-sheeted hyperboloid
(9.59). The coadjoint orbit $\Omega_Z$ is represented by the
two-sheeted hyperboloids
\begin{equation}
x^2+y^2=z^2-1 > 0,\quad p=q=r=0.
\end{equation}
The coadjoint $G^J$-orbit $\Omega_S$  of $S$ is
represented by the the cone
\begin{equation}
x^2+y^2=z^2 >0, \quad z>0, \quad p=q=r=0.
\end{equation}
 On the other hand, the coadjoint $G^J$-orbit
$\Omega_T$ of $T$ is represented by the cone
\begin{equation}
x^2+y^2=z^2>0, \quad z<0,\quad p=q=r=0.
\end{equation}
The coadjoint orbit $\Omega_P$ is represented by the variety
\begin{equation}
2pqx+(q^2-p^2)y+(p^2+q^2)z=0,\quad (p,q)\in\mathbb R^2-\{
(0,0)\},\quad r=0
\end{equation}
 in $\mathbb R^6$. The coadjoint orbit
orbit $\Omega_Q$ is represented by the variety (9.63) in
$\mathbb{R}^6.$ particular, we are interested in the coadjoint
orbits orbits $\Omega_{hR}\,(h\in\mathbb{R},\ h\not=0)$ of $hR$
which are represented by
\begin{equation}
x^2+y^2=z^2,\ x=h^{-1}pq,\ y+z=-h^{-1}p^2,\ y-z=h^{-1}q^2\ \text{and}\ r=h.
\end{equation}
For a fixed $h\not=0$, we note that
$\Omega_{hR}$ is two dimensional and satisfies the equation
(9.56). Indeed, from the above expression of $\Omega_{hR}$ and
(9.58), we have $X=\begin{pmatrix} x & y+z \\ y-z & -x \end{pmatrix}$ and
\begin{align*}
  x&= h(a\mu-b\lambda)(c\mu-d\lambda),\\ y+z&=-h(a\mu-b\lambda)^2,\\
y-z&=h(c\mu-d\lambda)^2,\\ p&=h(a\mu-b\lambda),\quad
q=h(c\mu-d\lambda),\quad r=h.
\end{align*}
 Hence these satisfy the
equation (9.56). An irreducible unitary representation $\pi_h$
that corresponds to a coadjoint orbit $\Omega_{hR}$ satisfies the
properties (9.57). In fact, $\pi_h$ is one of the irreducible
components of the so-called (Schr{\"o}dinger-)Weil representation
of $G^J$\,(cf.\,(9.47)). A coadjoint orbit $\Omega_{mR+\alpha X}$
or $\Omega_{mR+ \alpha Y}\,(m \in \mathbb{R}^{\times},\
\alpha\in\mathbb{R})$ is corresponded to a principal series
$\pi_{m,\alpha,{\frac 12}}$, the coadjoint orbit
$\Omega_{mR+kZ}\,(m\in \mathbb{R}^{\times},\ k\in \mathbb Z^+)$ of
$mR+kZ$ is attached to the discrete series $\pi_{m,k}^{\pm}$ of
$G^J$. There are no coadjoint $G^J$-orbits which correspond to the
 complimentary series $\pi_{m,\alpha,\nu}\,(m\in\mathbb
R^{\times},\,\alpha\in\mathbb R,\,\alpha^2<{\frac
12},\,\nu=\pm{\frac 12}).$ See \cite{B-S},\,pp.\,47-48. There are
no unitary representations of $G^J$ corresponding to the
$G^J$-orbits of $\alpha P_{\ast}+\beta Q_{\ast}$ with
$(\alpha,\beta)\neq (0,0).$

\medskip

Finally we mention that the coadjoint orbit $\Omega_{mR+\alpha X}$ or
$\Omega_{mR+\alpha Y}\,(m\in \mathbb R^{\times},\ \alpha\in\mathbb R)$ is
characterized by the variety
\begin{equation}
x^2+y^2-(z^2+\alpha^2)={2\over
m}pqx+{1\over m}(q^2-p^2)y+{1\over m}(p^2+q^2)z,\quad r=m.
\end{equation}
and the coadjoint orbit $\Omega_{mR+kZ} (m \in \mathbb R^{\times}, k \in \mathbb Z^{+})$
of $mR+kZ$ is represented by the
variety
\begin{equation}
x^2+y^2-(z^2-k^2)={2\over m}pqx+{1\over
m}(q^2-p^2)y+{1\over m}(p^2+q^2)z,\ z>0,\ r=m.
\end{equation}
 or
\begin{equation}
x^2+y^2-(z^2-k^2)={2\over m}pqx+{1\over m}(q^2-p^2)y+{1\over
m}(p^2+q^2)z,\ z<0,\ r=m
\end{equation}
 depending on the sign $\pm.$
\vspace{0.2in}\\

\footnotesize{

\end{document}